\documentclass[10pt]{article}
\usepackage[utf8]{inputenc}

\usepackage{amsmath,amsthm,amsfonts,amssymb,bm,mathdots}
\usepackage{xcolor}
\usepackage{geometry}
\usepackage{tikz,tikz-cd}
\usetikzlibrary{decorations.pathreplacing}
\usetikzlibrary{arrows}

\usepackage{ytableau}
\usepackage[noblocks]{authblk}
\usepackage{enumitem}
\usepackage{indentfirst}
\usepackage{url}
\usepackage{fdsymbol}
\usepackage{adjustbox}
\usepackage{stmaryrd}
\usepackage[hidelinks = true]{hyperref}
\usepackage{comment}

\theoremstyle{definition}
\newtheorem{thm}{Theorem}[section]
\newtheorem{prop}[thm]{Proposition}

\newtheorem{cor}[thm]{Corollary}
\newtheorem{lem}[thm]{Lemma}
\newtheorem{example}[thm]{Example}
\newtheorem{definition}[thm]{Definition}
\newtheorem{remark}[thm]{Remark}

\newcommand{\zigzagpath}[2][1]
{
\xdef\scale{#1}
\foreach \x [count=\i ] in {#2}
{
\ifnum\i=1
	\xdef\lastx{\x}
\else {
	\draw[thick] ((\i*\scale-2*\scale,\lastx*\scale)--(\i*\scale-\scale,\x*\scale);
	\xdef\lastx{\x}
}
\fi
}
}

\newcommand{\N}{\mathbb{N}}

\newcommand{\R}{\mathbb{R}}
\newcommand{\Z}{\mathbb{Z}}
\newcommand{\E}{\mathbb{E}}
\newcommand{\C}{\mathbb{C}}
\newcommand{\PP}{\mathbb{P}}
\newcommand{\bp}{\mathbf{p}}
\newcommand{\bk}{\mathbf{k}}
\newcommand{\bH}{\mathbf{H}}

\newcommand{\bb}{\boldsymbol{\upsilon}}
\newcommand{\bQ}{\mathbf Q}
\newcommand{\ttt}{\boldsymbol{\tau}}
\newcommand{\bL}{\mathbf L}

\newcommand{\bdel}{\boldsymbol{\delta}}

\begin{document}
\title{Airy limit for $\beta$-additions through Dunkl operators}
\author[1]{David Keating}
\affil[1]{Department of Mathematics, University of Illinois Urbana-Champaign, Champaign, USA\\dkeating@illinois.edu}
\author[2]{Jiaming Xu}
\affil[2]{Department of Mathematics, KTH Royal Institute of Technology, Stockholm, Sweden\\jxu0800@gmail.com}

\date{}

\maketitle

\begin{abstract}
 It is well known that the edge limit of Gaussian/Laguerre Beta-ensembles, as well as a large class of $\beta$-ensembles is given by the $\mathrm{Airy}(\beta)$ point process. We extend this universality result to a general class of additions of Gaussian and Laguerre ensembles, which were identified in \cite{AN} as projection of the ergodic measures of the $\beta$-corners process. In order to make sense of the $\beta$-addition, we introduce the Type-A Bessel function as the characteristic function of our matrix ensemble, following the approach of \cite{GM}, \cite{BCG}. Then we extract its moment information through the action of Dunkl operators, and obtain certain limiting functional expressed via conditional Brownian bridges for the Laplace transform of $\mathrm{Airy}(\beta)$. Our limit expression is universal up to proper rescaling among all of our additions, and agrees with the single-time Laplace transform expression from the concurrent work \cite{GXZ}. 
\end{abstract}

\setcounter{tocdepth}{2}
\tableofcontents

\section{Introduction}
\subsection{Overview}\label{sec:overview}
The focus of this article is the intersection of three topics in random matrix theory: addition of random matrices, $\beta$-ensembles, and edge universality. The studies of addition of independent large random matrices date back to \cite{Vo}, in which the author considers 
$$C(N)=A(N)+B(N),$$
where $A(N)$, $B(N)$ are two independent $N\times N$ real/complex self-adjoint matrices with deterministic eigenvalues $a_{1}(N)\ge \ldots\ge a_{N}(N)$, $b_{1}(N)\ge\ldots\ge b_{N}(N)$ and uniformly random eigenvectors. Voiculescu \cite{Vo} proves that if the empirical measures of $A(N)$, $B(N)$ satisfy

    $$\frac{1}{N}\sum_{i=1}^{N}\delta_{a_{i}(N)}\rightarrow \mu_{A},\quad  \frac{1}{N}\sum_{i=1}^{N}\delta_{b_{i}(N)}\rightarrow \mu_{B}$$
for some deterministic measures $\mu_{A}$, $\mu_{B}$ on $\R$, then 
$$\frac{1}{N}\sum_{i=1}^{N}\delta_{c_{i}(N)}\rightarrow \mu_{C},$$
where $$\mu_{C}=\mu_{A}\boxplus\mu_{B}$$ is also a deterministic measure on $\R$ known as the \emph{free convolution} of $\mu_{A}$ and $\mu_{B}$.
This classical result was an early breakthrough in the subject of free probability, which has various connections with different types of random matrix additions, see e.g. \cite{B1}, \cite{B2}, \cite{MSS}.

In random matrix literature, \emph{$\beta$-ensembles} are random $N$-tuple of real numbers $x_{1}\ge x_{2}\ge \ldots \ge x_{N}$ with density function
\begin{equation}\label{eq_betadensity}
    \frac{1}{Z_{N,\beta}}\prod_{1\le i<j\le N}(x_{i}-x_{j})^{\beta}\prod_{i=1}^{N}\exp\left[-V(x_{i})\right],
\end{equation}
where $Z_{N,\beta}$ is a normalizing constant and $V(x)$ is the potential function satisfying certain analytic conditions. When $\beta=1,2,4$, such density corresponds to invariant matrix ensembles with real/complex/quaternion entries, see e.g. \cite{D} for more details. For general $\beta>0$, the parameter is understood as the inverse temperature. Two of the most classical objects in this context are $\beta$-ensembles (G$\beta$E), when $V(x)=\frac{\beta}{2}x^{2}$, and Laguerre $\beta$-ensembles (L$\beta$E), when $V(x)=\left[\frac{\beta}{2}(a+1)-1\right]\log(x)-\frac{\beta}{2}x$. For G$\beta$E and L$\beta$E, Dumitriu and Edelman \cite{DE} give a concrete tridiagonal matrix realization of the density function. This construction is later extended to Jacobi ensembles by Killip and Nenciu \cite{KN}, then to a larger class of potentials $V(x)$ by Krishnapur, Rider, and Virag \cite{KRV}. The tridiagonal structure is a main tool in studying $\beta$-ensembles. 

To study the addition of invariant random matrices for general $\beta>0$, rather than using the density in Eqn.~\eqref{eq_betadensity}, an alternate algebraic approach was developed. This is defined in \cite{GM}, where the Type-A Bessel function, a multivariate special function with representation theoretic background, plays the role of the characteristic function for the matrix ensembles. On the level of Law of Large Numbers, \cite{CX} and \cite{Y} prove that for any fixed $\beta>0$, the asymptotic empirical measure of self-adjoint additions is still characterized by the free convolution as was known for $\beta=1,2$. In addition, \cite{GM}, \cite{Me}, \cite{MeP}, \cite{BCG} study the asymptotic empirical measures of the additions as $N\rightarrow\infty$ in the low and high temperature regimes (namely, when $\beta\rightarrow\infty$ and $\beta\rightarrow0$, respectively.), and obtain two other versions of convolutions that differ from the classical free convolution. Another similar story happens with the additions of $N\times M$ rectangular matrices with certain invariant properties: \cite{B1}, \cite{B2} study the asymptotic empirical measures of singular values for fixed $\beta=1,2$ when $N,M\rightarrow\infty$, \cite{Xu} generalize the additions to $\beta>0$ by using the Type-BC Bessel function, and \cite{GrM}, \cite{Gri}, \cite{Xu} obtained several limiting results about the rectangular empirical measures for low and high temperatures.

One common feature of the literature on matrix additions mentioned above is that they only consider global asymptotic results. In this text we study the edge behavior of the random spectrum. For a single matrix ensemble without any additions, much is known about the edge behavior. For example, it is known that the largest eigenvalue of a large class of random matrix ensembles, including G$\beta$E and $L\beta E$, converges weakly to the Tracy-Widom distribution ($\mathrm{TW}_{\beta}$) after proper rescaling. Moreover, one can consider the largest $k$ eigenvalues, where $k$ is any fixed positive integer, the same arguments show that the rescaled $k$-tuple converge to the largest $k$ particles of the $\mathrm{Airy}(\beta)$ process (abbreviated $\mathrm{Airy}(\beta)$). The universal edge limit of large random matrices was also extended to Wigner matrices by \cite{Sos}, and to generalized real/complex self-adjoint matrices in \cite{PS}, \cite{TV}, \cite{EYY}, \cite{LY}, \cite{BEY} and other works. For $\beta\ne 1,2,4$, when there are no invariant matrices, \cite{RRV}, \cite{GSh}, \cite{DG}, \cite{PaSh}, \cite{Sh}, \cite{BEY} and other works show the edge universality for a large class of $\beta$-ensembles defined by Eqn.(\ref{eq_betadensity}). We refer the readers to \cite[Section 1]{BEY} for a more detailed survey. 

One of the main contributions of this article is to extend the edge universality of the self-adjoint matrices to a class of additions consisting of finitely many G$\beta$E and L$\beta$E, for general $\beta>0$. For $\beta=1,2$, this problem was previously studied in \cite{Ahn} and \cite{JP} using different approaches. In these works the authors add invariant matrices of more general form, however, the results and proofs only hold for $\beta=1,2$. In the current work on $\beta$-additions, we use the same ensembles as in \cite{GM}, \cite{BCG}, whose definition rely on multivariate special functions. 

The proof of our main result relies on the moment method. Heuristically, when taking power sums
$$\sum_{i=1}^{N}\lambda_{i}(N)^{M}$$
of all eigenvalues, where $M\in \Z_{>0}$ is large enough compared to $N$, the statistic is dominated by the $\lambda_{i}(N)$s with largest absolute value. This allows one to study the edge behavior. A similar approach was initiated in \cite{Sos} (see also \cite{SS1}, \cite{SS2}), and later applied in \cite{Pa}, \cite{Sod}, \cite{GSh} to study the edge limit of various random matrix ensembles. One common feature of the above works is that the eigenvalue power sums are obtained by directly calculating powers of the corresponding Wigner matrices/tridiagonal realization of $\beta$-ensembles, and then taking the trace. Such an approach does not work in our setting, since to the best of our knowledge, a concrete matrix structure no longer exists for $\beta$-additions. Instead, we extract the moment information of the additions from their characteristic function using Type-A Dunkl operators, a class of differential-difference operators established in special function literature, see Section \ref{sec:pre} for details. The applications of Dunkl operators to random matrices was initiated by \cite{BCG}, using the same type of operators as in this text, and later \cite{Xu} extended such connection to rectangular matrices using Dunkl operators of type BC. Both these works study the expectation of  $\sum_{i=1}^{N}\lambda_{i}(N)^{k}$ for fixed $k\in \Z_{\ge 1}$, while in \cite{GXZ}, which is concurrent to this text, Type-A Dunkl operators are used to give explicit expressions of the Laplace transform of the $\mathrm{Airy}_{\beta}$ line ensemble, an infinite collection of continuous random curves, whose one-time marginal distribution gives $\mathrm{Airy}(\beta)$. 

\subsection{Matrix additions for general $\beta$}\label{sec:additiondef}
We start with the $\beta=2$ case which corresponds to additions of complex self-adjoint matrices. Fix $c,\alpha_0 \in \R_{\ge 0}$, and a sequence  $\{\alpha_i\}_{i\ge 1}$ with $\alpha_i \in \R$ and $\sum_{i=1}^\infty \alpha_i^2 <\infty$. We consider the addition of matrices of the form
\begin{align}\label{eq_olshanskiensemble}
c I + \sqrt{\alpha_0} \frac{X+X^*}{2} + \sum_{i=1}^\infty \alpha_i \left( \frac{1}{2} V_i V_i^* - I \right).
\end{align}
where $X$ is an $N \times N$ matrix with i.i.d. $\mathcal{N}(0,N)+i\mathcal{N}(0,N)$ entries and $V_i$ are i.i.d. columns of size $N$ with i.i.d. $\mathcal{N}(0,1)+i\mathcal{N}(0,1)$ entries. The matrix $\frac{X+X^*}{2}$ is the Gaussian Unitary ensemble. When $\alpha_{1}=\cdots=\alpha_{L}=a$, with all other $\alpha_i=0$, $i>L$, 
$$\sum_{i=1}^L   \frac{\alpha}{2} V_i V_i^*=\frac{\alpha}{2}YY^{*}$$ 
is the $N\times L$ Laguerre Unitary ensemble (with weight $\alpha$), where $Y$ is a $N\times L$ rectangular matrix with i.i.d. complex Gaussian entries. We note that Eqn.~(\ref{eq_olshanskiensemble}) is studied in asymptotic representation theory: it is the $N\times N$ corner of a random infinite complex Hermitian matrix, whose distribution is invariant under conjugation of $U(\infty)$ and ergodic (see \cite{OV}).

In this text we ignore the terms proportional to the identity by taking $c=\sum_{i=1}^\infty \alpha_{i}$, and study the $\beta$-version of above addition. Here $\beta>0$ is a parameter that is commonly viewed as the inverse temperature, and GUE, LUE are replaced by Gaussian $\beta$-ensemble (G$\beta$E) and Laguerre $\beta$-ensemble, (L$\beta$E) respectively. Since for $\beta\ne 1,2,4$ there is no concrete invariant matrix of $G\beta E$ and $L\beta E$, the $\beta$-addition is defined in terms of multiplication of the Type-A Bessel generating functions of the summands, as in \cite{GM}, \cite{BCG}. We discuss this in more detail in Section \ref{sec:pre}.

 We also highlight a connection of our ensembles and asymptotic representation theory. It is shown in \cite{OV} for $\beta=2$, and later proved in \cite{AN}, \cite{BR} for all $\beta>0$ that the class of $\beta$-additions  in Eqn.~(\ref{eq_olshanskiensemble}) can be identified with the (level $N$ projection of) the ergodic measures of the $\beta$-corners process. In particular, each of them gives a legitimate probability measure $P_{N}$ of an $N$-tuple $\lambda_1\ge\ldots\ge \lambda_N$ in $\R^{N}$. By classical free probability theory, the empirical measure of $\lambda_{i}$s converges to a deterministic measure $\mu$ on $\R$ as $N\rightarrow\infty$. The main focus of this text is the study of the upper edge behavior near $\mu_{+}$, the right endpoint of $\operatorname{supp}(\mu)$. This corresponds to the fluctuation of the largest eigenvalues $\lambda_{1}\ge \lambda_{2}\ge \ldots\ge \lambda_{k}$, for a finite fixed $k$.
 
\subsection{Main result}\label{sec:mainresult}
As a random point process on the real line parameterized by $\beta>0$, the $\mathrm{Airy}(\beta)$ process is known as the universal object characterizing the scaling limit of the largest eigenvalues of various random matrix ensembles. While for $\beta=1,2,4$ there have been several characterizations of the limiting point configuration and the corresponding Tracy-Widom distributions in the past literature, for general $\beta$ the explicit expressions for $\mathrm{Airy}(\beta)$ are less understood. In \cite{RRV}, the authors make sense of the stochastic Airy operator
\begin{equation}
SAO_{\beta}=-\frac{d^{2}}{dx^{2}}+x+\frac{2}{\sqrt{\beta}}B'(x)
\end{equation}
on $L^{2}(\R_{\ge 0})$ with Dirichlet boundary condition at zero, where $B(x)$ denotes the Brownian motion and the differential makes sense as Schwartz distribution, see \cite[Section 2]{RRV} for the precise analytic issues in defining such an operator. $SAO_{\beta}$ is known to have an orthonormal basis of eigenfunctions with eigenvalues $-\eta_{i}\ (i=1,2,\ldots)$, and the particles in the $\mathrm{Airy}(\beta)$ process are given by these eigenvalues.

On the other hand, \cite{GSh} introduces integral operators $\mathcal{U}(T)\ (T\ge 0)$, called the \emph{stochastic Airy semigroup}, such that by \cite[Corollary 2.12]{GSh},
\begin{equation}
    \mathcal{U}(T)\stackrel{d}{=}\exp\left(-\frac{T}{2}SAO_{\beta}\right)
\end{equation}
for each $T$. The right-hand side above is the unique stochastic operator with the same eigenfunctions as $SAO_{\beta}$ and eigenvalues $\exp\left(\frac{T\eta_{i}}{2}\right)$. In addition, $\mathcal{U}(T)$ is almost surely trace class (see \cite[Proposition 2.4]{GSh}), and the same reference gives an alternate characterization of $\mathrm{Airy}(\beta)$ in terms of its Laplace transform. In particular, we have the following expression.
\begin{prop}\cite[Proposition 2.14]{GSh}\label{prop:GS}
\begin{equation}\label{eq_GSfirstmoment}
    \begin{split}
    \E[\operatorname{Trace}(\mathcal{U}(T))]=&\E\left[\sum_{i=1}^{\infty}e^{\frac{T\eta_{i}}{2}}\right]\\
    =&\frac{2}{\pi}T^{-\frac{3}{2}}\E\left[\exp\left(-\frac{T^{\frac{3}{2}}}{2}\int_{0}^{1}\mathfrak{E}(t)dt+\frac{T^{\frac{3}{2}}}{2\beta}\int_{0}^{\infty}(l_{y})^{2}dy\right)\right],
    \end{split}
\end{equation}
where $\mathfrak{E}$ is a standard Brownian excursion on $[0,1]$, and each $l_{y}$ is the total local time of $\mathfrak{E}$ at level y.    
\end{prop}

Note that Proposition \ref{prop:GS} gives an explicit formula for the first moment of 
\begin{equation}\label{eq_laplace}
\sum_{i=1}^{\infty}e^{\frac{T\eta_{i}}{2}},
\end{equation}
the Laplace transform of $\mathrm{Airy}(\beta)$. Moreover, \cite{GSh} expresses the trace of $\mathcal{U}(T)$ via an explicit random kernel built from Brownian bridges conditioned to stay nonnegative. It might also be possible to extract explicit expressions for higher moments of Eqn.~(\ref{eq_laplace}) from the kernel although we do not pursue it here.

\vspace{0.2cm}
Our main result verifies that the upper edge limit of the following class of $\beta$-additions is given by $\mathrm{Airy}(\beta)$. We introduce the following notations necessary for stating the theorem.

Fix $\beta>0$. Let $X(N)$ be a $N\times N$ Gaussian $\beta$-ensemble, $$V_{N,L_{1}}(N),\ldots,V_{N,L_{k}}(N)$$ be independent $N\times L_{i}\ (L_{i}\ge N)$ Laguerre $\beta$-ensembles, and fix $\alpha_{1},\ldots,\alpha_{k}\in \R$. Assume that there exists $\gamma_{i}>0$ such that, as $N\rightarrow\infty$,
\begin{equation}\label{eq_gamma}
    \begin{split}
        \frac{L_{i}}{N}\rightarrow \gamma_{i}\ge 1
    \end{split}
\end{equation}
for each $i=1,2,\ldots,k$.

Fix $\alpha_0\ge 0$ and consider the $\beta$-addition 
\begin{equation}\label{eq_betaaddition}
    \sqrt{\alpha_0 } X(N)\boxplus_{N}^{\beta}\alpha_{1}V_{N,L_{1}}(N)\boxplus_{N}^{\beta}\ldots\boxplus_{N}^{\beta}\alpha_{k}V_{N,L_{k}}(N).
\end{equation}
For $\beta=1,2,4$ this is the usual additions of independent Gaussian and Wishart matrices, and for general $\beta>0$ we give a generalized definition in Section 2. Without loss of generality assume 
\begin{equation}\label{eq_alpha}
    |\alpha_{1}|>|\alpha_{2}|>\ldots>|\alpha_{k}|>0.
\end{equation}

For $l\in \Z_{\ge 1}$, let 
\begin{equation}\label{eq_k2}
    \kappa_{2}=\alpha_0 +\alpha^{2}_{1}\gamma_{1}+\ldots+\alpha^{2}_{k}\gamma_{k},\ \kappa_{2}(N)=\alpha_0 +\alpha^{2}_{1}\frac{L_{1}(N)}{N}+\ldots+\alpha^{2}_{k}\frac{L_{k}(N)}{N}, 
\end{equation}
and
\begin{equation}\label{eq_kl}
    \kappa_{l}=\alpha^{l}_{1}\gamma_{1}+\ldots+\alpha^{l}_{k}\gamma_{k},\ \kappa_{l}(N)=\alpha^{l}_{1}\frac{L_{1}(N)}{N}+\ldots+\alpha^{l}_{k}\frac{L_{k}(N)}{N}, 
\end{equation}
for $l\ne 2$.  We will only consider the case in which, for each $N=1,2,\ldots$, \begin{equation}\label{eq_positivekappa}
    \kappa_{l}(N)\ge 0
\end{equation} for all $l$. Let 
\begin{equation}\label{eq_vz}
    V(z)=\frac{1}{z}+\sum_{l=1}^{\infty}\kappa_{l}z^{l-1}
\end{equation}
be a meromorphic function on $\C$ locally near 0.
    \vspace{0.1cm}

\begin{thm}[Edge universality of $\beta$-ensemble additions]\label{thm:main}

Let $\lambda_{i}(N)\ (i=1,2,\ldots,N)$ denote the eigenvalues of our $N\times N$ addition in Eqn.~(\ref{eq_betaaddition}). Assume that (\ref{eq_gamma}), (\ref{eq_alpha}) and (\ref{eq_positivekappa}) hold,  then there exists a positive constant $\tilde{C}$ and a sequence of positive real numbers $\mu_{+}(N)$, such that
\begin{equation}\label{eq_mainresult}
    \lim_{N\rightarrow\infty}\left\{\lambda_{i}':=\frac{\lambda_{i}(N)-\mu_{+}(N)N}{N^{\frac{1}{3}}}
    \right\}_{i=1}^{N}\stackrel{d}{=} \tilde{C}\cdot\mathrm{Airy}(\beta),
\end{equation}
the $\mathrm{Airy}(\beta)$ process. 
\end{thm}

The following statement gives explicit descriptions of the parameters in Theorem \ref{thm:main}.
\begin{prop}\label{prop:main}
    Let $z_{c}(N)>0,\ z_{c}>0$ be the unique root of $V'_{N}(z),\ V'(z)$ on $(0,\frac{1}{a_{1}})$ respectively. Then $V_{N}(z_{c}(N))>0,\ V(z_{c})>0$, and  
    \begin{equation}\label{eq_c0}
    \begin{split}
        &\mu_{+}(N)=V_{N}(z_{c}(N))\rightarrow\mu_{+}=V(z_{c})\ \text{as}\ N\rightarrow\infty,\\
        &\tilde{C}=2^{-\frac{1}{3}}V''(z_{c})^{\frac{1}{3}}.
    \end{split}
    \end{equation}
\end{prop}

\begin{remark}
    A key restriction of our ensembles is that $\kappa_{l}\ge 0$ for all $l$. This is guaranteed when we have only additions of Laguerre ensembles, namely, all $\alpha_{i}>0$, $i\ge1$.  In general, by taking $a_{1}>0$, It is not hard to see that $\kappa_{l}>0$ for all $l$ large enough. However, what is not covered in our result is when the terms being subtracted contribute enough to make $\kappa_{l}<0$ for some $l$. We conjecture that our result are still true in that case, see Section \ref{sec:further} for more details.
\end{remark}

\begin{remark}
    Note that in the addition (\ref{eq_betaaddition}) we allow the absence of the Gaussian ensemble in the sum. It is possible that if $\alpha_0>0$ is large enough, the result in Theorem \ref{thm:main} can be (partially) obtained by the approach in \cite{AH}, which relies on the analysis of Dyson Brownian motion. However, their approach does not cover the cases $\beta\in (0,1)$, or $\alpha_0=0$.
\end{remark}

\begin{example}\label{ex:laguerre}
A special case of Theorem \ref{thm:main} gives the classical edge limit of a single Laguerre $\beta$-ensemble of variance 1. In this case we have 
\begin{equation}
    \begin{split}
        V(z)=\frac{1}{z}+\frac{\gamma}{1-z},\quad z_{c}=\frac{1}{\sqrt{\gamma}+1},\quad
        \mu_{+}=V(z_{c})=(\sqrt{\gamma}+1)^{2},
    \end{split}
\end{equation}
and \begin{equation}
    \tilde{C}=\frac{(\sqrt{\gamma}+1)^{\frac{4}{3}}}{\gamma^{\frac{1}{6}}}.
\end{equation}
Plugging in these parameters gives exactly \cite[Theorem 1.4]{RRV}.
\end{example}

\subsection{Proof ingredients}\label{sec:proofingredient}
Define $\lambda_{i}'=N^{-\frac{1}{3}}(\lambda_{i}(N)-\mu_{+}(N)N)$ as in Theorem \ref{thm:main}. The main target in this text is to calculate the asymptotic distribution of the $\lambda_{i}'$s in terms of the Laplace transform 
$$\sum_{i=1}^{N}\exp\left(\frac{T\lambda_{i}^{'}}{\mu_{+}}\right),$$
which is a collection of random variables parameterized by $T>0$.
More precisely, we compute the (mixed) $l^{th}$ moments of the Laplace transform \[\E\left[\prod_{j=1}^{l}\left(\sum_{i=1}^{N}\exp\left(\frac{\bk_{j}\lambda_{i}'}{\mu_{+}}\right)\right)\right]\] for all $l\in \Z_{\ge 1}$ by taking the asymptotics of $$\E\left[\prod_{j=1}^{l}\left(\sum_{i=1}^{N}\left(\frac{\lambda_{i}(N)}{\mu_{+}(N)N}\right)^{M_{j}}\right)\right],$$
where $M_{j}\approx \bk_{j}N^{\frac{2}{3}}$. In particular, we show that for each element in the class of additions we consider, the limiting moment expressions are all the same, after appropriate rescaling of $\lambda_{i}'$s. This extends the edge universality, since we know that among these $\beta$-additions the single Laguerre ensemble has an $\mathrm{Airy}(\beta)$ edge limit. 

The way we make sense of $\beta-$additions for general $\beta>0$ is based on the \emph{Type-A multivariate Bessel function} $$B_{(a_{1},\ldots,a_{N})}(x_{1},\ldots,x_{N};\beta),$$ which is a symmetric function in the variables $x_i$ and plays the role of characteristic function of our ensembles. We consider the actions of \emph{Type-A Dunkl operators}, a class of differential-difference operators $\mathcal{D}_{i}(i=1,2,\ldots,N)$, on $\E[B_{(\lambda_{1},\ldots,\lambda_{N})}(x_{1},\ldots,x_{N};\beta)]$. Since the power sums of $\mathcal{D}_{i}$ have Type-A Bessel functions as eigenfunctions, such actions extract the moment information of the $\lambda_{i}$s. 

The approach of using Dunkl actions to study random matrices follows the same line as \cite{BCG} and \cite{Xu}. But compared to these works, which take the global limit of the empirical measures, we focus on edge behavior and therefore we take large powers of $\lambda_{i}(N)$s of order $O(N^{\frac{2}{3}})$. This makes the asymptotic analysis of the explicit actions much more involved. A similar approach to this paper is taken in the concurrent work \cite{GXZ}, which takes Dunkl actions on $\beta$-Dyson Brownian motions and G$\beta$E corners processes, and establishes the first Laplace transform expression of Airy$_{\beta}$ line ensemble. Below we summarize the topics of the remaining sections, and simultaneously outline the common ingredients as well as the difference of these two works. 

In Section \ref{sec:pre}, we give the necessary backgrounds of Dunkl operators and Bessel generating functions (abbreviated as BGF). Thanks to the nice behavior of Gaussian and Laguerre ensembles, the BGF of our $\beta$-additions have an explicit form given in Eqn.~(\ref{eq_additionbgf}). The main object in \cite{GXZ}, however, is a single (multi-time/multi-level version of) G$\beta$E. In that case the BGF further simplifies to $$\exp\left(\frac{\alpha_{0}N}{\beta}\sum_{i=1}^{N}x_{i}^{2}\right).$$
While it is well known that under proper scaling, the global limit shape of G$\beta$E is the semicircle law of radius 2, the global limit shape of Eqn.~(\ref{eq_betaaddition}), and in particular $\mu_{+}$ in Theorem \ref{thm:main}, is obtained via a free convolution of one copy of G$\beta$E and $k$ copies of L$\beta$Es. We present the detailed connection to free probability in Section \ref{sec:freeprob}.

Section \ref{sec:dunklactions} studies the Dunkl action on our BGF for finite $N$ and arbitrary $\bk_{1},\ldots,\bk_{l}>0$. Similar to \cite{GXZ}, we introduce a probabilistic interpretation of the action as a weighted sum over certain (conditional) random walk bridges with i.i.d. increments. While the two works share the same class of operators, the expressions in \cite{GXZ} are in terms of Bernoulli walk bridges with increments $\pm1$ due to the simpler form of its BGF, while our walk bridges are known as \emph{\L{}ukasiewicz paths} in combinatorics and the increments take values in $\Z_{\ge -1}$. To more conveniently apply the probabilistic tools on these walk bridges, we specify the weights of the increments in such a way that they have mean 0. It turns out that the partition function of the increments then coincides with the scaling factor $\mu_{+}(N)$, which is essential in our characterization of the edge limit. In contrast to Bernoulli walks, there are much fewer explicit formulas and quantitative bounds for random \L{}ukasiewicz paths. We state a collection of such results in Section \ref{sec:technical} and Appendix A. To the best of our knowledge, these results have not been written down in the literature, and we believe that they are of independent interest. While our result covers all the finite additions of G$\beta$E and L$\beta$Es, part of the linear combinations involving minus signs seem out of reach from our probabilistic approach. We discuss such difficulties in more details in Section \ref{sec:further}.

After the above preparations, we present a self-contained framework that computes the asymptotics of the Dunkl actions as $N\rightarrow\infty$ in Section \ref{sec:asymptotics}. In particular, under proper scaling the random walk bridges converge weakly\footnote{Instead, \cite{GXZ} gives a quantitative bound between the pre-limit and limit expressions by coupling Bernoulli walks and Brownian bridges, see \cite[Lemma 4.40]{GXZ}.} to certain Brownian bridges. This observation, together with other significant technical efforts, lead to the convergence in moments in Theorem \ref{thm:mainconvergence}. Our framework has a number of ingredients in common with \cite[Section 4]{GXZ}, and therefore the limit expression in Theorem \ref{thm:mainconvergence} agrees with the single-time Laplace transform $\bL_\beta(\left(2P_{-1}\sigma\right)^{-2/3}\vec\bk, \vec\ttt=0)$ defined in \cite[Section 2]{GXZ}, up to a factor $\left(2P_{-1}\sigma\right)^{-2/3}$ that depends on the parameters in Eqn.~(\ref{eq_betaaddition}). While sharing some common techniques, the two works have different focuses in that one deals with multi-time/multi-level limit, and another identifies the limit expressions from a class of ensembles. Moreover, by combining the ideas from both works, one can prove that the edge limit of a class of dynamics/corners processes of $\beta$-additions is Airy$_{\beta}$ line ensemble.

Finally, in Section \ref{sec:momenttoweak} we pass from the moment convergence  to weak convergence of the rescaled eigenvalues. In particular, we show that the Laplace transform uniquely determines $\mathrm{Airy}(\beta)$.

This text (as well as \cite{GXZ}) gives a new explicit expression of the $\mathrm{Airy}(\beta)$ process. It is then natural to identify $\bL_\beta(\vec\bk, \vec\ttt=0)$ where $l=1$ with Eqn.~\eqref{eq_GSfirstmoment}, given in \cite{GSh}, or even identify $\bL_\beta(\vec\bk, \vec\ttt=0)$ for general $l\in \Z_{\ge 1}$ with some corresponding formulas of $\mathrm{Airy}(\beta)$ given in the literature, by separate combinatorial arguments. We leave this as an interesting open problem.

\subsection*{Acknowledgements} We thank Vadim Gorin and Lingfu Zhang for many stimulating discussions. We are also grateful to Yun Li for pointing out useful references about Airy processes, and to David Sivakoff, Andy Heeszel, Zihao Fang, Julien Sohier and Francesco Caravenna for pointing out useful references about random walks and conditional functional CLT. DK is supported by the NSF RTG grant DMS 1937241. 
\vspace{0.5cm}

\noindent\textbf{Notations.}
For any real numbers $a<b$, we denote $\llbracket a, b\rrbracket = [a, b]\cap\Z$. We also use $A\sim B$ to express $c_1\le A/B\le c_2$ for some constants $c_1<c_2$, uniform in the variables of $A$, $B$ according to the context. Similarly, $A\lesssim B$ expresses $A\le c_{3}B$ for some uniform constant $c_3>0$. For a real-valued function $f$ whose domain is a subset of $\R$, we use $f(x-)/f(x+)$ to denote the left/right limits of $f$ at $x\in \R$, when it exists. $f$ is \emph{C\'adl\'ag} if it is right continuous and has a left limit for each $x$ in the interior of its domain.

\section{Bessel generating functions}\label{sec:pre}

In this section we give a brief review of Type-A Dunkl operators and Bessel generating functions. See \cite{Ro}, \cite{A} and references therein for more details. In Section \ref{sec:bgf}, we use these generating function to define the notion of `addition' for general $\beta$. Throughout the whole text we take $\beta\in \R_{>0}$.  

\subsection{Dunkl operators}

Consider functions of $N$ variables, $x_1,\ldots,x_N$. Define a differential-difference operator acting on these functions, known as the \emph{Dunkl operator} (of Type-A), by
\begin{equation}\label{eq_dunkl}
\mathcal{D}_i := \frac{\partial}{\partial x_i} + \frac{\beta}{2} \sum_{j:j\ne i} \frac{1-s_{i,j}}{x_i-x_j}, \qquad i=1,\ldots,N
\end{equation}
where $s_{i,j}$ permutes the variables $x_i$ and $x_j$. A key property is that the Dunkl operators commute (see \cite{Du}), that is,
\begin{equation}
\mathcal{D}_i\mathcal{D}_j = \mathcal{D}_j\mathcal{D}_i, \qquad i,j=1,\ldots,N.
\end{equation}
We also define
\begin{equation}
\mathcal{P}_k:= \mathcal{D}_1^k + \cdots + \mathcal{D}_N^k, \qquad k>0,
\end{equation}
known as the \emph{symmetric Dunkl operators}.

\subsection{Multivariate Bessel functions}

Let $(a_1,a_2,\ldots,a_N)$ be an $N$-tuple of real numbers satisfying $a_1\ge a_2\ge \ldots \ge a_N$. The \emph{multivariate Bessel function}\footnote{Our notation slightly differs to the literature, where the parameter $\beta$ is replaced by $\theta:=\beta/2$.} $B_{(a_1,\ldots,a_N)}(x_1,\ldots,x_N;\beta)$ has many definitions. For example, it can be seen as the expectation of a certain observable of the $\beta$-corner process or as a limit of Jack symmetric polynomials, see e.g. \cite{OO}, \cite[Section 2.2]{BCG}.

Here we list several well-known properties of these functions:
\begin{itemize}
    \item They are symmetric in $x_1,\ldots,x_N$. This follows from their definition as a limit of Jack polynomials.
    \item They can be extended to an entire function of $2N$ variables $a_1,\ldots,a_N,x_1,\ldots,x_N$.
    \item They have the normalization $B_{(a_1,\ldots,a_N)}(0,\ldots,0;\beta)=1$.
    \item They are eigenfunctions of the $\mathcal{P}_k$ operators. In particular,
    \[
    \mathcal{P}_k B_{(a_1,\ldots,a_N)}(x_1,\ldots,x_N;\beta) = \left(\sum_{i=1}^N a_i^k \right) B_{(a_1,\ldots,a_N)}(x_1,\ldots,x_N;\beta).
    \]
\end{itemize}

\subsection{Bessel generating functions}\label{sec:bgf}
Let $W_N:=\{(a_1,a_2,\ldots,a_N)\in \R^N | a_1\ge a_2 \ge \ldots \ge a_N\}$. Given a probability measure $\mu$ on $W_N$, we define the \emph{Bessel generating function for $\mu$} by
\begin{equation}\label{eq:BGFint}
G_N(x_1,\ldots,x_N;\beta,\mu) = \int_{a_1\le \ldots \le a_N} \;B_{(a_1,\ldots,a_N)}(x_1,\ldots,x_N;\beta) \; d\mu(a_1,\ldots,a_N).
\end{equation}
Later in this text, we will drop $\mu$ from the notation of BGF when the measure is clear from the context. It follows from the symmetry of the Bessel functions that $G_N(x_1,\ldots,x_N;\beta)$ is also symmetric in the variables $x_1,\ldots,x_N$, and from the normalization of the Bessel functions we see that $G_N(0,\ldots,0;\beta)=1$.

We will be interested in measures $\mu$ that satisfy a certain decay condition. 
\begin{definition}\cite[Definition 2.7]{BCG} 
    We say that a measure $\mu$ on $W_N$ is \emph{exponentially decaying} with exponent $R>0$ if
    \[
    \int_{a_1\ge \ldots \ge a_N} e^{NR\max_i |a_i|} d\mu(a_1,\ldots,a_N) < \infty. 
    \]
\end{definition}
This guarantees that the integral \eqref{eq:BGFint} converges for all $(x_1,\ldots, x_N)$ in the domain
\[
\Omega_R:= \{(x_1,\ldots, x_N) \in \C^N: |\Re x_i|<R,i=1,\ldots,N \},
\]
see \cite[Lemma 2.9]{BCG}. We will only consider such measures on $W_N$.

When $\beta=1,2,4$, the Bessel function is known as the spherical transform for certain symmetric space (see \cite{Hel}). For example, when $\beta=2$, we have
\begin{equation}\label{eq:matrixint}
B_{(a_1,\ldots,a_N)}(x_1,\ldots,x_N;\beta) =\int_{U} \;\exp{\text{Tr}(UAU^*X)}\; dU
\end{equation}
where $A=\text{diag}(a_1,\ldots,a_N)$, $X= \text{diag}(x_1,\ldots,x_N)$ and $dU$ is the Haar measure on $N\times N$ Unitary matrices. One can view 
\[
M_A = UAU^*
\]
as a random hermitian matrix with fixed eigenvalues given by $a_1\ge\ldots \ge a_N$.

Suppose we take the sum of two such random matrices, $M_A$ and $M_B$, with fixed eigenvalues $a_1\ge\ldots \ge a_N$ and $b_1\ge \ldots \ge b_N$, respectively. Consider the eigenvalues $c_1,\ldots,c_N$ of $M_{A}+M_{B}$. By Eqn.~\eqref{eq:matrixint} and the linearity of the trace we have
\[
\mathbb{E}[B_{(c_1,\ldots,c_N)} (x_1,\ldots,x_N;\beta)| \vec a,\vec b] = B_{(a_1,\ldots,a_N)}(x_1,\ldots,x_N;\beta)B_{(b_1,\ldots,b_N)}(x_1,\ldots,x_N;\beta)
\]
It follows that 
\[
G_N (x_1,\ldots,x_N;\beta,\mu_C) = G_N(x_1,\ldots,x_N;\beta,\mu_A)G_N(x_1,\ldots,x_N;\beta,\mu_B).
\]

When $\beta\ne 1,2,4$ there is no underlying invariant matrices to sum, so we define the operation of addition in terms of the generating functions.

\begin{definition}\label{def:sum}
Let $\mu_A$ and $\mu_B$ be exponentially decaying measures on $W_N$. Define $\mu_C:= \mu_A \boxplus_{N}^{\beta} \mu_B$ (also written as $C=A\boxplus_{N}^{\beta}B$) to be the measure whose Bessel generating function is given by
\begin{equation}\label{eq:BGFsum}
G_N (x_1,\ldots,x_N;\beta,\mu_C) = G_N(x_1,\ldots,x_N;\beta,\mu_A)G_N(x_1,\ldots,x_N;\beta,\mu_B).
\end{equation}
\end{definition}

We note that it is a well-known open problem in the general case to construct a (unique) probability measure on $N$-tuple $(c_1,\ldots, c_N)$ for which Eqn.~\eqref{eq:BGFsum} holds. However, for the sums we consider this has been done in \cite{AN}. We denote this measure as $P_{N,\beta}$.

We will focus on two particular matrix ensembles: the Gaussian $\beta$-ensembles and the Laguerre $\beta$-ensembles.

\begin{definition}\label{GaussianandLaguerre}
    For $\beta>0$, an $N\times N$ \emph{Gaussian $\beta$-ensemble} is a probability measure on $N$-tuple $(\lambda_{1}\ge \ldots\ge \lambda_{N})\in \R^{N}$ with density function 
    \begin{equation}\label{eq_Gaussiandensity}    
   \frac{1}{Z_{N,\beta}}\prod_{1\le i<j\le N}(\lambda_{i}-\lambda_{j})^{\beta}\prod_{i=1}^{N}\exp\left(-\frac{\lambda_{i}^{2}}{2N}\right),
    \end{equation}
    and an $N\times L$ \emph{Laguerre $\beta$-ensemble} ($N\le L)$ is a probability measure on $N$-tuple $(\lambda_{1}\ge \ldots\ge \lambda_{N})\in \R^{N}$ with density function
    \begin{equation}\label{eq_Laguerre}
       \frac{1}{Z_{N,L,\beta}}\prod_{i=1}^{N}\left[\lambda_{i}^{\frac{\beta}{2}(L-N+1)-1}\exp\left(-\frac{\beta\lambda_{i}}{2} \right)\right]\prod_{1\le i\le j\le N}(\lambda_{i}-\lambda_{j})^{\beta}.
    \end{equation}
    $Z_{N,\beta}$ and $Z_{N,L,\beta}$ are certain normalizing constants.
\end{definition}

\begin{prop}\label{prop:GaussianLaguerrebgf}
For the $N\times N$ G$\beta$E, the Bessel generating function is given by
\[
G_N(x_1,\ldots,x_N;\beta)= \exp\left(\frac{N(x_1^2+\cdots + x_N^2)}{\beta}\right).
\]

For the $N\times L$ L$\beta$E, the Bessel generating function is given by
\[
G_N(x_1,\ldots,x_N;\beta)=\prod_{i=1}^{N}\left(1-\frac{2x_{i}}{\beta}\right)^{-L\beta/2}.
\]
\end{prop}
Sketches of proofs of these two facts are given in \cite{BCG} and \cite{Ne}, respectively.

Fix $\beta>0$, $\alpha_0\ge 0$, $\alpha_{1},\ldots,\alpha_{k}\in \R$. Consider the $\beta$-addition 
\begin{equation}\label{eq:addition}
\sqrt{\alpha_0 } X(N)\boxplus_{N}^{\beta}\alpha_{1}V_{N,L_{1}}(N)\boxplus_{N}^{\beta}\ldots\boxplus_{N}^{\beta}\alpha_{k}V_{N,L_{k}}(N).
\end{equation}

If $\lambda_1 \ge \ldots \ge \lambda_N$ are the random eigenvalues of this sum, then by Definition \ref{def:sum} and Proposition \ref{prop:GaussianLaguerrebgf}, their Bessel generating function is given by
\begin{align}
G_{N}(x_1,\ldots,x_N;\beta) = \E[B_{(\lambda_1,\ldots,\lambda_N)}(x_1,\ldots,x_N; \beta)] = \prod_{k=1}^N \mathcal{F}(x_k)
\end{align}
where
\begin{equation}\label{eq_Fz}
\begin{split}
\mathcal{F}(x) = & \exp \left(\frac{\alpha_0 N}{\beta}x^2 \right) \prod_{i=1}^k \left(1-\frac{2\alpha_i}{\beta} x \right)^{-L_{i}\beta/2}  \\
= & \exp\left(\frac{\alpha_0 N}{\beta}x^2 + \frac{\beta}{2} \sum_{l=1}^\infty \frac{2^l}{l}  \frac{\alpha_{1}^{l}L_{1}+\ldots+\alpha_{k}^{l}L_{k}}{\beta^l}x^l  \right).
\end{split}
\end{equation}
We see that
\begin{equation}\label{eq_additionbgf}
    G_{N}(x_{1},\ldots,x_{N};\beta)= \exp\left( \frac{\alpha_0 N}{\beta}p_2(x_1,\ldots,x_N) + \frac{\beta}{2} \sum_{l=1}^\infty  \frac{2^l}{l}  \frac{\alpha_{1}^{l}L_{1}+\ldots+\alpha_{k}^{l}L_{k}}{\beta^l}p_l(x_1,\ldots,x_N)  \right),
\end{equation}
where $p_k(x_{1},\ldots,x_{N})$ is the $k$-th power sum symmetric function.

\begin{remark}
    This class of matrices was chosen in large part because the multivariate Bessel generating function factors into the product of single variable functions $\mathcal{F}(z)$. This greatly simplifies the computation in what follows. However, the measure $P_{N,\beta}$ supporting this BGF is also of independent interest, known as (up to shift by a constant) the extremal consistent distribution on the infinite families of interlacing arrays under the Dixon-Anderson conditional
probability distribution, see \cite{AN}. The result in that paper generalizes the classification of ergodic unitary invariant infinite complex random matrices in \cite{OV} to all $\beta>0$.
\end{remark}

One of the key features of multivariate Bessel functions is that they are the (unique up to normalization) eigenfunctions of the Dunkl operators, see Proposition \ref{prop:eigenfunction} below. A probabilistic version of this result is given in Proposition \ref{prop:momentformula}. This will be the starting point of our moment calculations in the remaining sections.

 Let $p_{k}(x_{1},\ldots,x_{N})$ be the degree $k$ symmetric power sum of $x_{1},\ldots,x_{N}$, where $x_{i} (i=1,2,\ldots,N)$ are formal variables. 
\begin{prop}\cite{Op}\label{prop:eigenfunction}
Fix $\beta>0$, $N\in \mathbb{Z}_{\ge 1}$. For any $k=1,2,\ldots$, 
\begin{equation}
    p_{k}(\mathcal{D}_{1},\ldots,\mathcal{D}_{N})B_{(a_{1},\ldots,a_{N};\beta)}(x_{1},\ldots,x_{N})=p_{k}(a_{1},\ldots,a_{N})B_{(a_{1},\ldots,a_{N})}(x_{1},\ldots,x_{N};\beta).
\end{equation}
Here we evaluate the polynomial $p_k$ on the operators $\mathcal{D}_1, \dots, \mathcal{D}_N$.
\end{prop}

\begin{prop}\label{prop:momentformula}
 Fix $\beta>0$, $N\in \mathbb{Z}_{\ge 1}$. Define $G_{N}(x_{1},\ldots,x_{N};\beta)$ as in Eqn.~(\ref{eq_additionbgf}), and let $\lambda_{1}\ge \ldots\ge \lambda_{N}$ be its random eigenvalues distributed according to $P_{N,\beta}$. Then for any $k_{1},\ldots,k_{m}\in \Z_{\ge 1}$, we have
    \begin{equation}
        \E\left[\prod_{i=1}^{m}p_{k_{i}}(\lambda_{1},\ldots,\lambda_{N})\right]=\left(\prod_{i=1}^{m}p_{k_{i}}(\mathcal{D}_{1},\ldots,\mathcal{D}_{N})\right)G_{N}(x_{1},\ldots,x_{N};\beta)|_{x_1=\ldots=x_N=0}.
    \end{equation}
\end{prop}
\begin{proof}
This follows immediately from Proposition \ref{prop:eigenfunction} and the fact $B_{(a_{1},\ldots,a_{N})}(0,\ldots,0;\beta)=1$,  (see \cite[Proposition 2.10]{BCG} for more details), once we check that $P_{N,\beta}$ is exponentially decaying.

Consider the addition 
$$C=A\boxplus_{N}^{\beta}B,$$
where $A$, $B$ have deterministic eigenvalues $\vec{\lambda}(A)=(\lambda_{1}(A)\ge \ldots\ge \lambda_{N}(A))$, $\vec{\lambda}(B)=(\lambda_{1}(B)\ge\ldots\ge\lambda_{N}(B))$ and C has random eigenvalues $\vec{\lambda}(C)=(\lambda_{1}(C)\ge \ldots\ge \lambda_{N}(C))$. By \cite[Lemma 3.23]{A}, 
$$\max_{i}|\lambda_{i}(C)|\le \max_{i}|\lambda_{i}(A)|+\max_{i}|\lambda_{i}(B)|.$$
Moreover, it is clear form the density functions of G$\beta$E and L$\beta$E that for a given ensemble of the form Eqn.~(\ref{eq_betaaddition}), there exist constants $C>0$, $R>0$ such that $\PP[\max_{i}|\lambda_{i}|>x]<C\exp(-2NRx)$, where $\vec{\lambda}$ are the eigenvalues of a summand in Eqn.~(\ref{eq_betaaddition}) 
(note that here $N$ is a fixed parameter). A union bound then gives that Eqn.~(\ref{eq_betaaddition}) itself is exponentially decaying with exponent $R/(k+1)$.
\end{proof}

\subsection{Connection to free probability}\label{sec:freeprob}

We briefly recall some basic properties of free convolution that relate to the asymptotic behavior of self-adjoint matrix additions. Given a measure $\mu$ on $\R$, the free cumulants $\{\kappa_{l}(\mu)\}_{l=1}^{\infty}\in \R^{\infty}$ are a collection of quantities that linearize the free convolution, namely, for $\mu_{C}=\mu_{A}\boxplus\mu_{B}$
$$\kappa_{l}(\mu_{C})=\kappa_{l}(\mu_{A})+\kappa_{l}(\mu_{B})$$
for each $l=1,2,\ldots$ For a more detailed exposition, see e.g. \cite{No}, \cite{NS}.

Since in this text we study the additions of $\beta$-Hermitian matrices that generalize \cite{Vo}, it is not surprising that free cumulants again come into play. More precisely, they appear as the coefficients of the $\log G_{N}$ as we will now describe.

Recall that for each $N=1,2,\ldots$, we consider the addition
\begin{equation}
    \begin{split}
        \sqrt{\alpha_0} X(N)\boxplus_{N}^{\beta}\alpha_{1}V_{N,L_{1}}(N)\boxplus_{N}^{\beta}\ldots\boxplus_{N}^{\beta}\alpha_{k}V_{N,L_{k}}(N),
    \end{split}
\end{equation}
where $X(N)$ is a $N\times N$ Gaussian ensemble and $V_{N,L_{i}}(N)$ is a $N\times L_{i}(N)$ Laguerre ensemble. For simplicity, we will always assume $\alpha_{1},\ldots,\alpha_{k}$ 
are independent of $N$.

Define the finite-$N$ version of the free cumulants as follows. Let 
\begin{equation}\label{eq:Nfreecumulant}
\begin{split}
    &\kappa_{1}(N)=\frac{1}{N}\Big(\alpha_{1}L_{1}(N)+\ldots+\alpha_{k}L_{k}(N)\Big)\\
    &\kappa_{2}(N)=\frac{1}{N}\Big(\alpha_0 N+\alpha^{2}_{1}L_{1}(N)+\ldots+\alpha_{k}^{2}L_{k}(N)\Big)\\
    &\kappa_{l}(N)=\frac{1}{N}\Big(\alpha^{l}_{1}L_{1}(N)+\ldots+\alpha^{l}_{k}L_{k}(N)\Big)\qquad l\ge3.
\end{split}
\end{equation}

We assume that $ \gamma_{1}(N):=L_{1}(N)/N\to \gamma_{1}>0,\ldots,\gamma_{k}(N):=L_{k}(N)/N\to \gamma_{k}>0$ as $N\to \infty$, and without loss of generality assume that $\alpha_{1}>0$ (the case $\alpha_{1}<0$ only differs by a factor -1). In terms of the Bessel generating function, 
\begin{equation}
\begin{aligned}
G_{N}(x_1,\ldots,x_N;\beta) &\; = \exp\left( \frac{\alpha_0 N}{\beta}p_2(x_1,\ldots,x_N) + \frac{\beta}{2} \sum_{l=1}^\infty  \frac{2^l}{l}  \frac{\alpha_{1}^{l}L_{1}+\ldots+\alpha_{k}^{l}L_{k}}{\beta^l}p_l(x_1,\ldots,x_N)  \right)  \\
&\; = \exp \left(N\sum_{l=1}^{\infty}\frac{2^{l-1}}{l}\frac{\kappa_l(N)}{\beta^{l-1}} p_l(x_1,\ldots,x_N) \right).
\end{aligned}
\end{equation}
Therefore, \begin{align*}
    \frac{1}{N}\partial_i \ln G
   (x_1,\ldots,x_N;\beta)|_{x_1=\ldots=x_N=0} = & \kappa_{1}(N)\rightarrow \kappa_{1}, \\
    \frac{1}{N}\partial^2_i \ln G(x_1,\ldots,x_N;\beta)|_{x_1=\ldots=x_N=0} = & \frac{2\kappa_{2}(N)}{\beta}\rightarrow \frac{2\kappa_{2}}{\beta}, \\
   \frac{1}{N}\partial^l_i \ln G(x_1,\ldots,x_N;\beta)|_{x_1=\ldots=x_N=0} = & (l-1)! \frac{2^{l-1}\kappa_{l}(N)}{\beta^{l-1}}\rightarrow  (l-1)! \frac{2^{l-1}\kappa_{l}}{\beta^{l-1}}, \qquad l \ge 3,
\end{align*} where 
\begin{equation}\label{eq_cumulant}
    \begin{split}
&\kappa_{1}=\alpha_{1}\gamma_{1}+\ldots+\alpha_{k}\gamma_{k}\\
        &\kappa_{2}=\alpha_0+\alpha^{2}_{1}\gamma_{1}+\ldots+\alpha^{2}_{k}\gamma_{k}\\
        &\kappa_{l}=\alpha^{l}_{1}\gamma_{1}+\ldots+\alpha^{l}_{k}\gamma_{k}\qquad l\ge 3.\\
    \end{split}
\end{equation}
 Since $\alpha_{1}>0$, we see that $\kappa_{l}(N)>0$ for all $l$ large enough.

Under the setting above, by \cite[Theorem 1.1]{CX}, the empirical measure 
    $$\frac{1}{N}\sum_{i=1}^{N}\delta_{\lambda_{i}/N}$$
of our addition converges weakly in probability to a deterministic probability measure $\mu$ on $\R$ as $N\rightarrow\infty$. The $\mu_{+}$ in Theorem \ref{thm:main} and Proposition \ref{prop:main} is the right endpoint of $\operatorname{supp}(\mu)$. In addition, the $\kappa_{l}$ are precisely the free cumulants of $\mu$. They completely determine $\mu$ because of the one-to-one correspondence to moments of $\mu$ (see Eqn.~(\ref{eq_momentcumulant})), and the fact that $\mu$ is exponentially decaying.

\begin{example}[Gaussian ensemble]
    Letting $\alpha_i=0$ for $i\ge 1$ and $\alpha_0 = 1$  results in the G$\beta$E. Note that all the free cumulants are zero, except for $\kappa_2=1$. This implies that the limiting spectral measure is given by Wigner's semi-circle law, which has density
    \begin{align}
    d\mu(x) = \frac{1}{2\pi} \sqrt{4-x^2}, \qquad |x|\le 2
    \end{align}
    and $0$ otherwise. The non-zero moments are given by
    \begin{align}
        m_{2n} = \int_{-2}^{2} x^{2n} d\mu(x) = C_n 
    \end{align}
    where $C_n= \frac{1}{n+1} \binom{2n}{n}$ is the $n$-th Catalan number.
\end{example}

\begin{example}[Laguerre ensemble]
Suppose we take $\alpha_1=1$ and $\alpha_0,\alpha_{2},\alpha_{3},\ldots$ to be 0; then we get a $N\times L$ L$\beta$E. When $N,L\to \infty$ such that $L/N \rightarrow \gamma>1$, the limiting spectral measure is given by the Marchenko-Pastur distribution
\begin{align}
d\mu(x) = \frac{1}{2\pi} \frac{\sqrt{(x-\mu_-)(\mu_+-x)}}{ x}, \qquad \mu_-\le x \le  \mu_+
\end{align}
where $\mu_\pm = (\sqrt{\gamma}\pm 1)^2$. It is known that the Voiculescu $R$-transform characterizes the distribution, and for $\mu$ it is given by
\begin{align}
R(w) = \frac{\gamma}{1-w} = \sum_{n=0}^\infty \gamma w^n.
\end{align}
 The $R$-transform $R(w)$ is also the generating function of the free cumulants. In this case, we see that $\kappa_{l}=\gamma$ for each $l$.

\end{example}

\section{Dunkl actions and configurations}\label{sec:dunklactions}
\subsection{Walks and blocks }\label{sec:configuration}

Fix $l\in \Z_{\ge 1}$, $\bk_{1},\ldots,\bk_{l}\in \R_{>0}$, and let $M_{1}(N),\ldots,M_{l}(N)$ be sequences of positive integers (in the remainder of this text, we omit the dependence of $M_{j}$s on $N$ in the notation), such that for some constant $C>0$ uniform in $N$ and $j=1,2,\ldots,l$,
$$\left|M_{j}- \bk_{j}N^{2/3}\right|<C.$$
Let  $M=M_{1}+\ldots+ M_{l}$. In this section, we calculate the asymptotics of mixed moments of the form
\begin{equation}\label{eq_mixedmoment}
    \E\left[\prod_{j=1}^{l}\left(\sum_{i=1}^{N}\left(\frac{\lambda_{i}}{\mu_{+}(N)N}\right)^{M_{j}}\right)\right].
\end{equation}

Recall from Proposition \ref{prop:momentformula} that, for $l\in \Z_{\ge 1}$, the above can be obtained from the action 
\begin{equation}\label{eq_dunklpowersums}
\begin{split}
    &\frac{1}{\mu_{+}(N)^{M}}\cdot \prod_{j=1}^{l}\left[\left(\frac{\mathcal{D}_{1}}{N}\right)^{M_{j}}+\ldots+\left(\frac{\mathcal{D}_{N}}{N}\right)^{M_{j}}\right]G_{N}(\vec{x};\beta)\\
    &=\frac{1}{\mu_{+}(N)^{M}}\cdot\sum_{i_{l}=1}^{N}\cdots\sum_{i_{1}=1}^{N}\left(\frac{\mathcal{D}_{i_{l}}}{N}\right)^{M_{l}}\cdots \left(\frac{\mathcal{D}_{i_{1}}}{N}\right)^{M_{1}}G_{N}(\vec{x};\beta).
\end{split}
\end{equation}

Consider the action of a single term in the above expansion,
\begin{equation}\label{eq_keyaction}
    \frac{1}{\mu_{+}(N)^{M}}\left(\frac{\mathcal{D}_{i_{l}}}{N}\right)^{M_{l}}\cdots \left(\frac{\mathcal{D}_{i_{1}}}{N}\right)^{M_{1}}G_{N}(\vec{x};\beta)\Bigg|_{\vec{x}=0}.
\end{equation}

Recall that for the ensemble we consider,
\[
G_{N}(x_{1},\ldots,x_{N};\beta) = \exp\left(N\sum_{l=1}^{\infty} \frac{2^{l-1}}{l}\frac{\kappa_{l}(N)}{\beta^{l-1}} p_l(x_1,\ldots,x_N) \right).
\]
and, from Eqn.~(\ref{eq_dunkl}), that the Dunkl operator is given by 
\begin{equation}
    \frac{1}{N}\mathcal{D}_{i} = \frac{1}{N}\partial_{i}+\frac{\beta}{2N}\sum_{j=1}^{N}\frac{1-s_{ij}}{x_{i}-x_{j}}.
\end{equation}
 The Dunkl operator consists of two parts: a partial derivative term $\partial_{i}$ and the swapping term $\frac{1-s_{ij}}{x_{i}-x_{j}}$. We will consider the action of each term separately.

Let $c\in \R$, $n_{1},\ldots,n_{N}\in \Z_{\ge 0}$, and consider a monomial in $x_{1}$,\ldots,$x_{N}$ given by
$$c\cdot x_{1}^{n_{1}}\cdots x_{N}^{n_{N}}.$$ In order to understand the repeated action of the Dunkl operators on the Bessel generating function, we will compute the action of $N^{-1}\mathcal{D}_{i}$ on functions of the form 
\[
\left( c \prod_{i=1}^N x_i^{n_i} \right) G_{N}(\vec{x};\beta).
\]
By straightforward calculation, we have that the action of the partial derivative can be written as
\begin{equation}\label{eq_partial}
    \begin{split}
       \left( \sum_{l=1}^{\infty}\frac{2^{l-1}\kappa_{l}(N)}{\beta^{l-1}} x_i^{l-1}+\frac{n_{i}}{N}x_{i}^{-1}\mathbf{1}[n_i\ge 1]\right)\left( c \prod_{i=1}^N x_i^{n_i} \right) G_{N}(\vec{x};\beta) ,
    \end{split}
\end{equation}
while, for $j\ne i$, the action of the swapping term can be written
\begin{equation}\label{eq_swap}
\resizebox{0.9\textwidth}{!}{$
        \left(\frac{\beta}{2N}\cdot\frac{1-s_{ij}}{x_{i}-x_{j}}\right)\left[\left( c \prod_{i=1}^N x_i^{n_i} \right)G_N(\vec{x};\beta)\right]
        =\begin{cases}
            \frac{\beta}{2N}\left(\sum_{k=n_{j}}^{n_{i}-1}x_{i}^{k}x_{j}^{n_{i}+n_{j}-1-k}\right)\left(c\prod\limits_{k\ne i,j}x_{k}^{n_{k}}\right)G_{N}(\vec{x};\beta),\quad & n_{i}> n_{j}\\
            -\frac{\beta}{2N}\left(\sum_{k=n_{i}}^{n_{j}-1}x_{i}^{k}x_{j}^{n_{i}+n_{j}-1-k}\right)\left(c\prod\limits_{k\ne i,j}x_{k}^{n_{k}}\right) G_{N}(\vec{x};\beta),\quad & n_{i}< n_{j}\\
            0, & n_i=n_j.
        \end{cases}
    $}
\end{equation}
In either case, the Bessel generating function remains unchanged and only the form of the prefactor is altered.  We introduce the following linear operators on monomials of the form $\prod_{i=1}^N x_i^{n_i}$,  $n_1,\ldots,n_N\in \Z_{\ge 0}$.
\begin{itemize}
    \item For $1\le i\le N$, $a\in \Z_{\ge 1}$, $$U^{(a)}_{i}(x_{1}^{n_{1}}\cdots x_{N}^{n_{N}}):=\frac{2^{a-1}\kappa_a(N)}{\beta^{a-1}}x_{i}^{a-1}\cdot x_{1}^{n_{1}}\cdots x_{N}^{n_{N}}.$$
\item For $1\le i\le N$, $$d_{i}(x_{1}^{n_{1}}\cdots x_{N}^{n_{N}}):= n_{i}\mathbf{1}[n_i\ge 1]\cdot x_{1}^{n_{1}}\cdots x_{i-1}^{n_{i-1}}x_i^{n_i-1}x_{i+1}^{n_{i+1}}\cdots x_{N}^{n_{N}}.$$ 

\item For $1\le i\ne j\le N$,  $$d_{i,j}(x_{1}^{n_{1}}\cdots x_{N}^{n_{N}}):=\begin{cases}
        x_{1}^{n_{1}}\cdots x_{i-1}^{n_{i-1}}x_{i}^{n_{i}-1}x_{i+1}^{n_{i+1}}\cdots x_{N}^{n_{N}},& n_{i}>n_{j};\\
        0,& n_{i}\le n_{j}.
    \end{cases}$$
\item For $1\le i\ne j\le N$, 
$$\tilde{d}^{(a)}_{i,j}(x_{1}^{n_{1}}\cdots x_{N}^{n_{N}}):=\begin{cases}
    x_{1}^{n_{1}}\cdots x_{N}^{n_{N}}\cdot {x_i^{-a} x_j^{a-1}}, & \text{if } n_i>n_j,\ 1< a\le n_{i}-n_{j}\\
    -x_{1}^{n_{1}}\cdots x_{N}^{n_{N}}\cdot x_i^{a-1} x_j^{-a}, & \text{if } n_i < n_j,\ 1\le a\le n_{j}-n_{i}\\
    0, & \text{otherwise}.
\end{cases}$$
\end{itemize}
Here we have split the partial derivative term with respect to $x_i$ into
\[
\sum_{a=1}^\infty U^{(a)}_i + \frac{1}{N}d_i
\]
and swapping term into
\[
\frac{\beta}{2N} \sum_{j\ne i} \left(d_{i,j} + \sum_{a=1}^\infty \Tilde{d}_{i,j}^{(a)}\right).
\]
By Eqn.~(\ref{eq_partial}) and (\ref{eq_swap}),
\begin{equation}
\resizebox{0.9\textwidth}{!}{$
\begin{aligned}
    &\frac{1}{\mu_{+}(N)^{M}}\left(\frac{\mathcal{D}_{i_{l}}}{N}\right)^{M_{l}}\cdots \left(\frac{\mathcal{D}_{i_{1}}}{N}\right)^{M_{1}}G_{N}(\vec{x};\beta)
    \\
    &=\frac{1}{\mu_{+}(N)^{M}}\left(\sum_{a=1}^\infty U^{(a)}_{i_{l}}+\frac{1}{N}d_{i_{l}}+\frac{\beta}{2N}\sum_{j\ne i_{l}}\left(d_{i_l,j}+\sum_{a=1}^{\infty}\Tilde{d}^{(a)}_{i_{l},j}\right)\right)^{M_{l}}
\cdots \left(\sum_{a=1}^\infty U^{(a)}_{i_{1}}+\frac{1}{N}d_{i_{1}}+\frac{\beta}{2N}\sum_{j\ne i_{1}}\left(d_{i_1,j}+\sum_{a=1}^{\infty}\Tilde{d}^{(a)}_{i_{1},j}\right)\right)^{M_{1}} (1) \cdot G_{N}(\vec{x};\beta) \\
&=\frac{1}{\mu_{+}(N)^{M}}\left(\sum_{T(t):1\le t\le M }\prod_{t=1}^{M}T(t)\right)(1)\cdot G_{N}(\vec{x};\beta),
\end{aligned}
$}
\end{equation}
where the summation in the last line is over the finitely many choices of operators $T(t)$ for each $t\in \llbracket 1,M \rrbracket$. In particular, let 
\begin{align*}
C_0 = 0, \quad C_{p}=\sum_{j=1}^{p}M_{j}, \quad p=1,2,\ldots,l,
\end{align*}
then for $t \in \llbracket C_{p-1}+1, C_p \rrbracket$, $T(t)$ is chosen from the $U^{(a)}_{i_{p}}$, $\frac{1}{N}d_{i_{p}}$, $\frac{\beta}{2N} \tilde{d}^{(a)}_{i_{p},j}$ and $\frac{\beta}{2N} d_{i_{p},j}$. In this way, we may view each term in the sum as a discrete-time process of monomials in $x_1,\ldots,x_N$ starting at time zero from the monomial $P(0)=1$ with the monomial at time $t>0$ given by $P(t)=T(t)\cdots T(1)(1),\ t=1,2,\ldots,M$.

\begin{definition}
    For $N\in \Z_{\ge 1}$, let $P(\vec{x})=c\cdot x_{1}^{n_{1}}\cdots x_N^{n_{N}}$ be a monomial in $x_1,\ldots,x_N$. We denote 
    $$\operatorname{Deg}(P):=\sum_{i=1}^{N}n_{i},\quad \operatorname{Deg}_{i}(P):=n_{i}\quad \text{for}\ i=1,2,\ldots,N.$$
\end{definition}

 Observe that after the action of $T(t)\ (t\in \llbracket 1,M\rrbracket)$, if $P(t)\ne 0$, the degree has either
\begin{itemize}
\item increased by $a-1$ if $T(t) = U^{(a)}_{i}$ for some $i$, $a\in \Z_{\ge 1}$, or
\item decreased by 1 otherwise.
\end{itemize}
Moreover, only those sequences $\{P(t)\}_{t=0}^{M}$ for which $P(M)$ is a scalar will make a nontrivial contribution to Eqn.~(\ref{eq_keyaction}) as we always evaluate at $\vec x = \vec 0$. These observations motivate us to construct a correspondence between each process $\{P(t)\}_{t=0}^{M}$ that contributes to the action and a \L{}ukasiewicz path $E(t)$, which we also call a \emph{walk}.

\begin{definition}\label{def:walk}
     Fix $N,l\in \Z_{\ge 1}$ and $M_{1},\ldots,M_{l}\in \Z_{\ge 1}$, $M_1+\ldots+M_l=M$. For each sequence of monomials $\{P(s)\}_{s=0}^{M}$, coming form the action $\frac{1}{\mu_{+}(N)^{M}}\left(\frac{\mathcal{D}_{i_{l}}}{N}\right)^{M_{l}}\cdots \left(\frac{\mathcal{D}_{i_{1}}}{N}\right)^{M_{1}}G_{N}(x_1,\ldots,x_N;\beta)|_{\vec{x}=0}$, let $E(t):[0,M]\rightarrow \Z_{\ge 0}$ be a C\'adl\'ag function defined by 
        \begin{equation}
        E(t)=\operatorname{Deg}(P(\lfloor t\rfloor)).
        \end{equation} 
        We refer to $E(t)$ as a walk.
  \end{definition}

  While the walk records the total degree, we also need to keep track of the degree of each individual variable. Due to the actions of the $\tilde{d}_{i,j}$, each $\operatorname{Deg}_{i}(P(t))$ changes in certain nontrivial ways, which we record by a collection of \emph{blocks} $q_{i}(t)$. More precisely, we have the following. \begin{definition}\label{def:block}
        Fix $N,l\in \Z_{\ge 1}$ and $M_{1},\ldots,M_{l}\in \Z_{\ge 1}$, $M_1+\ldots+M_l=M$. For each sequence of monomials $\{P(s)\}_{s=0}^{M}$, coming form the action $\frac{1}{\mu_{+}(N)^{M}}\left(\frac{\mathcal{D}_{i_{l}}}{N}\right)^{M_{l}}\cdots \left(\frac{\mathcal{D}_{i_{1}}}{N}\right)^{M_{1}}G_{N}(x_1,\ldots,x_N;\beta)|_{\vec{x}=0}$, let $q_{i}(t):[0,M]\rightarrow \Z_{\ge 0}$ be a piecewise function, which is C\'adl\'ag on $[0,M]\setminus\{C_{0},\ldots,C_{l}\}$, defined by 
        \begin{equation}
            q_{i}(t)=\begin{cases}
                \operatorname{Deg}_{i}(P(\lfloor t\rfloor)),\quad & t\in [0,M]\setminus \bigcup\limits_{j:i_{j}=i}(C_{j-1},C_{j}]\\
                0,\quad &\text{otherwise}.
            \end{cases}
        \end{equation}
        Let $$q(t)=\sum_{i=1}^N q_{i}(t),\quad t\in [0,M].$$
         We call $\vec{q}(t)=(q_{1}(t),\ldots,q_{N}(t))$ the blocks of $\{P(s)\}_{s=0}^{M}$. 
  \end{definition}
  
  We call the pair $(E,\vec{q})$ a \emph{configuration} created by the action $\frac{1}{\mu_{+}(N)^{M}}\left(\frac{\mathcal{D}_{i_{l}}}{N}\right)^{M_{l}}\cdots \left(\frac{\mathcal{D}_{i_{1}}}{N}\right)^{M_{1}}G_{N}(\vec{x};\beta)|_{\vec{x}=0}$, and define $\mathcal{B}(i_{1},\ldots,i_{l})$ to be the set of all such configurations. Additionally, for $t\in [0,M]$, we call $x_{i}$ the \emph{main variable}, when $t\in \bigcup_{j:i_{j}=i}(C_{j-1},C_{j}]$. A key observation is that on $\bigcup_{j:i_{j}=i}(C_{j-1},C_{j}]$,
  \[
  E(t)-q(t)\in \Z_{\ge 0}=\operatorname{Deg}_i(P(\lfloor t\rfloor)).
  \]
  When $t\in [0,M]\setminus \bigcup_{j:i_{j}=i}(C_{j-1},C_{j}]$ (so $x_{i}$ is not the main variable), each $q_{i}(t)$ has finitely many discontinuities which we call \emph{jumps of $x_{i}$}. Each jump is corresponds to the action of one of the operators $\tilde{d}^{(a)}_{i,j}$, and each $q_{i}(t)$ is determined by the jumps 
$$\{(s_{m}(i),h_{m}(i))\}_{m=1,2,\ldots,\delta_{i},\ i=1,2,..N},$$ such that $(s_{m}(i),h_{m}(i))$ gives the \emph{occurrence time} $s_m(i)\in \llbracket 0,M\rrbracket\setminus \bigcup_{j:i_{j}=i}(C_{i_{j}-1},C_{i_{j}}]$, and the \emph{height} $$h_m(i)=q_{i}(s_{m}(i))-q_{i}(s_{m}(i)-)\in \Z\setminus \{0\}.$$  
Note that by Eqn.~(\ref{eq_swap}), for any jump of $x_{i}$ at time $t$, the jump height $h:=q_{i}(t)-q_{i}(t-)$ satisfies the following two rules:
\begin{equation}\label{eq_jumprules}
      \begin{split}
    &|h| \le |E(t-)-q(t-)-q_i(t-)-\mathbf{1}[h>0]|;\\
     &h\cdot (E(t-)-q(t-)-q_i(t-)) > 0.
    \end{split}
\end{equation}
We call $x_{i}$ an \emph{auxiliary variable} if $i\notin\{i_{1},\ldots,i_{l}\}$, and it has at least one jump (in other words, the degree of $x_i$ is nonzero at some point during the action).

By the above construction, each sequence of monomials $\{P( t)\}_{t=0}^M$ such that $P(M)|_{\vec x = 0} \ne 0$ corresponds to a unique configuration $(E,\vec q)$. We assign a weight $w(E,\vec{q})$ to each of these configurations as follows. Let $\delta=\sum_{i=1}^{N}\delta_{i}$ be the number of jumps in $(E,\vec{q})$, $\Delta\subset \llbracket 1,M\rrbracket$ be the set of jump times, $m_{i}(t)$ be the number of variables $x_{j}$ such that $\operatorname{Deg}_{j}(P(t))<\operatorname{Deg}_{i}(P(t))$, which gives the number of $d_{i,j}$ that act non-trivially on $P(t)$, and let $D(E)\subset \llbracket 1, M\rrbracket$ be the set of times\footnote{One can think of this as the set of the ``usual down-steps".} at which $E(t)-E(t-)=-1$, and $t\notin \Delta$. 

We then define the weight of a configuration according to
\begin{equation}\label{eq_weight}
    \begin{split}
        w(E,\vec{q})=\frac{1}{\mu_{+}(N)^{M}}&\prod_{l=1}^{\infty} \left(\frac{2^{l-1}\kappa_{l}(N)}{\beta^{l-1}}\right)^{\#\ t:\ E(t)-E(t-)=l-1}\\
        \times&  \frac{\beta^\delta \prod_{t\in \Delta}(-1)^{\mathbf{1}[q(t)<q(t-)]}}{2^\delta N^{\delta}}\prod_{t:\  t+1\in D(E)}\left(\frac{\beta m_{i}(t)}{2N}+\frac{E(t)-q(t)}{N}\right)\\
        =\frac{1}{\mu_{+}(N)^{M}}&\prod_{l=1}^{\infty} \kappa_{l}(N)^{\#\ t:\ E(t)-E(t-)=l-1}\\
        \times&  \frac{ \prod_{t\in \Delta}(-1)^{\mathbf{1}[q(t)<q(t-)]}}{N^{\delta}}\prod_{t:\  t+1\in D(E)}\left(\frac{ m_{i}(t)}{N}+2\frac{E(t)-q(t)}{\beta N}\right)
    \end{split}
\end{equation}
The second equality above holds, since the action of the three terms except for $\frac{1}{N}d_{i}$ in \[\sum_{a=1}^\infty U^{(a)}_i + \frac{1}{N}d_i+\frac{\beta}{2N} \sum_{j\ne i} \left(d_{i,j} + \sum_{a=1}^\infty \Tilde{d}_{i,j}\right)\]
satisfy that, after each step, the power of $(\beta/2)^{-1}$ introduced equals the change in height of the walk, and the total height change of $E$ after time $M$ is 0.

By the above discussions, we have
\begin{align*}
     \frac{1}{\mu_{+}(N)^{M}}\left(\frac{\mathcal{D}_{i_{l}}}{N}\right)^{M_{l}}&\cdots \left(\frac{\mathcal{D}_{i_{1}}}{N}\right)^{M_{1}}G_{N}(\vec{x};\beta)\Bigg|_{\vec{x}=0} =\sum_{(E,\vec{q})\in \mathcal{B}(i_{1},\ldots,i_{l})}w(E,\vec{q}).
\end{align*}

\subsection{A probabilistic point of view}\label{sec:conditionalrw}

Note that the first line of the right-hand side of the second equality in Eqn.~(\ref{eq_weight}) can be rewritten as \begin{equation}\label{eq_product}\prod_{t\in \llbracket 1,M\rrbracket}\frac{\kappa_{E(t)-E(t-)+1}(N)}{\mu_{+}(N)},\end{equation}
where we take $\kappa_{0}(N)$ to be 1. In this section, we introduce a probabilistic point of view for Eqn.~(\ref{eq_product}) and use it to specify $\mu_{+}(N)$.

 Given a sequence  $\{\kappa_{l}(N)\}_{l=1}^{\infty}$ defined by Eqn.~(\ref{eq:Nfreecumulant}), we treat them as the free cumulants of a probability measure, and let
\[
V_{N}(z) = \frac{1}{z} + \sum_{d= 1}^\infty \kappa_l(N) z^{l-1} 
\]
be their Voiculescu transform. Starting from this section, we assume that $$\kappa_{l}(N)\ge 0\ \ \text{for\ all}\ \ l\in \Z_{\ge 1}.$$ Then it is direct to check that 
\begin{itemize}
    \item $V_{N}'(z)$ has a unique positive zero $z_c(N)$ in its radius of convergence, and 
    \item  $V''_{N}(z_c(N))>0$.
\end{itemize}  
We then define a random variable $X(N)$ taking values in $\{-1,0,1,2,\ldots\}$ by
\begin{equation}\label{eq_rvx}
\begin{aligned}
    &P_{-1}(N):=\mathbb{P}(X(N)=-1) = \frac{1}{V_{N}(z_{c}(N))z_{c}(N)}, \\ &P_{l-1}(N):=\mathbb{P}(X(N)=l-1) = \frac{\kappa_{l}(N)z_{c}(N)^{l-1}}{V_{N}(z_{c}(N))}, \qquad l\ge 1.
\end{aligned}
\end{equation}

Denote the distribution of $X(N)$ by $\mu(N)$. A simple calculation shows that
\begin{equation}
\begin{aligned}
    \mathbb{E}[X(N)] = \frac{z_{c}(N)V_{N}'(z_{c}(N))}{V_N(z_{c}(N))}=0, \qquad \mathbb{E}[X(N)^2] = \frac{z_{c}(N)^2V_{N}''(z_{c}(N))}{V_{N}(z_{c}(N))}:=\sigma_{N}^{2}. 
\end{aligned}
\end{equation}

Moreover, by Eqn.~(\ref{eq_cumulant}),
\begin{equation}\label{eq_sigma}
\begin{split}
\mathbb{P}(X(N)=-1)\to P_{-1}:&=\frac{1}{V(z_{c})z_{c}},\quad\mathbb{P}(X(N)=l-1)\rightarrow P_{l-1}:=\frac{\kappa_{l}z_{c}^{l-1}}{V(z_{c})} ,\\
 \sigma_{N}^{2}\rightarrow \sigma^{2}:&=\frac{z_c^2V''(z_c)}{V(z_c)}\quad \text{as}\ N\rightarrow\infty. \end{split}\end{equation} 

We set \begin{equation}\mu_{+}(N):=V_{N}(z_{c}(N))\rightarrow \mu_{+}=V(z_{c})\quad \text{as}\; N\rightarrow\infty.\end{equation}

\subsection{Conditional walk bridges}\label{sec:walkbridges}
Using the interpretation in Section \ref{sec:conditionalrw}, we will consider $E(t)$ as a (conditional) walk bridge with i.i.d. increments of distribution $\mu(N)$. When $l=1$, $E(t)$ is simply a walk excursion of length $M=M_{1}$. However, for general $l\in \Z_{\ge 1}$, it will be more natural to view $E(t)$ as a concatenation of walk bridges with nonnegative heights and general boundary points. As a preparation, we introduce the following notions.

Let $L\in \Z_{\ge 1}, H_{1},H_{2}\in \Z$, $n\in \Z$. We say that $W(t):[n,n+L]\rightarrow \Z$ is a walk bridge of length $L$ from $H_{1}$ to $H_{2}$, if\begin{itemize}
    \item $W(t)$ is a C\'adl\'ag function, such that $W(t)=W(\lfloor t\rfloor)$,
    \item $W(n)=H_{1}$, $W(n+L)=H_{2}$.
\end{itemize}
We identify two walk bridges $W_{1}:[n_{1},n_{1}+L]\rightarrow\Z$ and $W_{2}:[n_{2},n_{2}+L]\rightarrow\Z$ if $W_{2}(t)=W_{1}(t+n_{2}-n_{1})$.

\begin{definition}\label{def:walkbridge}
    For $H_{1},H_{2}\in \Z_{\ge 0}$, $L\in \Z_{\ge 1}$, $N\in \Z_{\ge 1}$, $\mathcal{\tilde{W}}(H_{1},H_{2},L;N)$ is the set of walk bridges of length $L$  from $H_{1}$ to $H_{2}$. 
    $\mathcal{W}(H_{1},H_{2},L;N)$ is the set of walk bridges $W(t)$ of length $L$ from $H_{1}$ to $H_{2}$, such that $W(t)\ge0$ for all $t$ in its domain.
    
    For each element in $\mathcal{\tilde{W}}(H_{1},H_{2},L;N)$ or $\mathcal{W}(H_{1},H_{2},L;N)$, by viewing the increments of $W(t)$ as i.i.d. random variables with distribution $\mu(N)$, defined in Eqn.~(\ref{eq_rvx}), we assign a weight to each walk bridge as $$\text{weight}(W)=\prod_{l=-1}^{\infty}[P_{l}(N)V_{N}(z_{c}(N))]^{\# t:\ W(t)-W(t-)=l}.$$

    Let $$|\mathcal{\tilde{W}}(H_{1},H_{2},L;N)|=\sum_{W\in \mathcal{\tilde{W}}(H_{1},H_{2},L;N)}\text{weight}(W)$$
     be the \emph{partition function} of $\mathcal{\tilde{W}}(H_{1},H_{2},L;N)$. The partition function of $\mathcal{W}(H_{1},H_{2},L;N)$ is defined similarly.
\end{definition}

In Section \ref{sec:technical} and \ref{sec:asymptotics}, we will consider conditional walk bridges in $\mathcal{W}(H_{1},H_{2},L;N)$, and the asymptotics of their partition functions of all the four possible types: 1. $H_{1}=H_{2}=0$, 2. $H_{1}>0,\ H_{2}=0$, 3. $H_{1}=0,\ H_{2}>0$ and 4. $H_{1},H_{2}>0$. 
Each type above describes a segment of $E(t)$ in our configuration $(E,\vec{q})$. 

\begin{remark}
    When $P_{-1}(N)=P_{1}(N)=\frac{1}{2}$, $|\mathcal{\tilde{W}}(H_{1},H_{2},L;N)|$ can be counted by binomial theorem, and $|\mathcal{W}(H_{1},H_{2},L;N)|$ is given by the reflection principle. This corresponds to the moments of G$\beta$E, which is studied in great details in \cite{GXZ}. On the other hand, in our context the increments take value in $\Z_{\ge -1}$, which creates extra technicalities in the  computation of the partition functions. 
\end{remark}

Unlike for Bernoulli walks, there is an explicit formula of $|\mathcal{W}(H_{1},H_{2},L;N)|$ only when $H_{2}=0$.
\begin{lem}\label{lem:excursionpartition}
   For $H\in \Z_{\ge 0}$, $L\in \Z_{\ge 1}$, the partition function for the walk bridges with length $L$ from $H$ to $0$ is 
\begin{equation}\label{eq_contourwalkbridge}
  |\mathcal{W}(H,0,L;N)|=\frac{H+1}{L+1}|\tilde{W}(H+1,L+1;N)|z_{c}(N)=\frac{H+1}{L+1} \frac{1}{2\pi i} \oint\; \left(\frac{z}{z_c(N)}\right)^H V_{N}(z)^{L+1} \; dz
\end{equation}
where our contour only encloses the pole at zero.
\end{lem}
We prove Lemma \ref{lem:excursionpartition} in  Appendix A, by solving a ballot problem that generalizes the reflection principle for Bernoulli walks. Our approach, however, seems to no longer work when $H_{1}=0, H_{2}>0$. 

\begin{remark}
    When taking $H=0$, the walk bridge reduces to a walk excursion of length $L$, and we recover a classical result in free probability. It is well known (see e.g. \cite[Section 1.6]{No}) that for a probability measure $\mu$ with finite moments, its $L^{th}$ moment $m_{L}$ is given in terms of its free cumulants $\kappa_{l},\ l\in \Z_{\ge 1}$ by   \begin{equation}\label{eq_momentcumulant}
m_L = \sum_{\pi\in NC(L)} \prod_{B\in \pi} \kappa_{|B|},
\end{equation}
where the sum is over non-crossing partitions $\pi=B_{1}\sqcup\ldots\sqcup B_{n}$ of the set $\{1,2,\ldots,L\}$, $B\in \pi$ is a block of the partition, and $|B|$ is the size of the block. Another classical result (see e.g. \cite[Proposition 9.8]{NS}) gives a combinatorial bijection between a non-crossing partition of size $L$ and a \L{}ukasiewicz path of length $L$. One can also find the contour integral of $m_{L}$ in the literature, see e.g. \cite[Appendix]{BCG}.
\end{remark}

\section{Properties of random walk bridges}\label{sec:technical}

In this section, we provide several results regarding the walks sampled from $\mathcal{W}(H_{1},H_{2},L;N)$ or $\mathcal{\tilde{W}}(H_{1},H_{2},L;N)$.

\subsection{Conditional functional CLT's}\label{sec:conditionalfunctionalCLT}
Here we state a general result on the convergence of random walks to Brownian excursions due to \cite{CC} (see also \cite{S}). Let $W(n)=W(0)+X_{1}(N)+\ldots+X_{n}(N), n=0,1,2,\ldots$ be an integer-valued aperiodic random walk with i.i.d. increments $X(N)$ whose distribution depends on $N$, and $\E[X(N)]=0$. We assume that $W$ lies in the domain of attraction of the standard Gaussian law, namely, there exists $\sigma>0$ such that
\begin{equation}
    \frac{\sum_{i=1}^{L}X_{i}(N)}{\sigma\sqrt{L}}\implies \mathbf{N}(0,1)\quad \text{as}\ L\rightarrow\infty, \ N\rightarrow\infty.
\end{equation}

For $t>0$, $L,N\in \Z_{\ge1}$,
let $P_{t,L,N}^{*,x,y}$ be the law of $W(t)$ conditioned to be non-negative and satisfy $W(0)=x\ge 0, W(\lfloor tL\rfloor)=y\ge 0$. Define $\Omega$ to be the space of c\'adl\'ag functions on $\R_{\ge 0}$ endowed with the standard Skorohod topology and define $\psi_{t,L}: \Z^{\lfloor tL\rfloor} \to \Omega$ by
    \begin{equation}
    \psi_{t,L} (W(1),\ldots,W(\lfloor tL\rfloor))(s) = \left(\frac{1}{\sigma \sqrt{L}}W(\lfloor sL\rfloor)\right)_{0\le s\le t}.
    \end{equation}

\begin{thm}\cite[Theorem. 2.4 and Corollary 2.5]{CC}\label{thm:brownian1}
    As $L\to \infty$, $N\rightarrow\infty$, the following weak convergence holds in $\Omega$:
    \begin{equation}
    P_{t,L,N}^{*,0,0} \circ (\psi_{t,L})^{-1} \implies B_{e}(\cdot)
    \end{equation}
    where $B_{e}:\ [0,t]\rightarrow\R_{\ge 0}$ is the law of a Brownian excursion of length $t$.
\end{thm}

Theorem \ref{thm:brownian1} applies when we calculate the asymptotics of Eqn.~(\ref{eq_mixedmoment}) when $l=1$. 
For general $l\ge 1$, we also consider the conditional random walk bridges with different starting and ending heights. 

\begin{thm}\cite[Theorem. 2.4 and Corollary 2.5]{CC}\label{thm:brownian2}
    Suppose $H_{1}$ and $H_{2}$ are two sequences of integers such that 
    $$\frac{H_{1}}{\sigma\sqrt{L}}\rightarrow u>0,\qquad \frac{H_{2}}{\sigma\sqrt{L}}\rightarrow v>0$$
    as $L\rightarrow\infty$. Then as $L\rightarrow\infty$, $N\rightarrow\infty$, we have the weak convergence
    \begin{equation}
    P_{t,L,N}^{*,0,H_{1}} \circ (\psi_{t,L})^{-1} \implies B_{3}(\cdot),
    \end{equation}
    \begin{equation}
    P_{t,L,N}^{*,H_{2},0} \circ (\psi_{t,L})^{-1} \implies B_{3}(t-\cdot),
    \end{equation}
     where $B_{3}$ is a Brownian bridge of length $t$ from $0$ to $u$ conditioned on staying nonnegative, and $B_{3}(t-\cdot)$ is a Brownian bridge from $v$ to $0$ conditioned on staying nonnegative. 

    Moreover,
    \begin{equation}
    P_{t,L,N}^{*,H_{1},H_{2}} \circ (\psi_{t,L})^{-1} \implies B(\cdot),
    \end{equation}
   where $B$ is a Brownian bridge from $u$ to $v$ conditioned on staying nonnegative.
\end{thm}

\begin{figure}[h]
    \centering
\resizebox{\textwidth}{!}{
\begin{tabular}{ccc}
\begin{tikzpicture}[baseline = (current bounding box).center]
\draw (0,0)--(5,0);
\node[left] at (0,0) {$0$};
\zigzagpath[0.05]{0,1, 2, 3, 4, 5, 6, 7, 8, 9, 10, 11, 12, 11, 12, 13, 14, 13, 12, 13, 14, 15, 16, 17, 18, 19, 18, 17, 18, 19, 20, 21, 22, 23, 22, 21, 22, 23, 22, 21, 22, 23, 22, 21, 20, 21, 20, 19, 20, 19, 20, 19, 20, 19, 20, 19, 18, 17, 18, 19, 18, 17, 16, 15, 14, 15, 16, 17, 18, 19, 20, 19, 20, 21, 22, 21, 22, 23, 22, 23, 22, 21, 20, 21, 22, 21, 22, 21, 22, 23, 22, 21, 20, 21, 22, 21, 22, 21, 22, 21, 20};
\node[right] at (5,1) {$u$};
\end{tikzpicture} &
\begin{tikzpicture}[baseline = (current bounding box).center]
\draw (0,0)--(5,0);
\node[left] at (0,1) {$u$};
\zigzagpath[0.05]{20, 21, 20, 19, 20, 21, 20, 21, 20, 19, 18, 19, 18, 19, 20, 19, 18, 17, 16, 15, 14, 13, 14, 13, 12, 13, 14, 13, 12, 13, 14, 15, 14, 15, 14, 13, 14, 13, 12, 13, 14, 13, 12, 11, 10, 9, 8, 7, 6, 7, 6, 7, 6, 5, 6, 7, 8, 7, 6, 5, 6, 7, 8, 9, 10, 9, 10, 11, 12, 13, 12, 11, 10, 11, 12, 13, 12, 11, 12, 11, 12, 13, 12, 11, 10, 9, 8, 7, 8, 7, 6, 5, 6, 7, 6, 7, 6, 7, 8, 9, 10};
\node[right] at (5,0.5) {$v$};
\node at (0,0) {\phantom{0}};
\end{tikzpicture} &
\begin{tikzpicture}[baseline = (current bounding box).center]
\draw (0,0)--(5,0);
\node[left] at (0,0.5) {$v$};
\zigzagpath[0.05]{10, 9, 8, 9, 8, 7, 6, 7, 6, 5, 4, 3, 4, 5, 6, 7, 8, 9, 8, 7, 6, 5, 4, 5, 4, 5, 6, 7, 8, 7, 6, 7, 6, 7, 8, 7, 8, 9, 8, 9, 10, 9, 8, 9, 8, 7, 8, 9, 8, 7, 6, 5, 6, 5, 4, 5, 4, 3, 2, 3, 2, 3, 2, 1, 2, 1, 2, 3, 2, 3, 4, 3, 2, 3, 2, 3, 4, 5, 6, 5, 6, 7, 6, 7, 6, 7, 6, 5, 4, 5, 4, 3, 4, 3, 2, 3, 4, 3, 2, 1, 0};
\node[right] at (5,0) {$0$};
\node at (0,1) {};
\end{tikzpicture}
\end{tabular}
}
\caption{The three cases for the paths in Theorem \ref{thm:brownian2}.}\label{fig:splitting}
\end{figure}

\subsection{Conditional density functions}
We start with calculating the asymptotics of the partition function of a walk excursion. When the walk increments are Bernoulli random variables, as in \cite{GXZ}, this is given by (the asymptotics of) Catalan number.

\begin{lem}
    \label{lem:conditionalwalk}
     As $L\rightarrow\infty, N\rightarrow\infty$, we have

    $$ L^{\frac{3}{2}}\cdot \frac{|\mathcal{W}(0,0,L;N)|}{V_{N}(z_{c}(N))^{L}}\to\frac{1}{\sqrt{2\pi}}\frac{1}{P_{-1}\sigma}.$$
    
\end{lem}
\begin{proof}
   By Lemma \ref{lem:excursionpartition}, we have that 
    \begin{equation}\label{eq_contourexpression}
    |\mathcal{W}(0,0,L;N)|=\frac{1}{L+1} \frac{1}{2\pi i} \oint \; V_{N}(z)^{L+1} \; dz.
    \end{equation}
   The result then follows from a standard steepest descent analysis, which we sketch as follows.
   Denote $f_{N}(z)=\log V_{N}(z),$ which is analytic near $z=z_{c}(N)$ since $V_{N}(z)>0$ on $(0,\frac{1}{\alpha_{1}})$. Let $C^{\epsilon}_{z}(N)$ be the arc of the circle $C_{z}(N)$ with argument in $(-\frac{\epsilon}{\sqrt{L+1}},\frac{\epsilon}{\sqrt{L+1}})$. Then we have 
\begin{equation}\label{eq_steepestdescent}
    \begin{split}
        &\frac{1}{2\pi i}\oint_{C_{z}(N)}V_{N}(z)^{L+1}dz\stackrel{(a)}{\approx}\frac{1}{2\pi i}\oint_{C^{\epsilon}_{z}(N)}V_{N}(z)^{L+1}dz\\
        =&\frac{1}{2\pi i}\oint_{C^{\epsilon}_{z}(N)}\exp\Big[(L+1)f_{N}(z)\Big]dz\\
        \stackrel{(b)}{\approx}&\frac{1}{2\pi i}\oint_{C^{\epsilon}_{z}(N)}\exp\Big[(L+1)f_{N}(z_{c}(N))\Big]\exp\Big[(L+1)\frac{f''_{N}(z_{c}(N))}{2}(z-z_{c}(N))^{2}\Big]dz\\
        \stackrel{(c)}{\approx}&V_{N}(z_{c}(N))^{L+1}\cdot\frac{1}{2\pi i}\oint_{C^{\epsilon}_{z}(N)}\exp\Big[(L+1)\frac{f''(z_{c})}{2}(z-z_{c})^{2}\Big]dz\\
        \stackrel{(d)}{\approx}&\frac{V_{N}(z_{c}(N))^{L+1}}{\sqrt{L+1}}\cdot \frac{1}{2\pi}\int_{-\infty}^{\infty}\exp\Big[-\frac{|f''(z_{c})|}{2}t^{2}\Big]dt\\
       \approx&\frac{V(z_{c}(N))^{L}}{\sqrt{L+1}}\cdot \frac{1}{\sqrt{2\pi}}\frac{V(z_{c})}{\sqrt{|f''(z_{c})|}}.
    \end{split}
\end{equation}
The approximations (a)-(d) require some verification, whose details we omit. Since $$\frac{V(z_{c})}{\sqrt{|f''(z_{c})|}}=\frac{V(z_{c})^{3/2}}{\sqrt{V''(z_{c})}}=\frac{1}{P_{-1}\sigma}$$
by Eqn.~(\ref{eq_sigma}), the result then follows.
\end{proof}

\begin{definition}\label{def:densities}
    For $x>0$, $h,g> 0$, we define 
    \[\mathbf{F}(x;h,g)=\frac{1}{\sqrt{2\pi x}}\left(\exp\left(-\frac{(g-h)^{2}}{2x}\right)-\exp\left(-\frac{(g+h)^{2}}{2x}\right)\right).\]
    Note that (by reflection principle) $\mathbf{F}(x;h,g)$ is the transition density of a Brownian bridge from $g$ to $h$, conditioned on staying nonnegative.

    Let \[\mathbf{F}_{0}(x,h)=\frac{2h}{\sqrt{2\pi x^{3}}}\exp\left(-\frac{h^{2}}{2x}\right),\]
and let \[\mathbf{F}_{0,0}(x)=\frac{2}{\sqrt{2\pi x^{3}}}.\]
\end{definition}
Using Lemma \ref{lem:conditionalwalk}, we have a precise estimate of the partition functions of $\mathcal{W}(H,0,L;N)$ and $\mathcal{W}(0,H,L;N)$, in terms of $\mathbf{F}_{0}$ and $\mathbf{F}_{0,0}$. Note that the asymptotic expressions depend on the parameters $P_{-1}$ and $\sigma$ defined in Eqn.~(\ref{eq_sigma}), and differ in the two cases.
\begin{lem}\label{lem:conditionalwalk2}
    Fix $h>0$, $x>0$. As  $N\rightarrow\infty$, we have
    \begin{equation}\label{eq_conditionalwalk}
    \begin{split}
        &N^{2/3}\cdot\frac{|\mathcal{W}(\lfloor h\sigma N^{1/3}\rfloor,0,\lfloor xN^{2/3}\rfloor ;N)|}{V_{N}(z_{c}(N))^{\lfloor xN^{2/3}\rfloor}}\rightarrow (2P_{-1})^{-1}\cdot \mathbf{F}_{0}(x,h),\\
        &N^{2/3}\cdot\frac{|\mathcal{W}(0,\lfloor h\sigma N^{1/3}\rfloor,\lfloor xN^{2/3}\rfloor ;N)|}{V_{N}(z_{c}(N))^{\lfloor xN^{2/3}\rfloor}}\rightarrow \sigma^{-2}\cdot \mathbf{F}_{0}(x,h),\\
        &N\cdot \frac{|\mathcal{W}(0,0,\lfloor xN^{2/3}\rfloor;N)|}{V_{N}(z_{c}(N))^{\lfloor x N^{2/3}\rfloor}}\to(2P_{-1}\sigma)^{-1}\mathbf{F}_{0,0}(x).
        \end{split}
    \end{equation}
\end{lem}
\begin{proof}
Note that by our assumptions, the random walk law is supported by $\Z$ and aperiodic. Then by the classical local limit theorem (see e.g. \cite[Theorem 2]{DM}) and Lemma \ref{lem:excursionpartition}, 
\begin{equation}\label{eq_conditionalwalk1}
\begin{split}
  &\frac{|\mathcal{W}(\lfloor h\sigma N^{1/3}\rfloor,0,\lfloor xN^{2/3}\rfloor ;N)|}{V_{N}(z_{c}(N))^{\lfloor xN^{2/3}\rfloor}}\frac{1}{z_{c}(N)V_{N}(z_{c}(N))}\\
  =&\frac{|\mathcal{W}(\lfloor h\sigma N^{1/3}\rfloor,0,\lfloor xN^{2/3}\rfloor ;N)|z_{c}(N)^{-1}}{|\tilde{\mathcal{W}}(\lfloor h\sigma N^{1/3}\rfloor+1,0,\lfloor xN^{2/3}\rfloor+1 ;N)|} \frac{|\tilde{\mathcal{W}}(\lfloor h\sigma N^{1/3}\rfloor+1,0,\lfloor xN^{2/3}\rfloor+1 ;N)|}{V_{N}(z_{c}(N))^{\lfloor xN^{2/3}\rfloor+1}}\\
  =&\frac{\lfloor h\sigma N^{1/3}\rfloor+1}{\lfloor xN^{2/3}\rfloor+1}\cdot \frac{1}{\sigma  \sqrt{xN^{2/3}}}\left(\frac{1}{\sqrt{2\pi}}\exp\left(-\frac{h^{2}}{2x}\right)+o(1)\right).\\
 \end{split}
\end{equation}
The first line of Eqn.~(\ref{eq_conditionalwalk}) then follows from dividing both sides by $P_{-1}(N)$ and taking $N\rightarrow\infty$.

For the second line, consider a random walk excursion $W'(t)$ of length $2\lfloor xN^{2/3}\rfloor$,  with i.i.d. increments $X(N)$ as defined in Eqn.~(\ref{eq_rvx}), conditioned to stay non-negative. The midpoint density at $t=\frac{1}{2}$ of a standard Brownian excursion of length 1 is\footnote{One can see this by e.g. a limit transition from the Bernoulli walks.} $$ \frac{16}{\sqrt{2\pi}} y^{2}\exp\left(-2y^{2}\right).$$ So by Theorem  \ref{thm:brownian1} and the local limit theorem for non-negative walks (see e.g. \cite[Proposition 4.1]{CC}),
\begin{equation}\label{eq_conditionalwalk2}
\begin{split}
    &\frac{|\mathcal{W}(0,\lfloor h\sigma N^{1/3}\rfloor,\lfloor xN^{2/3}\rfloor ;N)|\cdot|\mathcal{W}(\lfloor h\sigma N^{1/3}\rfloor,0,\lfloor xN^{2/3}\rfloor ;N)|}{|\mathcal{W}(0,0,2\lfloor xN^{2/3}\rfloor;N)|}=\PP[W'(\lfloor xN^{2/3}\rfloor)=\lfloor h\sigma N^{1/3}\rfloor]\\
    =&\frac{1}{\sigma }\frac{1}{\sqrt{xN^{2/3}}}\left(\frac{4h^{2}}{\sqrt{\pi}x}\exp\left(-\frac{h^{2}}{x}\right)+o(1)\right), 
\end{split}
\end{equation}
By Lemma \ref{lem:conditionalwalk}, \begin{align*}
    \frac{|\mathcal{W}(0,0,2\lfloor xN^{2/3}\rfloor;N)|}{V_{N}(z_{c}(N))^{2\lfloor xN^{2/3}\rfloor}}=\frac{1}{\sqrt{2\pi}}\left(\frac{1}{P_{-1}\sigma}+o(1)\right)(2\lfloor xN^{2/3}\rfloor)^{-\frac{3}{2}}.
\end{align*} 
Multiplying this with Eqn.~(\ref{eq_conditionalwalk2}) gives 
\begin{align*}
    \frac{|\mathcal{W}(0,\lfloor h\sigma N^{1/3}\rfloor,\lfloor xN^{2/3}\rfloor ;N)|\cdot|\mathcal{W}(\lfloor h\sigma N^{1/3}\rfloor,0,\lfloor xN^{2/3}\rfloor ;N)|}{V_{N}(z_{c}(N))^{2\lfloor xN^{2/3}\rfloor}}=\frac{1}{\pi}\frac{1}{P_{-1}\sigma^{2} }\frac{1}{(xN^{2/3})^{2}}\left(\frac{h^{2}}{x}\exp\left(-\frac{h^{2}}{x}\right)+o(1)\right).
\end{align*}
The second line of Eqn.~(\ref{eq_conditionalwalk}) then follows from further dividing the above expression by the first line. The third line is a rephrase of Lemma \ref{lem:conditionalwalk}.
\end{proof}

\subsection{Time homogeneity of down-steps}
Recall that the only negative increment of our walk-bridges is -1, and we say that the walk takes a down-step at a certain time when this occurs. 
\begin{lem}\label{lem:timehom}
For $H_{1}, H_{2}\in \Z_{\ge 0},\ L\in\Z_{\ge 1}$, sample the walk $W(t)$ from $\mathcal{W}(H_{1},H_{2},L,N)$, with i.i.d.. increments as in Eqn.~(\ref{eq_rvx}), where $$\frac{H_{1}}{\sqrt{L}}\rightarrow u,\ \frac{H_{2}}{\sqrt{L}}\rightarrow v$$ for some given $u,v\ge0$.
     For any $0<t_{1}<t_{2}<1$, let $\rho_{-1}(t_1,t_2)$ denote the proportion of  down-steps in $\llbracket t_{1}L,t_{2}L\rrbracket$. Then we have 
    \begin{equation}
        \rho_{-1}(t_1,t_2)\rightarrow P_{-1}\quad \text{in\ probability}
    \end{equation}
    as $L, N\rightarrow\infty$, where $P_{-1}=\frac{1}{z_{c}V(z_{c})}$, as defined in Eqn.~(\ref{eq_sigma}).
\end{lem}
\begin{remark}\label{rem:p-1}
    By the definition of $X(N)$, we recall from Section \ref{sec:conditionalrw} that
    \begin{equation}
        P_{-1}=\lim_{N\rightarrow\infty}\PP[X(N)=-1].
    \end{equation}
\end{remark}
\begin{proof}
    Since $t_{1}>0, t_{2}<1$, $H_{1},H_{2}\sim \sqrt{L}$, under the scaling of $W(t)$ according to Theorem \ref{thm:brownian1} or \ref{thm:brownian2}, we have that with high probability, $$ \frac{1}{\sqrt{L}}W(\lfloor t_{1}L\rfloor),\ \frac{1}{\sqrt{L}}W(\lfloor t_{2}L\rfloor)>\epsilon>0$$
    for some $\epsilon(t_{1},t_{2})$ uniform in $L,N$ large enough. Therefore it suffices to take $t_{1}=0$, $t_{2}=1$, and $u,v>0$.

   Let $\tilde{W}(t)(t=0,1,2,\ldots,L)$ be a walk with the same increments and weights as $W(t)$, except that it does not have the conditioning to stay nonnegative. Observe that 
    \begin{equation}
        \PP[\tilde{W}(t)\ge 0\ \text{for\ all\ }t|\tilde{W}(0)=H_{1},\tilde{W}(L)=H_{2}]>\epsilon'
    \end{equation}
    for some $\epsilon'=\epsilon'(u,v)>0$ uniformly in $L,N$ large enough. 
    Therefore, we have for any $\delta>0$, 
    \begin{equation}
        \begin{split}
        &\PP\left[\left|\rho_{-1}-P_{-1}\right|>\delta\Big| \tilde{W}(t)\ge 0\ \text{for\ all\ }t\right] =\frac{\PP\left[\left|\rho_{-1}-P_{-1}\right|>\delta, \tilde{W}(t)\ge 0\ \text{for\ all\ }t\right]}{\PP\left[\tilde{W}(t)\ge 0\ \text{for\ all\ }t\right]}\\
        \le&\frac{1}{\epsilon}\PP\left[\left|\rho_{-1}-P_{-1}\right|>\delta\right].\\
        \end{split}
    \end{equation}
    It remains to show that the last line above goes to 0 as $L\rightarrow\infty$, $N\rightarrow\infty$. This is done by a straightforward calculation. 
    
    Denote $H_{2}-H_{1}$ by $H(L)$,  $L+H(L)$ by $H'(L)$. Then for $n_{-1},n_{0},n_{1}, n_{2},\ldots\in \Z_{\ge 0}$, 
    \begin{equation}\label{eq_constriant1}
    n_{-1}+n_{0}+n_{1}+\ldots+n_{H'(L)}=L,
    \end{equation} 
    where $J_i$ is the number of jumps of size $i$,  and 
    \begin{equation}\label{eq_constraint2}
        -n_{-1}+n_{1}+2n_{2}+\ldots+H'(L)n_{H'(L)}=H(L).
    \end{equation}
  We have 
    \begin{equation}\label{eq_distfunction}
    \begin{split}
        &\PP\left[\text{there\ are\ } n_{i}\ \text{jumps\ of\ size}\ i\ \text{in\ the\ walk, for each $i=-1,0,1,\ldots, H'(L)$ }\right]\\
        =&\binom{L}{n_{-1}\ n_{0}\ n_{1}\ \ldots\ n_{H'(L)}}p_{-1}(N)^{n_{-1}}p_{0}(N)^{n_{0}}p_{1}(N)^{n_{1}}\cdots p_{H'(L)}(N)^{n_{H'(L)}}\\
        \approx&(2\pi)^{-\frac{H'(L)}{2}}L^{L+\frac{1}{2}}n_{-1}^{-(n_{-1}+\frac{1}{2})}n_{0}^{-(n_{0}+\frac{1}{2})}\cdots n_{H'(L)}^{-(n_{H'(L)}+\frac{1}{2})}\\
        \qquad\quad&\cdot p_{-1}(N)^{n_{-1}}p_{0}(N)^{n_{0}}\cdots p_{H'(L)}(N)^{n_{H'(L)}}.
        \end{split}
    \end{equation}
   The third line of Eqn.~(\ref{eq_distfunction}) follows from Stirling's formula.

    By Eqn.~(\ref{eq_constriant1}) and (\ref{eq_constraint2}), we rewrite $n_{-1},n_{1}$ in terms of the other terms, such that 
    \begin{equation}\label{eq_k-1k1}
    \begin{split}
        n_{-1}=&\frac{1}{2}\left[L-n_{0}+n_{2}+2n_{3}+\ldots+(H'(L)-1)n_{H'(L)}\right]-\frac{1}{2}H(L)\\
        n_{1}=&\frac{1}{2}\left[L-n_{0}-3n_{2}-4n_{3}-\ldots-(H'(L)+1)n_{H'(L)}\right]+\frac{1}{2}H(L).\\
     \end{split}
    \end{equation}

    Our goal is to find $(n_{-1},n_{0},n_{1},\ldots)$ that maximize the probability in Eqn.~\eqref{eq_distfunction}. We take the logarithm of the right-hand side of Eqn.~(\ref{eq_distfunction}), plug in the expression in Eqn.~(\ref{eq_k-1k1}), and take the partial derivatives of the resulting expression in $n_{i}$. Then  after taking $\frac{\partial}{\partial n_{i}}\ (i=0,2,3,\ldots)$ we have:
    \begin{equation}\label{eq_derivative}
        \begin{split}
           0= &\log(p_{i}(N))-\frac{i+1}{2}\log(p_{1}(N))+\frac{i-1}{2}\log(p_{-1}(N))\\
            -&\left[\log(n_{i})+\frac{n_{i}+\frac{1}{2}}{n_{i}}-\frac{i+1}{2}\log(n_{1})-\frac{n_{1}-\frac{1}{2}}{n_{1}}\frac{i+1}{2}+\frac{i-1}{2}\log(n_{-1})+\frac{n_{-1}+\frac{1}{2}}{n_{-1}}\frac{i-1}{2}\right]\\
             =&\log(p_{i}(N))-\frac{i+1}{2}\log(p_{1}(N))+\frac{i-1}{2}\log(p_{-1}(N))\\
            -&\Big[\log(\frac{n_{i}}{L})+\frac{n_{i}+\frac{1}{2}}{n_{i}}-\frac{i+1}{2}\log(\frac{n_{1}}{L})-\frac{n_{1}
            -\frac{1}{2}}{n_{1}}\frac{i+1}{2}+\frac{i-1}{2}\log(\frac{n_{-1}}{L})+\frac{n_{-1}+\frac{1}{2}}{n_{-1}}\frac{i-1}{2}\Big]\\
            \approx&\log(p_{i}(N))-\frac{i+1}{2}\log(p_{1}(N))+\frac{i-1}{2}\log(p_{-1}(N))\\
            -&\left[\log(n_{i})-\frac{i+1}{2}\log(n_{1})+\frac{i-1}{2}\log(n_{-1})\right],
        \end{split}
    \end{equation}
where in the approximation in the last line we take $\frac{n_{i}+\frac{1}{2}}{n_{i}}\approx 1$ for $i=-1,0,1,2,\ldots$

There are $H'(L)-1$ (approximately linear) equations with $H'(L)-1$ unknowns $\log (n_{0})$, $\log(n_{2}),\ldots,\log(n_{H'(L)})$. We observe that for any fixed $i=0,2,\ldots$
$$\frac{n_{i}}{L}\approx p_{i}(N)$$ solves Eqn.~(\ref{eq_derivative}) as $L,N\rightarrow\infty$. Plugging these back to Eqn.~(\ref{eq_k-1k1}) gives $\frac{n_{-1}}{L}\approx P_{-1}$.
It is straightforward to see from Eqn.~(\ref{eq_distfunction}) that as $L\rightarrow\infty$, $N\rightarrow\infty$ the distribution function concentrates near its maximum, which in particular implies that for any $\delta>0$,
\begin{equation*}
    \PP\left[\left|\frac{n_{-1}}{H'(L)}-\PP\left[X_{1}(N)=-1\right]\right|>\delta\right]\rightarrow 0.
\end{equation*}
Since $\rho_{-1}=\frac{n_{-1}}{L}$, $X_{1}(N)\rightarrow X_{1}$ in distribution, the result follows.
    
\end{proof}

\subsection{Tail estimate of the walks}\label{sec:tailrw}

 We give several results about the extreme values of random walks, which will be  used in the tail estimate argument in Section \ref{sec:tailestimates}.

 Let $$\mathcal{\tilde{W}}(0,\cdot, L;N)=\bigcup_{H=0}^{\infty}\mathcal{\tilde{W}}(0,H,L;N),\quad \mathcal{\tilde{W}}(\cdot, H,L;N)=\bigcup_{H'=0}^{\infty}\mathcal{\tilde{W}}(H',H,L;N).$$ Here, $\mathcal{\tilde{W}}(0,\cdot, L;N)$ denotes the set of walks starting at $W(0)=0$ and ending at an arbitrary $W(L)=H \in \Z_{\ge 0}$. Similarly, $\mathcal{\tilde{W}}(\cdot, H,L;N)$ denotes the set of walks starting at an arbitrary $W(0)=H' \in \Z_{\ge 0}$ and ending\footnote{One can view $W$ as a time-reversed walk in this case.} at $W(L)=H$, for $H\in \Z_{\ge 0}$. These sets are equipped with the same increment weights as in Eqn.~(\ref{eq_rvx}).

\begin{lem} \label{lem:basicwalk}
    For $H\in \Z_{\ge 0},\ L\in \Z_{\ge 100}$, sample the walk $W(t)$ from $\mathcal{\tilde{W}}(0,\cdot,L;N)$ and  $\mathcal{\tilde{W}}(\cdot, 0,L;N)$ respectively.
    Then there exist constants $C>0$, $c>0$ uniform in $L,N$, such that for all $h\ge 0$,
    \begin{align*}
       \PP[W\ge 0, W(L)=\lfloor h\sqrt{L}\rfloor]\le \frac{C}{L}\quad &\text{if}\ \ W(t)\in \tilde{\mathcal{W}}(0,\cdot,L;N),\\
        \PP[W\ge 0, W(0)=\lfloor h\sqrt{L}\rfloor]\le \frac{C}{L}\exp(-ch^{2})\quad &\text{if}\ \ W(t)\in \tilde{\mathcal{W}}(\cdot ,0,L;N).
    \end{align*}
\end{lem}

\begin{proof}
    For $W(t)\in \tilde{\mathcal{W}}(0,\cdot,L;N)$ or $ \tilde{\mathcal{W}}(\cdot,0,L;N)$ we have\footnote{This can be obtained from the estimate of the  first ladder epoch $T_{1}$, see e.g. \cite[(2.3) and (2.5)]{C}.} $\PP[W(t)\ge 0, 0\le t\le L]\le \frac{C}{\sqrt{L}}$ for some $C>0$. Then the first bound of the lemma is a direct corollary of \cite[Theorem 2]{C}.  On the other hand, for  $W(t)\in \tilde{W}(\cdot,0,L;N)$, by Proposition \ref{prop:Ballot}  we have 
    \[\frac{\PP[W\ge 0,W(0)=\lfloor h\sqrt{L}\rfloor]}{\PP[W(0)=\lfloor h\sqrt{L}\rfloor]}\le \frac{Ch}{\sqrt{L}}.\] Therefore it suffices to prove 
    \[\PP[W(0)=\lfloor h\sqrt{L}\rfloor]\le \frac{C}{\sqrt{L}}\exp(-ch^{2}).\]
    Note that since the increment $X(N)\ge -1$, $W(0)\le L$. For $0\le h\le \rho \sqrt{L}$ where $\rho>0$ is some small enough constant, this follows directly from \cite[Theorem 2.3.11]{LL}. For $\rho \sqrt{L}\le h\le \sqrt{L}$, it suffices for us to prove  $$ \PP[W(0)=\lfloor h\sqrt{L}\rfloor]\le \frac{C}{\sqrt{L}}\exp(-cL).$$
    This can be proved by slightly adapting the classical characteristic function argument in \cite[Theorem 2.3.12]{LL}, whose details we omit.
\end{proof}

\begin{lem}\label{lem:walkmax}For $L\in \Z_{\ge 100}$, take discrete intervals $[\![X,Y]\!]$ contained in $[\![0,L]\!]$ where $Y-X\ge 1$.
Consider three types of random walks $W(t):[\![0,L]\!]\rightarrow \Z$ as follows. Let $V=W(X)-\min_{[\![X,Y]\!]}W+1$.

(1).    For $H\in \Z_{\ge 0}$, sample the walk  $W(t)$ in $\tilde{\mathcal{W}}(\cdot,H,L;N)$.
Then there exist constants $C>0$, $c>0$, which are uniform in $H,L,N$, such that for all $h\ge 0$, $\tilde{h}\ge 0$,
    $$\PP[W\ge 0,\ \max_{t} W(t)-H>h\sqrt{L},\ V>\tilde{h}\sqrt{Y-X}]\le\begin{cases}
         C\exp(-ch^{2})\exp(-c\tilde{h}^{2}),\quad &\text{if}\ \ H>0\\
         \frac{C}{\sqrt{L}}\exp(-ch^{2})\exp(-c\tilde{h}^{2}),\quad &\text{if}\ \ H=0.
    \end{cases}$$

(2). For $H\in \Z_{\ge 1}$, sample the walk  from $\tilde{\mathcal{W}}(\cdot,H,L,N)$.
Then there exist constants $C>0$, $c>0$, which are uniform in $L,N$, such that for all $h\ge 0$, 
    $$\PP[W(0)=0,\ W\ge 0,\ \max_{t} W(t)-H>h\sqrt{L},\ V>\tilde{h}\sqrt{Y-X}]\le\frac{C}{L}\exp(-ch^{2})\exp(-c\tilde{h}^{2}).$$

(3). Sample the walk from $\tilde{\mathcal{W}}(\cdot,0,L,N)$. 
Then there exist constants $C>0$, $c>0$, which are uniform in $L,N$, such that for all $h\ge 0$,
    $$\PP[W(0)=0,\ W\ge 0,\ \max_{t} W(t)>h\sqrt{L},\ V>\tilde{h}\sqrt{Y-X}]\le\frac{C}{L^{3/2}}\exp(-ch^{2})\exp(-c\tilde{h}^{2}).$$
\end{lem}
\begin{proof}
 In the proof we take $h\ge 1$, since for $0\le h\le 1$ the factor $\exp(-ch^{2})$ is absorbed in the constant $C$, and there is no need to introduce the quantity $T$ below. Let \begin{align*}\tilde{T}=&\max_{X\le t\le Y}\{\tilde{t}:\ W(\tilde{t})=\min_{X\le t\le Y}W(t)\},\\
 T=&\max\Big(\max_{0\le t\le X}\{t:\ W(t)=W(X)+\lfloor h\sqrt{L}/3\rfloor+1\},\max_{X+1\le t\le \tilde{T}}\{t:\ W(t)=W(\tilde{T})+\lfloor h\sqrt{L}/3\rfloor+1\},\\ 
 &\max_{\tilde{T}+1\le t\le L}\{t:\ W(t)=W(L)+\lfloor h\sqrt{L}/3\rfloor+1\},\Big),
 \end{align*}
 and define the events \[E_{T}=\{W\ge 0\ on\  \llbracket T,L\rrbracket \},\]
if they exist. 

 We first consider case (1), $H>0$.  We have \begin{align*}
     &\PP[W\ge 0,\ \max_{t} W(t)-H>h\sqrt{L},\ V>\tilde{h}\sqrt{Y-X}]  \\
     \le&\PP[X\le \tilde{T}<T\le L|E_{T}]+\PP[X\le T<\tilde{T}\le Y|E_{T}]+\PP[  T< X\le \tilde{T}\le L|E_{T}].
 \end{align*}
 Note that we currently omit the factor $\PP[E_{T}]$.
 It then suffice to get the desired bound for all the above three
 cases. Order the tuple $(X,\tilde{T},T)$ as $(X_{1}\le X_{2}\le X_{3})$, and let $W_{1}(t)=W(t)|_{0\le t\le X_{1}}$,  $W_{2}(t)=W(t)|_{X_{1}\le t\le X_{2}}$,  $W_{3}(t)=W(t)|_{X_{2}\le t\le X_{3}}$,  $W_{4}(t)=W(t)|_{X_{3}\le t\le L}$. 
Then we have \begin{align*}
   & \PP[ X\le \tilde{T}<T\le L|E_{T}]\le\sum_{\tilde{t}=X}^{Y}\sum_{X_{3}=X_{2}+1}^{Y}\PP[W_{2}(X_{1})-W_{2}(X_{2})> \tilde{h}\sqrt{Y-X}| X_{2}=\tilde{t},\ E_{T}]\\
    \times&\PP[X_{2}=\tilde{t}|E_{T}] \cdot \PP[W_{4}\le W_{4}(X_{3}),\ W_{4}(L)-W_{4}(X_{3})=-\lfloor h\sqrt{L}/3\rfloor|E_{T}]\\
    \le&\sum_{\tilde{t}=X}^{Y}\sum_{X_{3}=X_{2}+1}^{Y}\exp(-c\tilde{h}^{2})\cdot\PP[X_{2}=\tilde{t}|E_{T}] \cdot \frac{C}{L}\exp(-ch^{2})\le C\exp(-c\tilde{h}^{2})\exp(-ch^{2}),
\end{align*}
where we use Lemma \ref{lem:basicwalk} in the second inequality.

Similarly, \begin{align*}
   &\PP[X\le T<\tilde{T}\le Y|E_{T}]\\
   \le&\sum_{\tilde{t}=X}^{Y}\sum_{X_{2}=X}^{\tilde{t}-1}
 \Big( \PP[W_{3}(X_{2})\ge W_{3},\ W_{3}(X_{2})-W_{3}(X_{3})> \tilde{h}\sqrt{Y-X}/2,\\ 
 &\quad \quad\quad \quad \quad\quad \quad\quad \quad \quad\quad  W_{3}(X_{2})-W_{3}(X_{3})= \lfloor h\sqrt{L}/3\rfloor|X_{3}=\tilde{t},\ E_{T} ]\cdot \PP[X_{3}=\tilde{t}|E_{T}]\\
    +& \PP[ W_{2}(X_{1})-W_{2}(X_{2})> \tilde{h}\sqrt{Y-X}/2|E_{T}]
   \times\PP[W_{3}(X_{2})\ge W_{3},\ W_{3}(X_{2})-W_{3}(X_{3})= \lfloor h\sqrt{L}/3\rfloor|X_{3}=\tilde{t},\ E_{T} ]\\
   \times&\PP[X_{3}=\tilde{t}|E_{T}]\Big)\\
\le&\sum_{\tilde{t}=X}^{Y}\sum_{X_{2}=X}^{\tilde{t}-1}\frac{C}{L}\exp(-c\tilde{h}^{2})\exp(-ch^{2})\cdot \PP[X_{3}=\tilde{t}|E_{T}]+\sum_{\tilde{t}=X}^{Y}\sum_{X_{2}=X}^{\tilde{t}-1}\exp(-c\tilde{h}^{2})\cdot \frac{C}{L}\exp(-ch^{2})\cdot \PP[X_{3}=\tilde{t}|E_{T}]\\
    \le& C\exp(-c\tilde{h}^{2})\exp(-ch^{2}),
\end{align*}
and \begin{align*}
   &\PP[ T< X\le \tilde{T}\le L|E_{T}]\le\sum_{\tilde{t}=X}^{Y}\sum_{X_{1}=0}^{X-1}
  \PP[W_{2}\le W_{2}(X_{1}),\ W_{2}(X_{2})-W_{2}(X_{1}))=-\lfloor h\sqrt{L}/3\rfloor|W\ge0 ]\\
   \times&\PP[ W_{3}(X_{2})-W_{3}(X_{3})>\tilde{h}\sqrt{Y-X}|X_{3}=\tilde{t},\ E_{T}] \times \PP[X_{3}=\tilde{t}|E_{T}]\\
    \le&\sum_{\tilde{t}=X}^{Y}\sum_{X_{1}=0}^{X-1}\frac{C}{L}\exp(-ch^{2})\cdot \exp(-c\tilde{h}^{2})\cdot \PP[X_{3}=\tilde{t}|E_{T}]\le C\exp(-c\tilde{h}^{2})\exp(-ch^{2}).
    \end{align*}
    This finishes the proof for case (1), $H>0$.

    For case (1) where $H=0$, the same argument as above still holds, while we upgrade the conditioning event $E_{T}$ as $\{W\ge 0\ \text{on}\ \llbracket T,L\rrbracket| W(L)=0\}$, which gives one more factor \[\PP[W\ge 0\ \text{for}\ T\le t\le L| W(L)=0]\le \frac{C}{\sqrt{L}}.\] For case (2), we still condition on $E_{T}$, and there is an extra restriction on $W(t)|_{0\le t\le T}$ \[W(0)=0,\ W\ge 0\ \text{on}\ \llbracket 0,T\rrbracket\setminus\llbracket X,\tilde{T}\rrbracket,\] 
    which gives a factor $\le \frac{C}{L}$ to $\PP[W\ge 0,\ \max_{t} W(t)-H>h\sqrt{L},\ V>\tilde{h}\sqrt{Y-X}| E_{T}]$ by the second line of Lemma \ref{lem:basicwalk}.
    For case (3) we have both the upgraded conditioning event and the extra restriction on $W(t)|_{0\le t\le T}$.
    Combining the extra factors with the previous bounds  gives all the desired estimates.
\end{proof}

\begin{lem}\label{lem:max-min}
Take discrete intervals $[\![X_{1},Y_{1}]\!]$,\ldots,$[\![X_{k},Y_{k}]\!]$ contained in $[\![0,L]\!]$. Assume that $Y_{i}-X_{i}\ge 1$ and that the intervals are mutually disjoint. For each $i=1,2,\ldots,k$, denote $V_{i}=W(X_{i})-\min_{[\![X_{i},Y_{i}]\!]}W+1$. Take $a_{1},\ldots,a_{k}\in \Z_{>0}$ and let $a=a_{1}+\ldots+a_{k}$. Consider the following two cases:

(1). For $H\in \Z_{\ge 1},\ L\in \Z_{\ge 1}$, sample the walk  $W(t)$ in $\tilde{\mathcal{W}}(\cdot,H,L;N)$. 
Then there exist constants $C>0$, $c>0$ uniform in $H,L,N$, such that for all $h\ge 0$,
    \begin{align*}
        \E\left[\prod_{i=1}^{k}\frac{V_{i}^{a_{i}}}{a_{i}!}\mathbf{1}_{W\ge 0}\mathbf{1}_{\max W> h\sqrt{L}}\right] 
         <C^{a}\log(a+2)^{a+2}a^{-\frac{a}{2}}L^{\frac{a}{2}}\exp(-ch^{2}),\\
         \E\left[\prod_{i=1}^{k}\frac{V_{i}^{a_{i}}}{a_{i}!}\mathbf{1}_{W\ge 0,\ W(0)=0}\mathbf{1}_{\max W> h\sqrt{L}}\right] 
         <C^{a}\log(a+2)^{a+2}a^{-\frac{a}{2}}L^{\frac{a}{2}}\exp(-ch^{2})\cdot L^{-1}.
    \end{align*}

(2). Sample the walk  $W(t)$ in $\tilde{\mathcal{W}}(\cdot,0,L;N)$.
Then there exist constants $C>0$, $c>0$ uniform in $H,L,N$, such that for all $h\ge 0$,
    \begin{align*}
        \E\left[\prod_{i=1}^{k}\frac{V_{i}^{a_{i}}}{a_{i}!}\mathbf{1}_{W\ge 0}\mathbf{1}_{\max W> h\sqrt{L}}\right] 
         <C^{a}\log(a+2)^{a+2}a^{-\frac{a}{2}}L^{\frac{a}{2}}\exp(-ch^{2})\cdot L^{-1/2},\\
         \E\left[\prod_{i=1}^{k}\frac{V_{i}^{a_{i}}}{a_{i}!}\mathbf{1}_{W\ge 0,\ W(0)=0}\mathbf{1}_{\max W> h\sqrt{L}}\right] 
         <C^{a}\log(a+2)^{a+2}a^{-\frac{a}{2}}L^{\frac{a}{2}}\exp(-ch^{2})\cdot L^{-3/2}.
    \end{align*}  
  
\end{lem}
\begin{remark}
    A similar statement to this lemma is proved for Bernoulli walks, see \cite[Lemma A.7]{GXZ}. A key difference is that to deal with our walks, whose increment is unbounded in the upward direction, we define $V_{i}$ to be $W(X_{i})-\min_{\llbracket X_{i}, Y_{i}\rrbracket}W+1$ instead of the larger quantity $\max_{\llbracket X_{i},Y_{i}\rrbracket}W-\min_{\llbracket X_{i}, Y_{i}\rrbracket}W+1$, which is a  sufficient upper bound for the number of last jumps (defined in Section \ref{sec:classify})  in $\llbracket X_{i},Y_{i}\rrbracket$. We also make the observation that instead of considering the more general walk bridges in $\mathcal{W}(H_{1},H_{2},L;N)$ for arbitrary $H_{1}, H_{2}\in \Z_{\ge 0}$, if suffices to consider $H_{1}=0$, or consider a walk meander with undetermined left endpoint. This allows us to drop the drift term that appears in the upper bound of \cite[Lemma A.7]{GXZ}.
\vspace{5pt}
\end{remark}
 \begin{proof} 
We give the proof for the first bound in case (1), and the general cases follow from exactly the same argument.
Let $Z_{i}=\frac{V_{i}}{\sqrt{Y_{i}-X_{i}}}$ for $i=1,2,\ldots,k$ and let $Z=\max_{i}Z_{i}$. By Lemma \ref{lem:walkmax} and a union bound in $i$ we have $\PP[Z>h]<Ck\exp(-ch^{2})
$ for all $h>0$. Therefore summing over $h=0,1,2,\ldots$ gives $$\E[Z^{m}]<(C \log(k)m)^{m/2}\le (C\log(a)m)^{m/2}$$ for any $m\in \Z_{\ge 0}$.

Upper bound each $V_{i}$ by $Z\sqrt{Y_{i}-X_{i}}$. Then by the above inequality, there exists constant $C>0$ uniform in $H,L,N$, $W(Y_{i})$, such that 
\begin{align*}
    \E\left[\prod_{i=1}^{k}V_{i}^{a_{i}}\right]\le C^{a}\log(a+2)^{a+2}\prod_{i=1}^{k}a_{i}^{a_{i}/2}(Y_{i}-X_{i})^{a_{i}/2}.
\end{align*}
Using the bound $a_{i}!\ge c^{a_{i}}a_{i}^{a_{i}}$ for some uniform constant $c>0$, we have 
\begin{align*}
    \E\left[\prod_{i=1}^{k}\frac{V_{i}^{a_{i}}}{a_{i}!}\right]\le C^{a}\prod_{i=1}^{k}\left(\frac{Y_{i}-X_{i}}{a_{i}}\right)^{a_{i}/2}.
\end{align*}
An elementary calculation shows that $$\prod_{i=1}^{k}(Y_{i}-X_{i})^{a_{i}/2}\le \frac{\prod_{i=1}^{k}a_{i}^{a_{i}/2}}{a^{a/2}}L^{a/2}.$$
Therefore $$\E\left[\prod_{i=1}^{k}\frac{V_{i}^{a_{i}}}{a_{i}!}\right]\le C^{a}\log(a+2)^{a+2}a^{-a/2}L^{a/2}.$$
\end{proof}

\section{Asymptotics}\label{sec:asymptotics}

Recall from Section \ref{sec:configuration} that we are considering the action
\begin{align*}
    \frac{1}{\mu_{+}(N)^{M_{1}+\ldots+M_{l}}}\left(\frac{\mathcal{D}_{i_{l}}}{N}\right)^{M_{l}}\cdots \left(\frac{\mathcal{D}_{i_{1}}}{N}\right)^{M_{1}}G_{N}(\vec{x};\beta)\Bigg|_{\vec{x}=0}=\sum_{(E,\vec{q})\in \mathcal{B}(i_{1},\ldots,i_{l})}w(E,\vec{q}),
\end{align*}
where the right-hand side is a weighted sum over the configurations. In this section, we calculate the asymptotics of the right-hand side as $N\to\infty$, and then sum over all the tuples $(i_{1},\ldots,i_{l})$. Our limit expression is given as a functional of the conditional Brownian bridges, see Section \ref{sec:convergence}.

\subsection{Classifying local configurations}\label{sec:classify}

Fix $(i_{1},\ldots,i_{l})$. We will further classify the configurations in $\mathcal{B}(i_{1},\ldots,i_{l})$ according to their local behavior near each $C_{j}$, $j=1,2,\ldots,l-1$. 
The total weights of configurations with different local behaviors are of different orders when $N\rightarrow\infty$, and only a subset of them make a nontrivial contribution in the limit. In this section, we provide a list that completely classifies the possible local behaviors for our further analysis.
For $i\in\{i_{1},\ldots,i_{l}\}$, let \begin{equation}\label{eq_aibi}
a_{i}=\min\{j=1,2,\ldots,l:i_{j}=i\},\quad b_{i}=\max\{j=1,2,\ldots,l:i_{j}=i\}.\end{equation}
Let $\tau_{i}$ be a stopping time such that
    \begin{equation}\label{eq_taui}
        \tau_{i}=\min\{t: E(C_{a_{i}-1}+t)=E(C_{a_{i}-1})-1, 0<t\le\left\lfloor\frac{1}{4}(C_{a_{i}}-C_{a_{i}-1})\right\rfloor\},
    \end{equation}
    i.e., the waiting time after $C_{a_{i}-1}$ until the walk goes below $E(C_{a_{i}-1})$. If the  set is empty, we set $\tau_{i}=\varnothing$.
Let \begin{equation}\label{eq_tildetaui}
\tilde{\tau}_{i}=\min_{t>0}\{t:q(C_{a_{i}-1}+t)\ne q(C_{a_{i}-1}+t-)\},\end{equation}
i.e., the first jump time after $C_{a_{i}-1}$ at which a jump occurs. And we set $\tilde{\tau}_{i}=\varnothing$ if there is no jump after time $C_{a_{i}-1}$. These two quantities characterize the  behavior of the configuration locally on the right of $C_{a_{i}-1}$, or, roughly speaking, how $x_{i}$ enters the configuration.

More precisely, for each main variable $x_{i}$, we classify $x_{i}$ in the following types, according to the local behavior on the right of $C_{a_{i}-1}$:
\begin{itemize}
    \item Type I: $\tau_{i}\ne \varnothing$, and\footnote{This implies that $E(C_{a_{i}-1})>q(C_{a_{i}-1})$.} $\tau_{i}<\tilde{\tau}_{i}$ or $\tilde{\tau}_{i}=\varnothing$. 
    \item Type II: No jump of $x_{i}$ before time $C_{a_{i}-1}$, hence $E(C_{a_{i}-1})=q(C_{a_{i}-1})$. Moreover, $\tau_{i}=\tilde{\tau}_{i}<\infty$, and for some auxiliary variable $x_{j}$, $q_{j}(C_{a_{i}-1}+\tau_{i}-)>0$, $q_{j}(t)=0$ for all $t\ge C_{a_{i}-1}+\tau_{i}$.
    \item Type III: No jump of $x_{i}$ before time $C_{a_{i}-1}$, hence $E(C_{a_{i}-1})=q(C_{a_{i}-1})$. Moreover, $\tau_{i}=\tilde{\tau}_{i}<\infty$, and for some main variable $x_{j}$, $q_{j}(C_{a_{i}-1}+\tau_{i}-)>0$, $q_{j}(t)=0$ for all\footnote{In Type II and III, $x_{j}$ has its last jump at time $C_{a_{i}-1}+\tau_{i}$.} $t\ge C_{a_{i}-1}+\tau_{i}$.
    \item Type IV: $\tilde{\tau}_{i}<\tau_{i}$, or $\tilde{\tau}_{i}\ne \varnothing$ and $\tau_{i}=\varnothing$.
    \item Type V: $\tau_{i}=\tilde{\tau}_{i}=\varnothing$.

    Types I to V give a disjoint collection of local behavior on the right of $C_{a_{i}-1}$, while there are still some uncovered situations, which we collect in Type VI. More precisely, $x_{i}$ is in Type VI if $\tau_{i}=\tilde{\tau}_{i}\ne \varnothing$, and one of the following happens
    \item Type VI.1: No jump of $x_{i}$ before time $C_{a_{i}-1}$, hence $E(C_{a_{i}-1})=q(C_{a_{i}-1})$. Moreover, $\tau_{i}=\tilde{\tau}_{i}<\infty$, and for some main or auxiliary variable $x_{j}$, $q_{j}(C_{a_{i}-1}+\tau_{i}-)>q_{j}(C_{a_{i}-1}+\tau_{i})>0$. In other words, there is a jump of $x_{j}$ at time $C_{a_{i}-1}+\tau_{i}$, which does \textbf{not} send $q_{j}$ to 0.
    \item Type VI.2: $E(C_{a_{i}-1})>q(C_{a_{i}-1})$.
    \item Type VI.3: $E(C_{a_{i}-1})=q(C_{a_{i}-1})$, while there were jumps of $x_{i}$ before time $C_{a_{i}-1}$.
\end{itemize}
See Figure \ref{fig:types1} for an illustration of the different cases.

\begin{figure}
    \[
    \begin{tabular}{cc}
    \resizebox{0.35\textwidth}{!}{
    \begin{tikzpicture}[shorten >=-3pt,shorten <=-3pt,scale=0.8]
                \draw[fill=none,draw=none] (-2,-3) rectangle (18,13);

        \draw[dashed] (0,7)--(12,7);
        \draw[dashed] (0,0)--(12,0);
        \draw[dashed] (0,0)--(0,12);
        \draw[dashed] (7,0)--(7,12);
        \draw[dashed] (10,0)--(10,12);
        \draw[ultra thick,-o] (0,7)--(1,7);
        \draw[ultra thick,*-o] (1,9)--(2,9);
        \draw[ultra thick,*-o] (2,8)--(3,8);
        \draw[ultra thick,*-o] (3,10)--(4,10);
        \draw[ultra thick,*-o] (4,9)--(5,9);
        \draw[ultra thick,*-o] (5,8)--(6,8);
        \draw[ultra thick,*-o] (6,7)--(7,7);
        \draw[ultra thick,*-o] (7,6)--(8,6);
        \draw[ultra thick,*-o] (8,5)--(9,5);
        \draw[ultra thick,*-o] (9,9)--(10,9);
        \draw[ultra thick,*-o] (10,8)--(11,8);
        \draw[ultra thick,*-] (11,9)--(12,9);

        \draw[ultra thick,-o,blue] (0,2)--(10,2);
        \draw[ultra thick,*-,blue] (10,4)--(12,4);
        
        \node[below,font=\Huge,scale=1.5] at (0,-0.1) {$C_{a_i-1}$};
        \node[below,font=\Huge,scale=1.5] at (7,-0.1) {$\tau_i$};
        \node[below,font=\Huge,scale=1.5] at (10,-0.1) {$\tilde\tau_i$};

        \node[font=\Huge,scale=1.5] at (2,12) {I};
    \end{tikzpicture}
    } 
    &
    \resizebox{0.35\textwidth}{!}{
    \begin{tikzpicture}[shorten >=-3pt,shorten <=-3pt,scale=0.8]
                \draw[fill=none,draw=none] (-2,-3) rectangle (18,13);

        \draw[dashed] (0,7)--(12,7);
        \draw[dashed] (0,0)--(12,0);
        \draw[dashed] (0,0)--(0,12);
        \draw[dashed] (7,0)--(7,12);
        \draw[ultra thick,-o] (0,7)--(1,7);
        \draw[ultra thick,*-o] (1,9)--(2,9);
        \draw[ultra thick,*-o] (2,8)--(3,8);
        \draw[ultra thick,*-o] (3,10)--(4,10);
        \draw[ultra thick,*-o] (4,9)--(5,9);
        \draw[ultra thick,*-o] (5,8)--(6,8);
        \draw[ultra thick,*-o] (6,7)--(7,7);
        \draw[ultra thick,*-o] (7,6)--(8,6);
        \draw[ultra thick,*-o] (8,5)--(9,5);
        \draw[ultra thick,*-o] (9,9)--(10,9);
        \draw[ultra thick,*-o] (10,8)--(11,8);
        \draw[ultra thick,*-] (11,9)--(12,9);

        \draw[ultra thick,-o,blue] (0,6.9)--(7,6.9);
        \draw[ultra thick,*-,blue] (7,4.9)--(12,4.9);
        \draw[ultra thick,-o,red] (0,2)--(7,2);
        \draw[ultra thick,*-,red] (7,0)--(12,0);

        \node[below,font=\Huge,scale=1.5] at (0,-0.1) {$C_{a_i-1}$};
        \node[below,font=\Huge,scale=1.5] at (7,-0.1) {$\tau_i=\tilde\tau_i$};

        \node[font=\Huge,scale=1.5] at (2,12) {II$+$III};
    \end{tikzpicture}
    } 
    \end{tabular}
    \]
    \[
    \begin{tabular}{cc}
    \resizebox{0.35\textwidth}{!}{
    \begin{tikzpicture}[shorten >=-3pt,shorten <=-3pt,scale=0.8]
                \draw[fill=none,draw=none] (-2,-3) rectangle (18,13);

        \draw[dashed] (0,7)--(12,7);
        \draw[dashed] (0,0)--(12,0);
        \draw[dashed] (0,0)--(0,12);
        \draw[dashed] (7,0)--(7,12);
        \draw[dashed] (4,0)--(4,12);
        \draw[ultra thick,-o] (0,7)--(1,7);
        \draw[ultra thick,*-o] (1,9)--(2,9);
        \draw[ultra thick,*-o] (2,8)--(3,8);
        \draw[ultra thick,*-o] (3,10)--(4,10);
        \draw[ultra thick,*-o] (4,9)--(5,9);
        \draw[ultra thick,*-o] (5,8)--(6,8);
        \draw[ultra thick,*-o] (6,7)--(7,7);
        \draw[ultra thick,*-o] (7,6)--(8,6);
        \draw[ultra thick,*-o] (8,5)--(9,5);
        \draw[ultra thick,*-o] (9,9)--(10,9);
        \draw[ultra thick,*-o] (10,8)--(11,8);
        \draw[ultra thick,*-] (11,9)--(12,9);

        \draw[ultra thick,-o,blue] (0,2)--(4,2);
        \draw[ultra thick,*-,blue] (4,4)--(12,4);
        
        \node[below,font=\Huge,scale=1.5] at (0,-0.1) {$C_{a_i-1}$};
        \node[below,font=\Huge,scale=1.5] at (7,-0.1) {$\tau_i$};
        \node[below,font=\Huge,scale=1.5] at (4,-0.1) {$\tilde\tau_i$};

        \node[font=\Huge,scale=1.5] at (2,12) {IV};
    \end{tikzpicture}
    } 
    &
    \resizebox{0.35\textwidth}{!}{
    \begin{tikzpicture}[shorten >=-3pt,shorten <=-3pt,scale=0.8]
            \draw[fill=none,draw=none] (-2,-3) rectangle (18,13);

        \draw[dashed] (0,7)--(12,7);
        \draw[dashed] (0,0)--(12,0);
        \draw[dashed] (0,0)--(0,12);
        \draw[dashed] (12,0)--(12,12);
        \draw[ultra thick,-o] (0,7)--(1,7);
        \draw[ultra thick,*-o] (1,9)--(2,9);
        \draw[ultra thick,*-o] (2,8)--(3,8);
        \draw[ultra thick,*-o] (3,10)--(4,10);
        \draw[ultra thick,*-o] (4,9)--(5,9);
        \draw[ultra thick,*-o] (5,8)--(6,8);
        \draw[ultra thick,*-o] (6,7)--(7,7);
        \draw[ultra thick,*-o] (7,9)--(8,9);
        \draw[ultra thick,*-o] (8,8)--(9,8);
        \draw[ultra thick,*-o] (9,9)--(10,9);
        \draw[ultra thick,*-o] (10,8)--(11,8);
        \draw[ultra thick,*-] (11,11)--(12,11);

        \draw[ultra thick,-,blue] (0,3)--(12,3);
        
        \node[below,font=\Huge,scale=1.5] at (0,-0.1) {$C_{a_i-1}$};
        \node[below,font=\Huge,scale=1.5] at (12,-0.1) {$\frac{1}{4}(C_{a_i}-C_{a_i-1})$};

        \node[font=\Huge,scale=1.5] at (2,12) {V};
    \end{tikzpicture}
    } 
    \end{tabular}
    \]
    \[
    \begin{tabular}{ccc}
    \resizebox{0.32\textwidth}{!}{
    \begin{tikzpicture}[shorten >=-3pt,shorten <=-3pt,scale=0.8]
            \draw[fill=none,draw=none] (-2,-3) rectangle (18,13);

        \draw[dashed] (0,7)--(12,7);
        \draw[dashed] (0,0)--(12,0);
        \draw[dashed] (0,0)--(0,12);
        \draw[dashed] (7,0)--(7,12);
        \draw[ultra thick,-o] (0,7)--(1,7);
        \draw[ultra thick,*-o] (1,9)--(2,9);
        \draw[ultra thick,*-o] (2,8)--(3,8);
        \draw[ultra thick,*-o] (3,10)--(4,10);
        \draw[ultra thick,*-o] (4,9)--(5,9);
        \draw[ultra thick,*-o] (5,8)--(6,8);
        \draw[ultra thick,*-o] (6,7)--(7,7);
        \draw[ultra thick,*-o] (7,6)--(8,6);
        \draw[ultra thick,*-o] (8,5)--(9,5);
        \draw[ultra thick,*-o] (9,9)--(10,9);
        \draw[ultra thick,*-o] (10,8)--(11,8);
        \draw[ultra thick,*-] (11,9)--(12,9);

        \draw[ultra thick,-o,blue] (0,6.9)--(7,6.9);
        \draw[ultra thick,*-,blue] (7,4.9)--(12,4.9);
        \draw[ultra thick,-o,red] (0,3)--(7,3);
        \draw[ultra thick,*-,red] (7,1)--(12,1);

        \node[below,font=\Huge,scale=1.5] at (0,-0.1) {$C_{a_i-1}$};
        \node[below,font=\Huge,scale=1.5] at (7,-0.1) {$\tau_i=\tilde\tau_i$};

        \node[font=\Huge,scale=1.5] at (2,12) {VI.1};
    \end{tikzpicture}
    }  
    &
    \resizebox{0.32\textwidth}{!}{
    \begin{tikzpicture}[shorten >=-3pt,shorten <=-3pt,scale=0.8]
        \draw[fill=none,draw=none] (-2,-3) rectangle (18,13);

        \draw[dashed] (0,7)--(12,7);
        \draw[dashed] (0,0)--(12,0);
        \draw[dashed] (0,0)--(0,12);
        \draw[dashed] (7,0)--(7,12);
        \draw[ultra thick,-o] (0,7)--(1,7);
        \draw[ultra thick,*-o] (1,9)--(2,9);
        \draw[ultra thick,*-o] (2,8)--(3,8);
        \draw[ultra thick,*-o] (3,10)--(4,10);
        \draw[ultra thick,*-o] (4,9)--(5,9);
        \draw[ultra thick,*-o] (5,8)--(6,8);
        \draw[ultra thick,*-o] (6,7)--(7,7);
        \draw[ultra thick,*-o] (7,6)--(8,6);
        \draw[ultra thick,*-o] (8,5)--(9,5);
        \draw[ultra thick,*-o] (9,9)--(10,9);
        \draw[ultra thick,*-o] (10,8)--(11,8);
        \draw[ultra thick,*-] (11,9)--(12,9);

        \draw[ultra thick,-o,blue] (0,5)--(7,5);
        \draw[ultra thick,*-,blue] (7,2)--(12,2);

        \node[below,font=\Huge,scale=1.5] at (0,-0.1) {$C_{a_i-1}$};
        \node[below,font=\Huge,scale=1.5] at (7,-0.1) {$\tau_i=\tilde\tau_i$};

        \node[font=\Huge,scale=1.5] at (2,12) {VI.2};
    \end{tikzpicture}
    } 
    &
    \resizebox{0.32\textwidth}{!}{
    \begin{tikzpicture}[shorten >=-3pt,shorten <=-3pt,scale=0.8]
        \draw[fill=none,draw=none] (-5,-3) rectangle (15,13);
        \draw[dashed] (-4,7)--(10,7);
        \draw[dashed] (-4,0)--(10,0);
        \draw[dashed] (0,0)--(0,12);
        \draw[dashed] (7,0)--(7,12);
        \draw[dashed] (-3,0)--(-3,12);

        \draw[ultra thick,-o] (-4,9)--(-3,9);
        \draw[ultra thick,*-o] (-3,8)--(-2,8);
        \draw[ultra thick,*-o] (-2,9)--(-1,9);
        \draw[ultra thick,*-o] (-1,8)--(0,8);
        
        \draw[ultra thick,*-o] (0,7)--(1,7);
        \draw[ultra thick,*-o] (1,9)--(2,9);
        \draw[ultra thick,*-o] (2,8)--(3,8);
        \draw[ultra thick,*-o] (3,10)--(4,10);
        \draw[ultra thick,*-o] (4,9)--(5,9);
        \draw[ultra thick,*-o] (5,8)--(6,8);
        \draw[ultra thick,*-o] (6,7)--(7,7);
        \draw[ultra thick,*-o] (7,6)--(8,6);
        \draw[ultra thick,*-o] (8,5)--(9,5);
        \draw[ultra thick,*-] (9,9)--(10,9);

        \draw[ultra thick,-o,blue] (-4,7.9)--(-3,7.9);
        \draw[ultra thick,*-o,blue] (-3,6.9)--(7,6.9);
        \draw[ultra thick,*-,blue] (7,4.9)--(10,4.9);
        \draw[ultra thick,-o,red] (-4,1)--(-3,1);
        \draw[ultra thick,*-,red] (-3,0)--(10,0);

        \node[below,font=\Huge,scale=1.5] at (0,-0.1) {$C_{a_i-1}$};
        \node[below,font=\Huge,scale=1.5] at (7,-0.1) {$\tau_i=\tilde\tau_i$};

        \node[font=\Huge,scale=1.5] at (2,12) {VI.3};
    \end{tikzpicture}
    } 
    \end{tabular}
    \]
    \caption{An illustration of the classification of main variables into types I-VI.  In each figure the walk $E(t)$ is drawn in black, the sum of blocks $q(t)$ is drawn in blue, and a particular block $q_j$ is drawn in red if relevant. 
    \textbf{Top:} The left shows Type I in which $\tau_i <\tilde\tau_i$. The right shows Types II and III in which $\tau_i=\tilde\tau_i$. The difference between the two types depends on whether the block $q_j$, drawn in red, belongs to an auxiliary variable (Type II) or a main variable (Type III).
    \textbf{Middle:} The left shows Type IV in which $\tilde\tau_i<\tau_i$. The right shows Type V in which both $\tau_i$ and $\tilde \tau_i$ are $\varnothing$. Note that the illustrations show that $E(C_{a_i-1}) > q(C_{a_i-1})$ but $E(C_{a_i-1}) = q(C_{a_i-1})$ is also allowed for these types.
    \textbf{Bottom:} The illustrations for the three subtypes of Type VI. In each of the subtypes $\tau_i=\tilde\tau_i$. In Type VI.1, the block $q_j$ drawn in red can be any main or auxiliary variable. In Type VI.3 the red block $q_i$ represents variable $x_i$.}
    \label{fig:types1}
\end{figure}

We will introduce a cancellation between Type I and Type II in Section \ref{sec:cancellation}, to resolve the blow-up issue related to the weight given by these two types. Types III, IV, and V give nontrivial contributions to the limit, but to simplify the limiting expression, we perform another cancellation of terms of Type III in Section \ref{sec:secondcancellation}. One can think of each of Type VI.1-3 as certain ``atypical" deformations of Types I to V, and we show in Section \ref{sec:convergence} that the contribution of all of these vanishes in the limit. 

For $i\in\{i_{1},\ldots,i_{l}\}$, we also give another classification about how $x_{i}$ is removed during the action. More precisely, we say that $x_{i}$ is of Type A, B, or C, according to the following rules:
\begin{itemize}
    \item Type A: $q_{i}(t)=0$ for all $t\ge C_{b_{i}}$. 
    \item Type B: for some main variable $x_{j}$ of Type III, and $a_{j}-1\ge b_{i}$, $q_{i}(C_{a_{j}}+\tau_{j}-)>0$, $q_{i}(t)=0$ for all $t\ge C_{a_{j}}+\tau_{j}$.
    \item Type C: all the remaining cases.
\end{itemize}
See Figure \ref{fig:types2} for an illustration. 

We will show that Type A gives nontrivial contributions to the limit in Section \ref{sec:convergence}, while the total weight of configurations with some main variable of Type C vanishes as $N\rightarrow\infty$. Each $x_{i}$ of Type B is paired with some $x_{j}$ of Type III, so the weight of a configuration with any main variable of Type B is canceled out in Section \ref{sec:secondcancellation}.

\begin{figure}
    \centering
    \[
    \begin{tabular}{cc}
       \resizebox{0.4\textwidth}{!}{
    \begin{tikzpicture}[shorten >=-3pt,shorten <=-3pt,scale=0.8]
        \draw[fill=none,draw=none] (-5,-3) rectangle (15,13);
        \draw[dashed] (-4,7)--(10,7);
        \draw[dashed] (-4,0)--(10,0);
        \draw[dashed] (0,0)--(0,12);

        \draw[ultra thick,-o] (-4,9)--(-3,9);
        \draw[ultra thick,*-o] (-3,8)--(-2,8);
        \draw[ultra thick,*-o] (-2,9)--(-1,9);
        \draw[ultra thick,*-o] (-1,8)--(0,8);
        
        \draw[ultra thick,*-o] (0,7)--(1,7);
        \draw[ultra thick,*-o] (1,9)--(2,9);
        \draw[ultra thick,*-o] (2,8)--(3,8);
        \draw[ultra thick,*-o] (3,10)--(4,10);
        \draw[ultra thick,*-o] (4,9)--(5,9);
        \draw[ultra thick,*-o] (5,8)--(6,8);
        \draw[ultra thick,*-o] (6,7)--(7,7);
        \draw[ultra thick,*-o] (7,6)--(8,6);
        \draw[ultra thick,*-o] (8,5)--(9,5);
        \draw[ultra thick,*-] (9,9)--(10,9);

        \draw[ultra thick,-o,blue] (-4,7.9)--(-3,7.9);
        \draw[ultra thick,*-o,blue] (-3,6.9)--(0,6.9);
        \draw[ultra thick,*-o,blue] (0,4.9)--(8,4.9);
        \draw[ultra thick,*-,blue] (8,2)--(10,2);

        \draw[ultra thick,-,red] (-4,0)--(10,0);

        \node[below,font=\Huge,scale=1.5] at (0,-0.1) {$C_{b_i}$};

        \node[font=\Huge,scale=1.5] at (2.4,12) {Type A};
    \end{tikzpicture}
    } 
    & 
       \resizebox{0.4\textwidth}{!}{
        \begin{tikzpicture}[shorten >=-3pt,shorten <=-3pt,scale=0.8]
                \draw[fill=none,draw=none] (-2,-3) rectangle (18,13);

        \draw[dashed] (0,7)--(12,7);
        \draw[dashed] (0,0)--(12,0);
        \draw[dashed] (0,0)--(0,12);
        \draw[dashed] (7,0)--(7,12);
        \draw[ultra thick,-o] (0,7)--(1,7);
        \draw[ultra thick,*-o] (1,9)--(2,9);
        \draw[ultra thick,*-o] (2,8)--(3,8);
        \draw[ultra thick,*-o] (3,10)--(4,10);
        \draw[ultra thick,*-o] (4,9)--(5,9);
        \draw[ultra thick,*-o] (5,8)--(6,8);
        \draw[ultra thick,*-o] (6,7)--(7,7);
        \draw[ultra thick,*-o] (7,6)--(8,6);
        \draw[ultra thick,*-o] (8,5)--(9,5);
        \draw[ultra thick,*-o] (9,9)--(10,9);
        \draw[ultra thick,*-o] (10,8)--(11,8);
        \draw[ultra thick,*-] (11,9)--(12,9);

        \draw[ultra thick,-o,blue] (0,6.9)--(7,6.9);
        \draw[ultra thick,*-,blue] (7,4.9)--(12,4.9);
        \draw[ultra thick,-o,red] (0,2)--(7,2);
        \draw[ultra thick,*-,red] (7,0)--(12,0);

        \node[below,font=\Huge,scale=1.5] at (0,-0.1) {$C_{a_j-1}$};
        \node[below,font=\Huge,scale=1.5] at (7,-0.1) {$\tau_j=\tilde\tau_j$};

        \node[font=\Huge,scale=1.5] at (2.4,12) {Type B};
    \end{tikzpicture}
    }
    \end{tabular}
    \]
    \caption{An illustration of Types A and B. In both cases the block $q_i$ is drawn in red.}
    \label{fig:types2}
\end{figure}

Let $(j_{1},\ldots,j_{m})$ be the $m$-tuple that gives the indices in $\{i_{1}\}\cup\ldots\cup\{i_{l}\}$ without repetition and ordered according to their first appearance in $(i_{1},\ldots,i_{l})$. Let $J=\{j_{1},\ldots,j_{m}\}$ be the set of main variable indices, and let $J_{I},\ldots,J_{VI.3},\ J_{A},\ldots,J_{C}$ be the set of main variable indices of the corresponding type. We call the tuple \begin{equation}\label{eq_index}I(E,\vec{q})=(I_{1},I_{2})=(I_{1,1},\ldots,I_{1,m},I_{2,1},\ldots,I_{2,m})\in \{I,II,\ldots,VI\}^{m}\times \{A,B,C\}^{m}\end{equation} the \emph{index} of $(E,\vec{q})$, where each $I_{1,p},I_{2,p}$ gives the type of $x_{j_{p}}$.

As a preparation, we define the following notions for the jumps in a given configuration. Let $U$ be the set of auxiliary variable indices. For $i\in U$ or $i\in J_{B}\cup J_{C}$, let 
\begin{equation}\label{eq_lastjump1}l_{i}=\max\{t:q_{i}(t-)>0\}.\end{equation}
For $i\in J_{VI.3}$, let 
\begin{equation}\label{eq_lastjump2}l'_{i}=\max\{t\le C_{a_{i}-1}: q_{i}(t)>0\}.\end{equation} 
We call the jumps at time $l_{i}$ and $l'_{i}$ the \emph{last jumps} of $x_{i}$, and denote the set of last jump times of $x_{i}s,\ i\in J\cup U$ by $\mathcal{L}$. By definition, $|\mathcal{L}|=|U|+|J_{B}|+|J_{C}|+|J_{VI.3}|$.

Similarly,  let 
$$f_{i}=\begin{cases}
\min\{t:\ q_{i}(t)>0 \},\quad & i\in U;\\
\min \{t\le C_{a_{i}-1}: q_{i}(t)>0 \},\quad & i\in J.
\end{cases}$$
We call $f_{i}$ a \emph{first jump} of  $x_{i}$, and denote the set of first jump times of the variables $x_{i}$, $i\in J\cup U$ by $\mathcal{F}$. Note that for $i\in J$, $x_{i}$ might not have a first jump. Denote the total number of first jumps of $x_{i}s,\ i\in J$ by $f_{J}$, so that by definition $|\mathcal{F}|=|U|+f_{J}$. We say a jump  is a \emph{free jump} if it is not a first or last jump, and denote the set of free jump times by $\mathcal{R}$.

Recall that $\Delta$ denotes the set of all jump times. If $x_{i}$ is of Type A, let 
\begin{equation}\label{etai}
    \eta_{i}=\max\{C_{b_{i}}-\left\lfloor\frac{1}{4}(C_{b_{i}}-C_{b_{i}-1})\right\rfloor\le t< C_{b_{i}}: t\in \Delta\setminus\mathcal{L}\}.
\end{equation}
If the set is empty, we set $\eta_{i}=\varnothing$.

The classifications we give in this section are the same as those in the concurrent work \cite[Section 4.4]{GXZ}. While in this text we deal with more complicated walk-bridges, the only negative increments of the walk are still $-1$, and their scaling limits are still Brownian bridges. This allows the two works to share some common techniques.

 \subsection{Blowing up terms and cancellations}\label{sec:cancellation}

In this section, we introduce an equivalence relation between configurations created by the same Dunkl action as in Eqn.~(\ref{eq_keyaction}), but with different indices. 
The reason for doing this is that, if the configuration has at least one main variable of Type I or II, then the total weight of configurations $(E,\vec{q})$ with possibly different $E$ and a fixed index will blow up in the limit. This will be clearer by looking at the following example.

\begin{example}\label{ex:scalingissue1}
    Let us consider the case when $l=2$. We are interested in the action of
    \[
    \frac{N(N-1)}{\mu_+(N)^{2M}} \left(\frac{\mathcal{D}_{2}}{N}\right)^{M} \left(\frac{\mathcal{D}_{1}}{N}\right)^{M}
    \]
    on the BGF. First, let us consider the case in which there are no jumps. In this case, the index $I_{1}=(V,V)$, and the walk must hit zero at the $M$-th step as shown in the figure below with $M=14$.
    \[
    \resizebox{0.8\textwidth}{!}{
    \begin{tikzpicture}[shorten >=-3pt,shorten <=-3pt]
        \draw[ultra thin, opacity = 0.5] (0,0) grid (28,7);

        \draw[ultra thick,*-o] (0,0)--(1,0);
        \draw[ultra thick,*-o] (1,2)--(2,2);
        \draw[ultra thick,*-o] (2,1)--(3,1);
        \draw[ultra thick,*-o] (3,4)--(5,4);
        \draw[ultra thick,*-o] (5,5)--(6,5);
        \draw[ultra thick,*-o] (6,4)--(7,4);
        \draw[ultra thick,*-o] (7,3)--(9,3);
        \draw[ultra thick,*-o] (9,2)--(10,2);
        \draw[ultra thick,*-o] (10,1)--(11,1);
        \draw[ultra thick,*-o] (11,2)--(13,2);
        \draw[ultra thick,*-o] (13,1)--(14,1);
        \draw[ultra thick,*-o] (14,0)--(15,0);
        \draw[ultra thick,*-o] (15,3)--(16,3);
        \draw[ultra thick,*-o] (16,4)--(18,4);
        \draw[ultra thick,*-o] (18,3)--(19,3);
        \draw[ultra thick,*-o] (19,5)--(20,5);
        \draw[ultra thick,*-o] (20,4)--(22,4);
        \draw[ultra thick,*-o] (22,3)--(23,3);
        \draw[ultra thick,*-o] (23,2)--(24,2);
        \draw[ultra thick,*-o] (24,1)--(25,1);
        \draw[ultra thick,*-o] (25,0)--(26,0);
        \draw[ultra thick,*-o] (26,2)--(27,2);
        \draw[ultra thick,*-o] (27,1)--(28,1);
        \draw[ultra thick,*-] (28,0);
        
        \draw[ultra thick] (14,0)--(14,7);
    \end{tikzpicture}
    }
    \]
    Now consider the case in which we have exactly one jump. There are two possibilities. 
    The jump could add powers of $x_2$ in the section where $x_1$ is the main variable, in which case the index is (V,I), and the height of the $x_2$ block must be equal to the height of the walk at time $M$. Or we could have a jump removing leftover powers of $x_1$ in the section where $x_2$ is the main variable, in which case the index is (V,III). Examples of the two possibilities are given below with the block drawn in dashed blue.
    \[
    \begin{aligned}
    & 
    \resizebox{0.8\textwidth}{!}{
    \begin{tikzpicture}[shorten >=-3pt,shorten <=-3pt]
        \draw[ultra thin, opacity = 0.5] (0,0) grid (28,8);

        \draw[ultra thick,*-o] (0,0)--(1,0);
        \draw[ultra thick,*-o] (1,2)--(2,2);
        \draw[ultra thick,*-o] (2,1)--(3,1);
        \draw[ultra thick,*-o] (3,4)--(5,4);
        \draw[ultra thick,*-o] (5,5)--(6,5);
        \draw[ultra thick,*-o] (6,4)--(7,4);
        \draw[ultra thick,*-o] (7,3)--(8,3);
        \draw[ultra thick,*-o] (8,6)--(9,6);
        \draw[ultra thick,*-o] (9,7)--(11,7);
        \draw[ultra thick,*-o] (11,6)--(12,6);
        \draw[ultra thick,*-o] (12,7)--(13,7);
        \draw[ultra thick,*-o] (13,6)--(14,6);
        \draw[ultra thick,*-o] (14,5)--(15,5);
        \draw[ultra thick,*-o] (15,7)--(18,7);
        \draw[ultra thick,*-o] (18,6)--(19,6);
        \draw[ultra thick,*-o] (19,5)--(20,5);
        \draw[ultra thick,*-o] (20,4)--(22,4);
        \draw[ultra thick,*-o] (22,3)--(23,3);
        \draw[ultra thick,*-o] (23,2)--(24,2);
        \draw[ultra thick,*-o] (24,1)--(25,1);
        \draw[ultra thick,*-o] (25,0)--(26,0);
        \draw[ultra thick,*-o] (26,2)--(27,2);
        \draw[ultra thick,*-o] (27,1)--(28,1);
        \draw[ultra thick,*-] (28,0);
        
        \draw[ultra thick,blue,dashed,*-o] (0,0)--(11,0);
        \draw[ultra thick,blue,dashed,*-*] (11,5)--(14,5);
        
        \draw[ultra thick] (14,0)--(14,8);
    \end{tikzpicture}
    }
    \\
    & 
    \resizebox{0.8\textwidth}{!}{
    \begin{tikzpicture}[shorten >=-3pt,shorten <=-3pt]
        \draw[ultra thin, opacity = 0.5] (0,0) grid (28,8);

        \draw[ultra thick,*-o] (0,0)--(1,0);
        \draw[ultra thick,*-o] (1,2)--(2,2);
        \draw[ultra thick,*-o] (2,1)--(3,1);
        \draw[ultra thick,*-o] (3,4)--(5,4);
        \draw[ultra thick,*-o] (5,5)--(6,5);
        \draw[ultra thick,*-o] (6,4)--(7,4);
        \draw[ultra thick,*-o] (7,3)--(8,3);
        \draw[ultra thick,*-o] (8,6)--(9,6);
        \draw[ultra thick,*-o] (9,7)--(11,7);
        \draw[ultra thick,*-o] (11,6)--(12,6);
        \draw[ultra thick,*-o] (12,7)--(13,7);
        \draw[ultra thick,*-o] (13,6)--(14,6);
        \draw[ultra thick,*-o] (14,5)--(15,5);
        \draw[ultra thick,*-o] (15,7)--(18,7);
        \draw[ultra thick,*-o] (18,6)--(19,6);
        \draw[ultra thick,*-o] (19,5)--(20,5);
        \draw[ultra thick,*-o] (20,4)--(22,4);
        \draw[ultra thick,*-o] (22,3)--(23,3);
        \draw[ultra thick,*-o] (23,2)--(24,2);
        \draw[ultra thick,*-o] (24,1)--(25,1);
        \draw[ultra thick,*-o] (25,0)--(26,0);
        \draw[ultra thick,*-o] (26,2)--(27,2);
        \draw[ultra thick,*-o] (27,1)--(28,1);
        \draw[ultra thick,*-] (28,0);

        \draw[ultra thick,blue,dashed,*-o] (14,5)--(20,5);
        \draw[ultra thick,blue,dashed,*-*] (20,0)--(28,0);
        
        \draw[ultra thick] (14,0)--(14,8);
    \end{tikzpicture}
    }
    \end{aligned}
    \]
    One can check that as $M= O(N^{2/3}), N\rightarrow\infty$, the total weight of configurations with index $I_{1}=(V,I)$  is of order $O(\sqrt{M})$, and the total weight of configurations with $I_{1} =(V,III)$ is of order $O(1)$. The former weight blows up in the limit.
    \end{example}
    
To resolve this issue, we say that a configuration $$(E,\vec{q})\in \mathcal{I} \quad\text{if}\quad |J_{II}|=0.$$ A configuration $(E,\vec{q})$ is \emph{grouped with $(E,\vec{q}')\in \mathcal{I}$}, written as $(E,\vec{q})\sim(E,\vec{q}')$, if  one can obtain $(E,\vec{q})$ from $(E,\vec{q}')$ by the following procedure. 
For all the main variables $x_{i}$ with $i\in J_{II}$,  let $x_{j}$ denote the unique variable such that $q_{j}(C_{a_{i}-1}+\tau_{i}-)>0$, $q_{j}(t)=0$ for all $t\ge C_{a_{i}-1}+\tau_{i}$. Remove the last jump of $x_{j}$ at time $C_{a_{i}-1}+\tau_{i}$, and replace all of the remaining jumps of $x_{j}$ by jumps of $x_{i}$. Note that $i\in J_{II}$ in $(E,\vec{q}')$ implies that $i\in J_{I}$ in $(E,\vec{q})$.

From the above algorithm, one can see that each $(E,\vec{q})\in \mathcal{I}$ provides an equivalence class of configurations  $(E,\vec{q}')\sim (E,\vec{q})$, indexed by the  $2^{|J_{I}|}$ elements in $\{I,II\}^{|J_{I}|}$. Moreover, each $(E,\vec{q}')$ is grouped with a unique element in $\mathcal{I}$. If a main variable $x_{i}$ is of Type II in $(E,\vec{q}')$, the extra jump at time $C_{a_{i}-1}+\tau_{i}$ gives a minus sign in the weight, so the weights of the grouped configurations are of opposite sign. Because of this, the inner cancellation of weights in each equivalence class resolves the blow-up issue. We give a more detailed discussion in the next section.

\begin{example}\label{ex:scalingissue2}

   The configurations with $I_{1}=(V,I)$ in Example \ref{ex:scalingissue1} will be grouped with configurations that have two jumps, with  $I_{1}=(V,II)$, which looks as follows:
    \[
    \resizebox{0.8\textwidth}{!}{
    \begin{tikzpicture}[shorten >=-3pt,shorten <=-3pt]
        \draw[ultra thin, opacity = 0.5] (0,0) grid (28,8);

        \draw[ultra thick,*-o] (0,0)--(1,0);
        \draw[ultra thick,*-o] (1,2)--(2,2);
        \draw[ultra thick,*-o] (2,1)--(3,1);
        \draw[ultra thick,*-o] (3,4)--(5,4);
        \draw[ultra thick,*-o] (5,5)--(6,5);
        \draw[ultra thick,*-o] (6,4)--(7,4);
        \draw[ultra thick,*-o] (7,3)--(8,3);
        \draw[ultra thick,*-o] (8,6)--(9,6);
        \draw[ultra thick,*-o] (9,7)--(11,7);
        \draw[ultra thick,*-o] (11,6)--(12,6);
        \draw[ultra thick,*-o] (12,7)--(13,7);
        \draw[ultra thick,*-o] (13,6)--(14,6);
        \draw[ultra thick,*-o] (14,5)--(15,5);
        \draw[ultra thick,*-o] (15,7)--(18,7);
        \draw[ultra thick,*-o] (18,6)--(19,6);
        \draw[ultra thick,*-o] (19,5)--(20,5);
        \draw[ultra thick,*-o] (20,4)--(22,4);
        \draw[ultra thick,*-o] (22,3)--(23,3);
        \draw[ultra thick,*-o] (23,2)--(24,2);
        \draw[ultra thick,*-o] (24,1)--(25,1);
        \draw[ultra thick,*-o] (25,0)--(26,0);
        \draw[ultra thick,*-o] (26,2)--(27,2);
        \draw[ultra thick,*-o] (27,1)--(28,1);
        \draw[ultra thick,*-] (28,0);

        \draw[ultra thick, dashed, red,*-o] (0,0)--(11,0);
        \draw[ultra thick, dashed, red,*-o] (11,5)--(20,5);
        \draw[ultra thick, dashed, red,*-*] (20,0)--(28,0);
        
        \draw[ultra thick] (14,0)--(14,8);
    \end{tikzpicture}
    }
    \]
    Here the red dashed line indicates the block for some $x_i$, $i\ne 1,2$. Note that to the left of $C_1$ we have kept the same jump location and height (but have changed the index) and we have added a new jump to the right of $C_1$ at the unique location that can send the block to zero. In addition, the weights of these two types of configurations are of opposite signs, and their sum is of order $O(1)$ when $N\rightarrow\infty$.

    \end{example}

\subsection{Tail estimates}\label{sec:tailestimates}

In this section, we prove the main technical lemma of this text, which plays two roles. Firstly, it shows that the weight contribution of a configuration  decays exponentially, hence it bounds the tail of Eqn.~\eqref{eq_weight} asymptotically. Secondly, the exponent of $N$ shows that the contribution of configurations with Type VI local behavior is of smaller order in $N$, while configurations with Type I local behavior cause a blow-up issue.

For a configuration $(E,\vec{q})\in \mathcal{B}(i_{1},\ldots,i_{l})$, we refer to the following collection as the \emph{characters} of $(E,\vec q)$.:
\begin{itemize}
    \item $I(E,\vec{q})$, the index of $(E,\vec{q})$ defined in Eqn.~(\ref{eq_index}).
    \item $\tau_{j_{1}},\ldots,\tau_{j_{m}}\in \Z_{\ge 1}\cup \varnothing$.
    \item $\eta_{j_{p}}$ defined in Eqn.~\eqref{etai} for $j_{p}\in J_{A}$, and $J_{\varnothing,A}$, the set of $j_{p}$ where $\eta_{j_{p}}=\varnothing$.
\end{itemize}

\begin{lem}\label{lem:tail}
    There exists a constant $C>1$ that depends on $\bk_{1},\ldots,\bk_{l}$ and $\sigma$, and is uniform in $N$, such that for all $N\ge 100$,
    \begin{equation}\label{eq_tail}
    \begin{split}
        \sum_{(E,\vec{q})\in \mathcal{B}(i_{1},\ldots,i_{l})\cap\mathcal{I}}\left|w(E,\vec{q})\right|\le&C \prod_{j_{p}\in J_{A}:\eta_{j_{p}}\ne \varnothing}(C_{b_{j_{p}}}-\eta_{j_{p}}+1)^{-1/2}\cdot \prod_{p:j_{p}\in J_{IV}}\tau_{j_{p}}^{-1/2}\prod_{p:j_{p}\in J_{I}\cup J_{III}\cup J_{VI}}\tau_{j_{p}}^{-3/2}\\
       \times& N^{-m}N^{-\frac{1}{3}(|J_{A}|-|J_{\varnothing,A}|)}N^{-\frac{1}{3}|J_{IV}|}N^{\frac{1}{3}|J_{I}|-\frac{1}{3}|J_{VI}|-\frac{1}{3}|J_{VI.3}|-\frac{1}{3}|J_{C}|},
        \end{split}
    \end{equation}
    where the sum is over all configurations with given characters. 
\end{lem}

\begin{proof}
We divide the proof into the following three parts.
\medskip

\noindent{\textbf{Part I. Enumerations}}
  
 We enumerate the total weight of the configurations with the given characters and the following fixed data:\begin{itemize}
       \item  $\delta$, the total number of jumps,
       \item  $|U|$, the number of auxiliary variables, and
       \item  $k_{p}$, the number of jumps of $x_{j_{p}}$ for $p=1,2,\ldots,m$.
   \end{itemize}  Recall that the contribution of each walk/block configuration can be both positive or negative, and by counting without regard to the sign we get a summation of their absolute values. To do this we bound the count of all the data that must be determined in order to specify a walk/block configuration consistent with the above data. This is done in several steps.

  \noindent 0. Choose $|U|$ auxiliary variables from $x_{1},\ldots,x_{N}$. The number of choices is bounded by \begin{equation}\label{eq_taillemma0}
      \binom{N}{|U|}\le \frac{N^{|U|}}{|U|!}.
  \end{equation}
  
  \noindent 1. Determine all the first and free jump times. In particular, for the main variables we determine the following.
  \begin{itemize}
     \item When $j_{p}\in J_{VI}$, we determine whether there is a jump of $x_{i},\ i\in J$ at $C_{a_{j_{p}}-1}+\tau_{j_{p}}$ (if yes, then its jump time is fixed), and we denote the total number of such jumps by $|J_{m,VI}|$. 
     \item When $j_{p}\in J_{A}$, we determine whether there is a jump of $x_{i},\ i\in J$ at $\eta_{j_{p}}$  (if yes, then its jump time is fixed), and we denote the total number of such jumps by $|J_{m,A}|$.  
 \end{itemize}  
 Let $k'_{p}\le k_{p}$ be the number of first and free jump times of $x_{j_{p}}$ not fixed by the above items. 

 For each auxiliary variable $x_{i}$, $i\in U$, we determine the first jump time of $x_{i}$. Then we determine the free jump times of all the auxiliary variables $x_{i}$ without specifying their corresponding variables. In particular, for $p=1,2,\ldots,m$ we determine the following.
 \begin{itemize}
     \item When $j_{p}\in J_{VI}$, we determine whether there is a jump of $x_{i},\ i\in U$ at $C_{a_{j_{p}}-1}+\tau_{j_{p}}$ (if yes, then its jump time is fixed), and we denote the total number of such jumps by $|J_{a,VI}|$. By definition, $|J_{m,VI}|+|J_{a,VI}|=|J_{VI}|$.
     \item When $j_{p}\in J_{A}$, we determine whether there is a jump of $x_{i},\ i\in U$ at $\eta_{j_{p}}$  (if yes, then its jump time is fixed), and we denote the total number of such jumps by $|J_{a,A}|$. By definition, $|J_{m,A}|+|J_{a,A}|=|J_{A}|-|J_{\varnothing,A}|$.
 \end{itemize}  
 There are no more than $(l+1)^{2l}$ possible cases in the above four items, so we just consider an arbitrary case of them. Also note that for $j_{p}\in J_{IV}$, there is at least one first/free jump between $C_{a_{j_{p}}-1}$ and $C_{a_{j_{p}}-1}+\tau_{j_{p}}$. Therefore\footnote{There are also restrictions that \begin{itemize}
      \item There is no jump in $\llbracket C_{a_{j_{p}}-1}+1,C_{a_{j_{p}}-1}+\tau_{j_{p}}\rrbracket$ if $x_{j_{p}}$ is of Type I and III. There is no jump in $\llbracket C_{a_{j_{p}}-1}+1,C_{a_{j_{p}}-1}+\tilde{\tau}_{j_{p}}-1\rrbracket$ if $x_{j_{p}}$ is of Type VI. There is no jump in $\llbracket C_{a_{j_{p}}-1}+1, C_{a_{j_{p}}-1}+\left\lfloor\frac{1}{4}(C_{a_{j_{p}}}-C_{a_{j_{p}}-1})\right\rfloor\rrbracket$ if $x_{j_{p}}$ is of Type V.
     \item There is no jump in $\llbracket \eta_{j_{p}}+1,C_{b_{j_{p}}}\rrbracket$.
      \end{itemize} But we do not take these restrictions into account when giving the bound. }, 
  the total number of choices of the first and free jump times not fixed by the above is bounded by
      \begin{equation}\label{eq_taillemma2}
      \begin{split}&\prod_{p:\ j_{p}\in J_{IV}}\tau_{j_{p}}\cdot \binom{M}{\underbrace{1,\ldots,1}_{|U|},k_{1}',\ldots, k_{m}',\delta-\sum_{p=1}^{m}k'_{p}-2|U|-|J_{IV}|-|J_{VI}|-|J_{m,A}|-|J_{a,A}|-|J_{B}|-|J_{C}|}\\
      \le&\prod_{p:\ j_{p}\in J_{IV}}\tau_{j_{p}}\cdot\frac{M^{\delta-|J_{A}|+|J_{\varnothing,A}|-|J_{B}|-|J_{C}|-|U|-|J_{IV}|-|J_{VI}|}}{(k_{1}-l-1)!\cdots (k_{m}-l-1)!(\delta-\sum_{p=1}^{m}k_{p}-2|U|-2l)!},
      \end{split}
      \end{equation}
      where here and below we take $n!=1$ for $n\in \Z_{\le 0}$.
  
  \noindent 2. For $p=1,2,\ldots,m$, let $$X_{j_{p}}=\begin{cases}
      C_{a_{j_{p}}-1}+\tau_{j_{p}}\quad\quad&\text{if}\ \ \tau_{j_{p}}\ne \varnothing,\\
      C_{a_{j_{p}}-1}+\left\lfloor\frac{1}{4}(C_{a_{j_{p}}}-C_{a_{j_{p}}-1})\right\rfloor\quad &\text{if}\ \ \tau_{j_{p}}=\varnothing
  \end{cases}$$

  and let $$Y_{j_{p}}=\begin{cases}
      \eta_{j_{p}},&\text{if}\ \eta_{j_{p}}\ne \varnothing,\\
     C_{b_{j_{p}}} -\left\lfloor\frac{1}{4}(C_{b_{j_{p}}}-C_{b_{j_{p}}-1})\right\rfloor,\quad &\text{if}\ \eta_{j_{p}}=\varnothing.
  \end{cases}$$

 Split $\llbracket 0,M\rrbracket$ into disjoint segments by  $C_{0},\ldots,C_{l}$, the points  $X_{j_{p}}$ and the points $Y_{j_{p}}$.  Denote the segments by $S_{1},\ldots,S_{n}$, ordered backward in time. More precisely, all the above positions (except $C_{0}=0$) serve as the right boundary of the segment on its left. Note that $n\le 3l$. 
  
Split $\llbracket 0,M\rrbracket$ into finer disjoint intervals by $C_{0},\ldots,C_{l}$,  $X_{j_{p}}$, $Y_{j_{p}}$, all the first and free jump times of $x_{i}$, $i\in J\cup U$. Denote these intervals by $\tilde{S}_1,\ldots,\tilde{S}_{\tilde{n}}$, and the left endpoint of $\tilde{S}_{\tilde{k}}$ by $t_{\tilde{k}}$. More precisely, all the above positions (except $C_{0}=0$) serve as the right boundary of the interval on its left. Note that the intervals are subdivisions of the segments $S_1,\ldots,S_n$, and $$\delta-|U|-l\le \tilde{n}\le \delta-|U|+4l.$$

 Recall the definitions of last jump from Eqn.~(\ref{eq_lastjump1}) and (\ref{eq_lastjump2}).  We make two determinations about the last jump\footnote{The last jump times should occur after the corresponding first jumps, but we do not take this restriction into account in this proof.} times. First, we determine the number of last jumps in each segment $\tilde{S}_{k}$. Second, we determine the order of these last jump times. The total number of choices is bounded by 
  \begin{equation}\label{eq_taillemma3}
      \binom{\tilde{n}-1+|U|+2l}{|U|+2l}(|U|+2l)!\le \frac{(\delta+6l-1)!}{(\delta-|U|+4l-1)!}.
  \end{equation}
  
  \noindent 3.  For $k=1,2,\ldots,n$, starting from $k=1$ (the segment ending at $M$), sample a walk $E_{k}(t): \{E(t)\}_{t\in S_{k}}$ backward on $S_{k}$, which is compatible with all the ``jump time data" determined in the previous steps. We view each segment of $E(t)$ as a random walk  with fixed right endpoint, and use Lemma \ref{lem:walkmax} to bound the partition function $$Par(k):=\sum_{E_{k}(t)}\prod_{l\ge 0}P_{l}(N)^{\# t:\ E_{k}(t)-E_{k}(t-)=l-1}$$ of the walk bridge. More precisely, let $C>0$ be a large enough constant uniform in $N$, 
  \begin{itemize}

      \item  $E(t)\ge E(C_{C_{b_{j_{p}}}})$ on $[Y_{j_{p}},C_{b_{j_{p}}}]$ for all $x_{j_{p}}$ of Type A. Then by Lemma \ref{lem:walkmax} (1),
      \begin{equation}\label{eq_conditionalprob1}
         Par(k)\le C\cdot (C_{b_{j_{p}}}-Y_{j_{p}}+1)^{-1/2}\le \begin{cases}
             C\cdot N^{-1/3},\quad& j_{p}\in J_{\varnothing,A};\\
             C\cdot (C_{b_{j_{p}}}-Y_{j_{p}}+1)^{-1/2},\quad& j_{p}\in J\setminus J_{\varnothing,A}
         \end{cases}.
     \end{equation}

      \item  $E(t)\ge E(C_{a_{j_{p}}-1})$ on $[C_{a_{j_{p}}-1}, C_{a_{j_{p}}-1}+\left\lfloor\frac{1}{4}(C_{a_{j_{p}}}-C_{a_{j_{p}}-1})\right\rfloor]$ for all $x_{j_{p}}$ of Type V. Then by  Lemma \ref{lem:walkmax} (2),
      \begin{equation}\label{eq_conditionalprob3}
      \begin{split}
          Par(k)\le\begin{cases}
              C\cdot N^{-2/3},\quad&\text{if\ there\ is\ no\ first\ jump\ of\ } x_{j_{p}};\\
              C\cdot N^{-1/3},\quad&\text{otherwise}.
          \end{cases}
          \end{split}
      \end{equation}
      The first inequality holds\footnote{In particular, this is the case for $Par_{n}$ since $E(0)=0$.} since, in this case, the value of $E(C_{a_{j_{p}}-1})$ is fixed to be $q\left(C_{a_{j_{p}-1}}+\left\lfloor\frac{1}{4}(C_{a_{i}}-C_{a_{i}-1})\right\rfloor\right)$. Denote the number of $S_{k}$ such that the first case above holds by $\delta_{1}$.

      \item  $E(C_{a_{j_{p}}-1})=E( C_{a_{j_{p}}-1}+\tau_{j_{p}}-1)=E( C_{a_{j_{p}}-1}+\tau_{j_{p}})+1$ when $\tau_{i}\ne\varnothing$. Then by Lemma  \ref{lem:walkmax} (3),
      \begin{equation}\label{eq_conditionalprob2}
          Par(k)\le C\cdot \tau_{j_{p}}^{-3/2}.
      \end{equation}
      \item $E(t)=E(t-1)-1$ at each jump time $t$. This further bounds the partition functions by a factor less than 1, which we ignore in this proof.
      
  \end{itemize}.

  \noindent 4.   Moving backward in time, determine  the height of all last jumps, and the height of all first and free jumps of auxiliary variables, in each segment $\tilde{S}_{\tilde{k}}\subset S_{k}$.

  \noindent 5. Repeat Steps 3 and 4 in segments $S_{k+1}$ until $k=n$. Note that the walks on each segment are mutually independent, conditioned on the boundary heights.
  
Now we bound the total number of choices of jump heights, given the walk $E(t)$.  On each $\tilde{S}_{\tilde{k}}$, denote the number of last jumps by $u(\tilde{k})$. Then the number of choices of last jump heights on each $\tilde{S}_{\tilde{k}}$ is bounded by $$\frac{[E(t_{\tilde{k}})-\min_{t\in \tilde{S}_{\tilde{k}}} E(t)]^{u'(\tilde{k})}}{u'(\tilde{k})!},$$
where $u'(\tilde{k})=u(\tilde{k})$ or $u(\tilde{k})-1$ is the ``degree of freedom" of the last jump heights. $u'(\tilde{k})=u(\tilde{k})-1$ when one of the last jump heights is fixed by the restriction $E(C_{a_{j_{p}}-1})=q(C_{a_{j_{p}}-1})$ for some $j_{p}\in J_{III}\cup J_{VI.3}$.

As for the choices of each first or free jump height and its corresponding variable, there are three cases: \begin{itemize}
    \item no freedom if it is a first jump;
    \item if it is not a first jump, its number of choices is bounded by $\max_{0\le t\le M} E(t)$; and
    \item in some of the configurations, one of the free jump heights is fixed by the restriction $E(C_{a_{j_{p}}-1})=q(C_{a_{j_{p}}-1})$ for some $j_{p}\in J_{III}\cup J_{IV}\cup J_{V}\cup J_{VI.1}\cup J_{VI.3}$. Then there are at most $|U|$ choices of this jump, given by the $|U|$ choices of its jump variable.
\end{itemize}

Let $\delta_{2}$ be the degree of freedom loss of the free jumps, $\delta_{3}=\sum_{\tilde{k}=1}^{\tilde{n}}u(\tilde{k})-u'(\tilde{k})$ be the total degree of freedom loss of the last jumps. Note that by definition $\delta_{2}\le m\le l$, and \begin{equation}\label{eq_equality}
    \delta_{1}+\delta_{2}+\delta_{3}+f_{J}=m+|J_{VI.3}|.
\end{equation} Then we obtain the factor:

\begin{equation}\label{eq_step8factor}
    |U|^{\delta_{2}}[\max_{0\le t\le M}E(t)]^{|\mathcal{R}|-\delta_{2}} \prod_{\tilde{k}=1}^{\tilde{n}}\frac{[ E(t_{\tilde{k}})-\min_{t\in S_{\tilde{k}}} E(t)]^{u'(\tilde{k})}}{u'(\tilde{k})!}.
\end{equation}
  
\vspace{0.5cm}

\noindent{\textbf{Part II. Taking expectations}}
Recall the weight of each configuration from Eqn.~(\ref{eq_weight}). First note that for a walk $E(t)$ such that $$N^{-\frac{1}{3}}\max_{0\le t\le M}E(t)\in [h,h+1),$$ we have \begin{equation}\label{eq_taillemma1}
   N^{-\delta}\prod_{t:\  t+1\in D(E)}\left(\frac{m_{i}(t)}{N}+2\frac{E(t)-q(t)}{\beta N}\right)\le N^{-\delta}\left(1+\frac{\max_{0\le t\le M}E(t)}{N}\right)^{M}\le N^{-\delta} \exp(Ch).
\end{equation} 

Now we take the weighted sum over $E(t)$ as an expectation on the conditional walk bridges described above. Denote $\max_{t\in S_{k}}E(t)=h_{k}N^{\frac{1}{3}}$, then by Lemma \ref{lem:walkmax}, among any possible conditioning at $C_{i}$, there exists some uniform constant $C>0,\ c>0$, such that $\PP[h\le h_{k}<h+1]\le C\exp(-ch^{2})$. To deal with the second part of Eqn.~(\ref{eq_step8factor}), we apply\footnote{The purpose of introducing this lemma is to get an upper bound uniformly in $\delta$ and $\tilde{n}$, because the bound only depends on $n$, the number of segments $S_{k}$, which is still uniformly bounded. }  Lemma \ref{lem:max-min}. Denote $u'(k)=\sum_{\tilde{k}:\ \tilde{S}_{\tilde{k}}\subset S_{k}}u'(\tilde{k})$ the degree of freedom of the last jump heights in $S_{k}$. We combine the terms corresponding to segments $\tilde{S}_{\tilde{k}}$ in the same segment $S_{k}$, and give an upper bound of Eqn.~(\ref{eq_step8factor}) 
\begin{equation}\label{eq_taillemma5}
    |U|^{\delta_{2}}(h N^{1/3})^{|\mathcal{R}|-\delta_{2}}\prod_{k=1}^{n}C^{u'(k)}\log(u'(k)+2)^{u'(k)+2}u'(k)^{-u'(k)/2}(N^{1/3})^{u'(k)}\cdot C\exp(-ch^{2}),
\end{equation}
where we sum over the walks with $N^{-\frac{1}{3}} \max_{0\le t\le M}E(t)\in [h,h+1)$.

Further, to upper bound the product 
\begin{equation}\label{eq_lastjumpchoicefactor}
  \prod_{k=1}^{n}C^{u'(k)}\log(u'(k)+2)^{u'(k)+2}u'(k)^{-u'(k)/2},
\end{equation}
note that at least one $u'(k)>\frac{|U|}{3l}$, and for this $u'(k)$ we have\footnote{We get 1 when the exponent $|U|$ is 0.}

\begin{equation}\label{eq_taillemma7}
u'(k)^{-u'(k)/2} \le \left(\frac{3l}{|U|}\right)^{\frac{|U|}{6l}}. 
\end{equation}
For other $u'(k)$, simply bound $u'(k)^{-u'(k)/2}$ by 1.

Therefore, using the fact that  $\sum_{k=1}^{n}u'(k)+2\le 2|J|+|U|$, we can bound Eqn.~(\ref{eq_lastjumpchoicefactor}) by (taking $h\ge \Z_{\ge 0}$)
\begin{equation}\label{eq_taillemma6}
C^{|U|}\log(2|J|+|U|)^{2|J|+|U|} \left(\frac{3l}{|U|}\right)^{\frac{|U|}{6l}}.
\end{equation}  
Using Eqn.~(\ref{eq_taillemma6}), we multiply Eqn.~(\ref{eq_taillemma5}) by Eqn.~(\ref{eq_taillemma1}) then sum over $h\in \Z_{\ge 1}$. To simplify the notations, let  $$u_{1}=\frac{|U|}{6l},\ \text{and}\ u_{2}=|\mathcal{R}|-\delta_{2}\in \Z_{\ge 0}.$$
 Replace $\sum_{k=1}^{n}u'(k)$ by $\mathfrak{L}-\delta_{3}$, where $|\mathcal{L}|$ is the total number of last jumps. Then the product is a sum of three terms of the following form:
\begin{equation}\label{eq_taillemma4}
\begin{split}
&|U|^{\delta_{2}}(N^{1/3})^{|\mathcal{R}|-\delta_{2}+|\mathcal{L}|-\delta_{3}}C^{|U|}\log(2|J|+|U|)^{2|J|+|U|}(3l)^{u_{1}}\sum_{h\ge0}h^{u_{2}}(|U|)^{-u_{1}}\exp(Ch)\exp(-ch^{2})\\
\lesssim& |U|^{\delta_{2}}(N^{1/3})^{|\mathcal{R}|-\delta_{2}+|\mathcal{L}|-\delta_{3}}\log(2|U|)^{2|U|}C^{u_{2}}\sqrt{(u_{2})!}(|U|)^{-u_{1}}(3l)^{u_{1}} \\
\lesssim&|U|^{\delta_{2}}(N^{1/3})^{|\mathcal{R}|-\delta_{2}+|\mathcal{L}|-\delta_{3}}\log(2|U|)^{2|U|} C^{u_{2}}\left(\frac{u_{2}}{e}\right)^{u_{2}/2}(u_{2})^{\frac{1}{4}}(|U|)^{-u_{1}}(3l)^{u_{1}}.
\end{split}
\end{equation}

\medskip

\noindent{\textbf{Part III. Combining the factors}} 

We multiply all the remaining factors obtained from Steps 0 to 5, namely, multiplying Eqn.~(\ref{eq_taillemma4}) by  Eqn.~(\ref{eq_taillemma0}), (\ref{eq_taillemma2}), (\ref{eq_taillemma3}), and the proper copies of Eqn.~(\ref{eq_conditionalprob1}), (\ref{eq_conditionalprob3}) and (\ref{eq_conditionalprob2}). Then up to a constant independent of the characters, this gives an upper bound 
\begin{equation}\label{eq_sumoverh}
    \begin{split}
        &N^{-\delta}\frac{N^{|U|}}{|U|!}\prod_{j_{p}\in J_{A}:\eta_{j_{p}}\ne\varnothing}(C_{b_{j_{p}}}-\eta_{j_{p}}+1)^{-1/2}\cdot \prod_{p:j_{p}\in J_{IV}}\tau_{j_{p}}^{-1/2}\cdot\prod_{p:j_{p}\in J_{I}\cup J_{III}\cup J_{VI}}\tau_{j_{p}}^{-3/2}\\ 
        \times& N^{-\frac{1}{3}\delta_{1}}N^{-\frac{1}{3}|J_{V}|}N^{-\frac{1}{3}|J_{\varnothing,A}|} N^{\frac{1}{3}(|\mathcal{R}|-\delta_{2}+|\mathcal{L}|-\delta_{3})}\\
        \times&\frac{M^{\delta-|J_{A}|+|J_{\varnothing,A}|-|J_{B}|-|J_{C}|-|U|-|J_{IV}|-|J_{VI}|}}{(k_{1}-l-1)!\cdots (k_{m}-l-1)!(\delta-\sum_{p=1}^{m}k_{p}-2|U|-2l)!}\frac{(\delta+6l-1)!}{(\delta-|U|+4l-1)!}\\
        \times&|U|^{\delta_{2}}\log(2|U|)^{2|U|} C^{u_{2}}\left(\frac{u_{2}}{e}\right)^{u_{2}/2}(u_{2})^{\frac{1}{4}}(|U|)^{-u_{1}}(3l)^{u_{1}}.
    \end{split}
\end{equation}

It remains to upper bound Eqn.~(\ref{eq_sumoverh}), and sum over $k_{1},\ldots,k_{m}$, $u_{1}$ and $u_{2}$.
By definition, \begin{equation}
|U|=6l\cdot u_{1},\quad u_{2}\le\delta-2|U|\le u_{2}+3l.
\end{equation}
Then up to a constant independent of the characters, we have \begin{align*}&\frac{1}{(|U|)!}\frac{(\delta+6l-1)!}{(\delta-|U|+4l-1)!}\le\frac{(\delta+4l-1)!}{|U|!(\delta-|U|+4l-1)!}(\delta+6l)^{2l}\\
\lesssim& 2^{\delta}\delta^{2l}\lesssim 2^{u_{2}+2|U|}(u_{1}+u_{2})^{2l}\lesssim 2^{12l(u_{1}+u_{2})}(u_{1}+u_{2})^{2l}.\end{align*} And for fixed $\delta$, $|U|$ summing over $k_{1},\ldots,k_{p}$, $$\sum_{k_{1},\ldots,k_{p}}\frac{1}{(k_{1}-l-1)!\cdots (k_{m}-l-1)!(\delta-\sum_{p=1}^{m}k_{p}-2|U|-2l)!}\lesssim \frac{(m+1)^{u_{2}}}{(u_{2}-(m+2)l-m)!}.$$
So, up to a constant independent of the characters, the sum of $|w(E,\vec{q})|$ with fixed characters, $\delta$ and $|U|$ is upper bounded by
\begin{equation}
\begin{split}
    &\prod_{j_{p}\in J_{A}:\eta_{j_{p}}\ne \varnothing}(C_{b_{j_{p}}}-Y_{j_{p}}+1)^{-1/2}\cdot \prod_{p:j_{p}\in J_{IV}}\tau_{j_{p}}^{-1/2}\prod_{p:j_{p}\in J_{I}\cup J_{III}\cup J_{VI}}\tau_{j_{p}}^{-3/2}\\
    &N^{-\delta}N^{|U|}\cdot N^{-\frac{1}{3}\delta_{1}}N^{-\frac{1}{3}|J_{V}|}N^{-\frac{1}{3}|J_{\varnothing,A}|}\cdot N^{\frac{2}{3}\left(\delta-|J_{A}|+|J_{\varnothing,A}|-|J_{B}|-|J_{C}|-|U|-|J_{IV}|-|J_{VI}|\right)}N^{\frac{1}{3}(|\mathcal{R}|-\delta_{2}+|\mathcal{L}|-\delta_{3})} \\ \times&\log(12l\cdot u_{1})^{12l\cdot u_{1}}(6l\cdot u_{1})^{\delta_{2}}C^{u_{2}}2^{12l\cdot u_{2}}\frac{(m+1)^{u_{2}}}{(u_{2}-(m+2)l-m)!}u_{1}^{l}(u_{2})^{u_{2}/2}(u_{2})^{\frac{1}{4}}u_{1}^{-u_{1}}(3l)^{u_{1}}(u_{1}+u_{2})^{2l}.
\end{split}
\end{equation}
By Eqn.~(\ref{eq_equality}) and $\delta=|\mathcal{R}|+|\mathcal{L}|+|\mathcal{F}|$, $|\mathcal{F}|=|U|+f_{J}$, $|J|=|J_{A}|+|J_{B}|+|J_C|=|J_{I}|+|J_{III}|+|J_{IV}|+|J_{V}|+|J_{VI}|$, $|J_{III}|=|J_{B}|$, \begin{align*}
    &N^{-\delta}N^{|U|}N^{-\frac{1}{3}\delta_{1}}N^{-\frac{1}{3}|J_{V}|}N^{-\frac{1}{3}|J_{\varnothing,A}|}\cdot N^{\frac{2}{3}\left(\delta-|J_{A}|+|J_{\varnothing,A}|-|J_{B}|-|J_{C}|-|U|-|J_{IV}|-|J_{VI}|\right)}N^{\frac{1}{3}(|\mathcal{R}|-\delta_{2}+|\mathcal{L}|-\delta_{3})} \\
    =&N^{-m}N^{\frac{1}{3}|J_{\varnothing,A}|}N^{-\frac{2}{3}|J_{IV}|-\frac{2}{3}|J_{VI}|}N^{-\frac{1}{3}|J_{V}|-\frac{1}{3}|J_{VI.3}|}\\
    =&N^{-m}N^{-\frac{1}{3}(|J_{A}|-|J_{\varnothing,A}|)}N^{-\frac{1}{3}|J_{IV}|}N^{\frac{1}{3}|J_{I}|-\frac{1}{3}|J_{VI}|-\frac{1}{3}|J_{VI.3}|-\frac{1}{3}|J_{C}|}.
\end{align*} 
Moreover, $$\sum_{u_{1}}\sum_{u_{2}}\log(12l\cdot u_{1})^{12l\cdot u_{1}}(6l\cdot u_{1})^{\delta_{2}}C^{u_{2}}2^{12l\cdot u_{2}}\frac{(m+1)^{u_{2}}}{(u_{2}-(m+2)l-m)!}u_{1}^{l}(u_{2})^{u_{2}/2}(u_{2})^{\frac{1}{4}}u_{1}^{-u_{1}}(3l)^{u_{1}}(u_{1}+u_{2})^{2l}\lesssim 1,$$
since the summand clearly decays exponentially as $u_{1}$ or $u_{2}\rightarrow\infty$. Eqn.~(\ref{eq_tail}) then follows and the proof is finished.

\end{proof}

\begin{remark}
    The reason why we count the jump heights of auxiliary variables backward in time is that, if we count them forward in time, then given $\Vec{q}(t-1) (t\in \llbracket 0,M\rrbracket)$, the number of choices of $\Vec{q}(t)$ would only be bounded by $|U| E(t-1)\le |U|\max_{t}E(t)$ instead of $\max_{t}E(t)$. This is because there are at most $|U|$ choices of jump variables, and each has freedom bounded by the current height of the walk.
\end{remark}

We now combine Lemma \ref{lem:tail} with the cancellation in Section \ref{sec:cancellation}, and give the following estimate that resolves the blow-up issue of Type I local behavior.

\begin{lem}\label{lem:cancellation1}
   There exists a constant $C>1$ that depends on $\bk_{1},\ldots,\bk_{l}$ and $\sigma$, and is uniform in $N$, such that for all $N\ge 100$,
    \begin{equation}
    \begin{split}
        \sum_{(E,\vec{q})\in \mathcal{I}}\left|\sum_{(E,\vec{q}')\sim (E,\vec{q})}w(E,\vec{q}')\right|\le&C \prod_{j_{p}\in J_{A}:\eta_{j_{p}}\ne \varnothing}(C_{b_{j_{p}}}-Y_{j_{p}}+1)^{-1/2}\cdot \prod_{p:j_{p}\in J_{IV}}\tau_{j_{p}}^{-1/2}\prod_{p:j_{p}\in  J_{III}\cup J_{VI}}\tau_{j_{p}}^{-3/2}\\
       \times&\prod_{p:j_{p}\in J_{I}}\tau_{j_{p}}^{-1/2}\cdot N^{-m}N^{-\frac{1}{3}(|J_{A}|-|J_{\varnothing,A}|)}N^{-\frac{1}{3}|J_{IV}|}N^{-\frac{1}{3}|J_{I}|-\frac{1}{3}|J_{VI}|-\frac{1}{3}|J_{VI.3}|-\frac{1}{3}|J_{C}|}
        \end{split}
    \end{equation}
    where the outer sum is over all configurations with given characters $ I(E,\vec{q)},\ \tau_{j_{1}},\dots \tau_{j_{m}},\ \eta_{j_{p}}'s:\ j_{p}\in J_{A}$ of $(E,\vec{q})\in \mathcal{I}$.
\end{lem}
\begin{proof}
    By the discussion in Section \ref{sec:cancellation}, for a given $(E,\vec{q})\in \mathcal{I}$, let $|J|$ and $|U|$ be the number of main and auxiliary variables in $(E,\vec{q})\in \mathcal{I}$, respectively, we have
    \begin{equation}\label{eq_weightcancellation1}
    \resizebox{0.9\textwidth}{!}{$
    \begin{aligned}
        &\left|\sum_{(E,\vec{q}')\sim (E,\vec{q})}w(E,\vec{q}')\right| \le \prod_{p\ j_{p}\in J_{I}}\Bigg[\prod_{\substack{1\le t\le \tau_{j_{p}} \text{ s.t}\\  E(C_{a_{j_{p}}-1}+t)-E(C_{a_{j_{p}}-1}+t-1)=-1}}\left(1+\frac{E(C_{a_{j_{p}}-1}+t)}{N}\right)\\-&\frac{N-|U|-2|J|+1}{N}\prod_{\substack{1\le t< \tau_{j_{p}} \text{ s.t}\\ E(C_{a_{j_{p}}-1}+t)-E(C_{a_{j_{p}}-1}+t-1)=-1}}\left(\frac{N-E(C_{a_{j_{p}}-1})}{N}+\frac{E(C_{a_{j_{p}}-1}+t)-E(C_{a_{j_{p}-1}})}{N}\right)\Bigg]
        \cdot |w(E,\vec{q})|.
    \end{aligned}
    $}
    \end{equation}
    Indeed, for configurations in the same equivalence class, on the interval $[C_{a_{j_p}-1},C_{a_{j_p}-1}+\tau_{j_p}]$ the quantity $q(t)$ satisfies
$
0\le q(t)\le E(C_{a_{j_p}-1})
$
when $x_{j_p}$ is of Type I, and $q(t)\equiv E(C_{a_{j_p}-1})$ when $x_{j_p}$ is of Type II. Also note that the number of variables $x_{j}$ with nonzero degree at a certain time $t$ is at most $q(t)$, which is $E(C_{a_{j_p}-1})$ in the second case above. Hence the factors
\[
1+\frac{E(C_{a_{j_p}-1}+t)}{N}
\quad\text{and}\quad
\frac{N-E(C_{a_{j_{p}}-1})}{N}+\frac{E(C_{a_{j_p}-1}+t)-E(C_{a_{j_p}-1})}{N}
\]
appearing in Eqn.~\eqref{eq_weightcancellation1} provide upper (lower resp.) bounds for the corresponding terms
in Eqn.~\eqref{eq_weight}. Finally, for $x_{j_{p}}$ of Type II, there is one more jump of $x_{j}:\ j\notin J\cup U$ at $C_{a_{j_{p}}-1}+\tau_{j_{p}}$, and $N-|U|-2|J|+1$ gives a lower bound of the number of jump variables  available here. This is how Eqn.~\eqref{eq_weightcancellation1} follows from Eqn.~\eqref{eq_weight}.

    For $E(t):\ N^{-1/3}\max_{0\le t\le M}E(t)\in [h,h+1)$ and $p:\ j_{p}\in J_{I}$, by the triangle inequality we have 
    \begin{align*}
       & \prod_{\substack{1\le t\le \tau_{j_{p}} \text{ s.t}\\ E(C_{a_{j_{p}}-1}+t)-E(C_{a_{j_{p}}-1}+t-1)=-1}}\left(1+\frac{E(C_{a_{j_{p}}-1}+t)}{N}\right)\\
       &\quad\quad\quad-\frac{N-|U|-2|J|+1}{N}\prod_{\substack{1\le t< \tau_{j_{p}} \text{ s.t}\\ E(C_{a_{j_{p}}-1}+t)-E(C_{a_{j_{p}}-1}+t-1)=-1}}\left(\frac{N-E(C_{a_{j_{p}}-1})}{N}+\frac{E(C_{a_{j_{p}}-1}+t)-E(C_{a_{j_{p}-1}})}{N}\right)\\
        \le& C (\tau_{j_{p}} N^{-2/3}+|U|N^{-1})(h+1)\exp(Ch)
    \end{align*}
    for some constant $C>0$ uniform in $N$. When summing over $(E,\vec{q})\in \mathcal{I}$ in the same way as in the proof of Lemma \ref{lem:tail}, the (at most $|J_{I}|$ copies of) extra factors   $C(h+1)\exp(Ch)$ and\footnote{So that the first term dominates the sum.} $C|U|(h+1)\exp(Ch)$ is absorbed in Eqn.~(\ref{eq_taillemma4}), which gives the desired bound.
 \end{proof}

Let $\mathcal{C}$ denote the set of configurations $(E,\vec{q})$ such that $|J_{VI}|=|J_{C}|=0$.
\begin{cor}\label{cor:typeIVandC}
    There exists a constant $C>1$ that depends on $\bk_{1},\ldots,\bk_{l}$ and $\sigma$, and is uniform in $N$, such that for all $N\ge 100$,
    \begin{equation}
    \begin{split}
    \left|\sum_{(E,\vec{q})}w(E,\vec{q})\right|\le& C\cdot N^{-m};\\
    \left|\sum_{(E,\vec{q})\notin \mathcal{C}}w(E,\vec{q})\right|\le& C\cdot N^{-m-1/3}.
   \end{split}
    \end{equation}
\end{cor}
\begin{proof}
     Note that  $1\le \tau_{j_{p}}\le M\le (\sum_{i=1}^{l}T_{i}+1)N^{2/3}$ for $\tau_{j_{p}}\ne \varnothing$, and  $1\le C_{b_{j_{p}}}-Y_{j_{p}}+1\le M\le (\sum_{i=1}^{l}T_{i}+1)N^{2/3}$. Then this follows directly from  Lemma \ref{lem:cancellation1} by summing over the characters.
\end{proof}
By Corollary \ref{cor:typeIVandC}, in the remainder of this section, it suffices for us to consider only $(E,\vec{q})\in \mathcal{C}$.

\subsection{A second type of cancellation}\label{sec:secondcancellation}
In this section, we introduce a second  cancellation mechanism. Fix $N\ge 100$ and $\bk_{1},\ldots,\bk_{l}>0$, and consider only the configurations in the set $\mathcal{C}$, i.e., $|J_{VI}|=|J_{C}|=0$. 

For $(i_{1},\ldots,i_{l})\in \llbracket1,N\rrbracket$, recall that a configuration $(E,\vec{q})\in \mathcal{B}(i_{1},\ldots,i_{l})$ if it is created by the action above with this set of indices $i_{j}$. We write $(E,\vec{q})\in \tilde{\mathcal{B}}(i_{1},\ldots,i_{l})\subset \mathcal{B}(i_{1},\ldots,i_{l})$ if none of its main variables is of Type III. We then say that $(E,\vec{q}')\in \mathcal{B}(i'_{1},\ldots,i'_{l})$ \emph{deforms to} $(E,\vec{q})\in \tilde{\mathcal{B}}(i_{1},\ldots,i_{l})$, and write $R(E,\vec{q}')=(E,\vec{q})$, if  $(E,\vec{q}')$ has at least one main variable of Type III, and one can obtain $\vec{q}$ and $(i_{1},\ldots,i_{l})$ from $\vec{q}'$ and $(i'_{1},\ldots,i'_{l})$ by the following algorithm:
\begin{itemize}
    \item  Let $i'$ denote the element in $J'_{III}$, the set of Type III main variables of $(E,\vec{q}')$, that has the largest $a_{i'}$ (defined in Eqn.~(\ref{eq_aibi})). Let $j$ be the index in $J_{B}$ where $q_{j}(C_{a_{i'}-1}+\tau_{i'}-)>0$ and $q_{j}(t)=0$ for $t\ge C_{a_{i'}-1}+\tau_{i'}$, namely, the last jump of $x_{j}$ occurs at time $C_{a_{i'}-1}+\tau_{i'}$. Then remove the last jump of $x_{j}$ at $C_{a_{i'}-1}+\tau_{i'}$, replace all the main variables $x_{i'_{k}}\ (1\le k\le l)$  by $x_{j}$ for all $i'_{k}=i'$, and replace all the jumps of $x_{i'}$ by jumps of $x_{j}$. 
    \item Repeat the above process until $|J'_{III}|=0$.
\end{itemize}

\begin{example}
    Take $i_{1}=i_{2}=i_{4}=1$, $i_{3}=2$. Then $R(E,\vec{q}')=(E,\vec{q})\in \tilde{\mathcal{B}}(1,1,2,1)$  for those $(E,\vec{q}')\in \mathcal{B}(1,i_{2}',2,i_{4}')$ where $i_{2}'\ne i_{4}'$ are elements of $\llbracket 3,N\rrbracket$, $x_{i_2'}$ has its last jump at $C_{3}+\tau_{i_{4}'}$, and $x_{1}$ has its last jump at $C_{1}+\tau_{i_{2}'}$. Figure \ref{fig:cancel2} gives an example of the corresponding configurations.
\end{example}

\begin{figure}
    \centering
    \[
    \resizebox{0.8\textwidth}{!}{
    \begin{tikzpicture}[shorten >=-3pt,shorten <=-3pt]
        \draw[ultra thin, opacity = 0.5] (0,0) grid (28,8);
        \draw[ultra thin,dashed] (0,-0.5)--(0,8.5);
        \draw[ultra thin,dashed] (7,-0.5)--(7,8.5);
        \draw[ultra thin,dashed] (14,-0.5)--(14,8.5);
        \draw[ultra thin,dashed] (21,-0.5)--(21,8.5);
        \draw[ultra thin,dashed] (28,-0.5)--(28,8.5);

        \draw[ultra thick,*-o] (0,0)--(1,0);
        \draw[ultra thick,*-o] (1,2)--(2,2);
        \draw[ultra thick,*-o] (2,1)--(3,1);
        \draw[ultra thick,*-o] (3,4)--(5,4);
        \draw[ultra thick,*-o] (5,5)--(6,5);
        \draw[ultra thick,*-o] (6,4)--(9,4);
        \draw[ultra thick,*-o] (9,3)--(10,3);
        \draw[ultra thick,*-o] (10,7)--(11,7);
        \draw[ultra thick,*-o] (11,6)--(12,6);
        \draw[ultra thick,*-o] (12,5)--(13,5);
        \draw[ultra thick,*-o] (13,6)--(14,6);
        \draw[ultra thick,*-o] (14,5)--(15,5);
        \draw[ultra thick,*-o] (15,8)--(16,8);
        \draw[ultra thick,*-o] (16,7)--(17,7);
        \draw[ultra thick,*-o] (17,6)--(18,6);
        \draw[ultra thick,*-o] (18,5)--(19,5);
        \draw[ultra thick,*-o] (19,4)--(20,4);
        \draw[ultra thick,*-o] (20,3)--(22,3);
        \draw[ultra thick,*-o] (22,4)--(23,4);
        \draw[ultra thick,*-o] (23,3)--(24,3);
        \draw[ultra thick,*-o] (24,2)--(25,2);
        \draw[ultra thick,*-o] (25,3)--(26,3);
        \draw[ultra thick,*-o] (26,2)--(27,2);
        \draw[ultra thick,*-o] (27,1)--(28,1);
        \draw[ultra thick,*-] (28,0);

        \draw[ultra thick, blue,*-o] (0,0.2)--(7,0.2);
        \draw[ultra thick, blue,*-o] (7,4.1)--(9,4.1);
        \draw[ultra thick, blue,*-o] (9,0.2)--(14,0.2);
        \draw[ultra thick, blue,*-o] (14,5.1)--(16,5.1);
        \draw[ultra thick, blue,*-o] (16,2.1)--(18,2.1);
        \draw[ultra thick, blue,*-o] (18,3.1)--(24,3.1);
        \draw[ultra thick, blue,*-*] (24,0.2)--(28,0.2);

        \draw[ultra thick, red,*-o] (0,0.1)--(7,0.1);
        \draw[ultra thick, red,*-o] (7,3.9)--(9,3.9);
        \draw[ultra thick, red,*-*] (9,0.1)--(28,0.1);

        \draw[ultra thick, green,*-o] (0,0)--(14,0);
        \draw[ultra thick, green,*-o] (14,4.9)--(16,4.9);
        \draw[ultra thick, green,*-o] (16,2)--(18,2);
        \draw[ultra thick, green,*-o] (18,2.9)--(24,2.9);
        \draw[ultra thick, green,*-*] (24,0)--(28,0);

        \node[below] at (0,-0.5) {$C_0$};
        \node[below] at (7,-0.5) {$C_1$};
        \node[below] at (14,-0.5) {$C_2$};
        \node[below] at (21,-0.5) {$C_3$};
        \node[below] at (28,-0.5) {$C_4$};
        \node[] at (3.5,8.5) {$1$};
        \node[] at (10.5,8.5) {$i_2'$};
        \node[] at (17.5,8.5) {$2$};
        \node[] at (24.5,8.5) {$i_4'$};

        \draw[<-] (21.4,-0.5)--(22.1,-0.5);
        \draw[->] (22.9,-0.5)--(23.6,-0.5);
        \node at (22.5,-0.5) {$\tau_{i_4'}$};

        \draw[<-] (7.4,-0.5)--(7.6,-0.5);
        \draw[->] (8.4,-0.5)--(8.6,-0.5);
        \node at (8,-0.5) {$\tau_{i_2'}$};
    \end{tikzpicture}
    }
\]

\[
    \resizebox{0.8\textwidth}{!}{
    \begin{tikzpicture}[shorten >=-3pt,shorten <=-3pt]
        \draw[ultra thin, opacity = 0.5] (0,0) grid (28,8);
        \draw[ultra thin,dashed] (0,-0.5)--(0,8.5);
        \draw[ultra thin,dashed] (7,-0.5)--(7,8.5);
        \draw[ultra thin,dashed] (14,-0.5)--(14,8.5);
        \draw[ultra thin,dashed] (21,-0.5)--(21,8.5);
        \draw[ultra thin,dashed] (28,-0.5)--(28,8.5);

        \draw[ultra thick,*-o] (0,0)--(1,0);
        \draw[ultra thick,*-o] (1,2)--(2,2);
        \draw[ultra thick,*-o] (2,1)--(3,1);
        \draw[ultra thick,*-o] (3,4)--(5,4);
        \draw[ultra thick,*-o] (5,5)--(6,5);
        \draw[ultra thick,*-o] (6,4)--(9,4);
        \draw[ultra thick,*-o] (9,3)--(10,3);
        \draw[ultra thick,*-o] (10,7)--(11,7);
        \draw[ultra thick,*-o] (11,6)--(12,6);
        \draw[ultra thick,*-o] (12,5)--(13,5);
        \draw[ultra thick,*-o] (13,6)--(14,6);
        \draw[ultra thick,*-o] (14,5)--(15,5);
        \draw[ultra thick,*-o] (15,8)--(16,8);
        \draw[ultra thick,*-o] (16,7)--(17,7);
        \draw[ultra thick,*-o] (17,6)--(18,6);
        \draw[ultra thick,*-o] (18,5)--(19,5);
        \draw[ultra thick,*-o] (19,4)--(20,4);
        \draw[ultra thick,*-o] (20,3)--(22,3);
        \draw[ultra thick,*-o] (22,4)--(23,4);
        \draw[ultra thick,*-o] (23,3)--(24,3);
        \draw[ultra thick,*-o] (24,2)--(25,2);
        \draw[ultra thick,*-o] (25,3)--(26,3);
        \draw[ultra thick,*-o] (26,2)--(27,2);
        \draw[ultra thick,*-o] (27,1)--(28,1);
        \draw[ultra thick,*-] (28,0);

        \draw[ultra thick, blue,*-o] (0,0.2)--(7,0.2);
        \draw[ultra thick, blue,*-o] (7,4.1)--(9,4.1);
        \draw[ultra thick, blue,*-o] (9,0.2)--(14,0.2);
        \draw[ultra thick, blue,*-o] (14,5.1)--(16,5.1);
        \draw[ultra thick, blue,*-o] (16,2.1)--(18,2.1);
        \draw[ultra thick, blue,*-o] (18,3.1)--(21,3.1);
        \draw[ultra thick, blue,*-*] (21,0.2)--(28,0.2);

        \draw[ultra thick, red,*-o] (0,0.1)--(7,0.1);
        \draw[ultra thick, red,*-o] (7,3.9)--(9,3.9);
        \draw[ultra thick, red,*-*] (9,0.1)--(28,0.1);

        \draw[ultra thick, green,*-o] (0,0)--(14,0);
        \draw[ultra thick, green,*-o] (14,4.9)--(16,4.9);
        \draw[ultra thick, green,*-o] (16,2)--(18,2);
        \draw[ultra thick, green,*-o] (18,2.9)--(21,2.9);
        \draw[ultra thick, green,*-*] (21,0)--(28,0);

        \node[below] at (0,-0.5) {$C_0$};
        \node[below] at (7,-0.5) {$C_1$};
        \node[below] at (14,-0.5) {$C_2$};
        \node[below] at (21,-0.5) {$C_3$};
        \node[below] at (28,-0.5) {$C_4$};
        \node[] at (3.5,8.5) {$1$};
        \node[] at (10.5,8.5) {$i_2'$};
        \node[] at (17.5,8.5) {$2$};
        \node[] at (24.5,8.5) {$i_2'$};

        \draw[<-] (7.4,-0.5)--(7.6,-0.5);
        \draw[->] (8.4,-0.5)--(8.6,-0.5);
        \node at (8,-0.5) {$\tau_{i_2'}$};
    \end{tikzpicture}
    }
\]

\[
    \resizebox{0.8\textwidth}{!}{
    \begin{tikzpicture}[shorten >=-3pt,shorten <=-3pt]
        \draw[ultra thin, opacity = 0.5] (0,0) grid (28,8);
        \draw[ultra thin,dashed] (0,-0.5)--(0,8.5);
        \draw[ultra thin,dashed] (7,-0.5)--(7,8.5);
        \draw[ultra thin,dashed] (14,-0.5)--(14,8.5);
        \draw[ultra thin,dashed] (21,-0.5)--(21,8.5);
        \draw[ultra thin,dashed] (28,-0.5)--(28,8.5);

        \draw[ultra thick,*-o] (0,0)--(1,0);
        \draw[ultra thick,*-o] (1,2)--(2,2);
        \draw[ultra thick,*-o] (2,1)--(3,1);
        \draw[ultra thick,*-o] (3,4)--(5,4);
        \draw[ultra thick,*-o] (5,5)--(6,5);
        \draw[ultra thick,*-o] (6,4)--(9,4);
        \draw[ultra thick,*-o] (9,3)--(10,3);
        \draw[ultra thick,*-o] (10,7)--(11,7);
        \draw[ultra thick,*-o] (11,6)--(12,6);
        \draw[ultra thick,*-o] (12,5)--(13,5);
        \draw[ultra thick,*-o] (13,6)--(14,6);
        \draw[ultra thick,*-o] (14,5)--(15,5);
        \draw[ultra thick,*-o] (15,8)--(16,8);
        \draw[ultra thick,*-o] (16,7)--(17,7);
        \draw[ultra thick,*-o] (17,6)--(18,6);
        \draw[ultra thick,*-o] (18,5)--(19,5);
        \draw[ultra thick,*-o] (19,4)--(20,4);
        \draw[ultra thick,*-o] (20,3)--(22,3);
        \draw[ultra thick,*-o] (22,4)--(23,4);
        \draw[ultra thick,*-o] (23,3)--(24,3);
        \draw[ultra thick,*-o] (24,2)--(25,2);
        \draw[ultra thick,*-o] (25,3)--(26,3);
        \draw[ultra thick,*-o] (26,2)--(27,2);
        \draw[ultra thick,*-o] (27,1)--(28,1);
        \draw[ultra thick,*-] (28,0);

        \draw[ultra thick, blue,*-o] (0,0.1)--(14,0.1);
        \draw[ultra thick, blue,*-o] (14,5.1)--(16,5.1);
        \draw[ultra thick, blue,*-o] (16,2.1)--(18,2.1);
        \draw[ultra thick, blue,*-o] (18,3.1)--(21,3.1);
        \draw[ultra thick, blue,*-*] (21,0.1)--(28,0.1);

        \draw[ultra thick, red,*-o] (0,0)--(14,0);
        \draw[ultra thick, red,*-o] (14,4.9)--(16,4.9);
        \draw[ultra thick, red,*-o] (16,2)--(18,2);
        \draw[ultra thick, red,*-o] (18,2.9)--(21,2.9);
        \draw[ultra thick, red,*-*] (21,0)--(28,0);

        \node[below] at (0,-0.5) {$C_0$};
        \node[below] at (7,-0.5) {$C_1$};
        \node[below] at (14,-0.5) {$C_2$};
        \node[below] at (21,-0.5) {$C_3$};
        \node[below] at (28,-0.5) {$C_4$};
        \node[] at (3.5,8.5) {$1$};
        \node[] at (10.5,8.5) {$1$};
        \node[] at (17.5,8.5) {$2$};
        \node[] at (24.5,8.5) {$1$};

    \end{tikzpicture}
    }
\]
    \caption{An example of the deformation algorithm given in Section \ref{sec:secondcancellation}. In each configuration black is the walk, blue is the sum of the blocks, red is the block for $x_1$, and green is the block for $x_{i_2'}$. The labels above each interval denote the index of the corresponding main variable. \textbf{Top:}  We start with a configuration in $\mathcal{B}(1,i_{2}',2,i_{4}')$ where $i_{2}'\ne i_{4}'$ are two distinct elements in $\llbracket 3,N\rrbracket$. Note that $x_{i_2'}$ is of Type III with $x_1$ the corresponding variable of Type B, and $x_{i_4'}$ is of Type III with $x_{i_2'}$ the corresponding variable of Type B. To deform the configuration, we start with the rightmost Type III variable and delete the jump of $x_{i_2'}$ between $C_3$ and $C_4$, then change the main variable for this interval to $x_{i_2'}$. \textbf{Middle:} After the first deformation, $x_{i_2'}$ is still of Type III  with $x_1$ the corresponding variable of Type B. We deform the configuration by deleting the jump of $x_1$ between $C_1$ and $C_2$ and changing all instances of the variable $x_{i_2'}$ to $x_1$. This changes the main variable in the intervals from $C_1$ to $C_2$ and $C_3$ to $C_4$ from $x_{i_2'}$ to $x_1$, as well as changing the block of $x_{i_2'}$ to a block of $x_1$. \textbf{Bottom:} The end result is a configuration in $\mathcal{B}(1,1,2,1)$ without any Type III variables.}
    \label{fig:cancel2}
\end{figure}

We group $(E,\vec{q})$ with all $(E,\vec{q}')$ such that $(E,\vec{q})=R(E,\vec{q}')$. Then, similarly to Section \ref{sec:cancellation}, the weights of configurations in the same group have opposite signs and cancel out. More precisely, we have the following estimate.

\begin{lem}\label{lem:cancellation2}
    Fix a tuple $(i_{1},\ldots,i_{l})$ that is not composed of mutually distinct indices, and define the set $J$ from $(i_{1},\ldots,i_{l})$ as before. Then there exists a constant $C>1$ that depends on $\bk_{1},\ldots,\bk_{l}$ and $\sigma$, and is uniform in $N$, such that for all $N\ge 100$, 
    \begin{equation}
    \sum_{(E,\vec{q})\in \mathcal{\tilde{B}}(i_{1},\ldots,i_{l})}\left|w(E,\vec{q})+\sum_{(E,\vec{q}'):\ R(E,\vec{q}')=(E,\vec{q})}w(E,\vec{q}'
        )\right|\le C\cdot N^{-|J|} N^{-1/12}.
    \end{equation} 
\end{lem}
\begin{proof}
From the definition, one can see that $(E,\vec{q})=R(E,\vec{q}')$ for some $(E,\vec{q}')\in \mathcal{B}(i_{1}',\ldots,i_{l}')$ with mutually distinct $i'_{1},\ldots,i_{l}'$, unless\footnote{Here, the exponent $1/4$ can be taken as any real number in $(0,\frac{1}{3})$.}
\begin{itemize}
    \item $i_{1},\ldots,i_{l}$ are mutually distinct. 
    \item There exists some index $i\notin\{a_{j}:j\in J\}$, such that $E(C_{i-1}+t)\ge E(C_{i-1})$ for all $0\le t\le N^{1/4}$, or
    \item there is a jump at $C_{i-1}+t$ for some $1\le t\le N^{1/4}$.
\end{itemize}
    We briefly explain how the second and third situations reduce the weight contribution of the configurations, recalling the argument in the proof of Lemma \ref{lem:tail}. Again, viewing $E(t)$ as a conditional random walk bridge, there exists a constant $C>1$ that depends on $\bk_{1},\ldots,\bk_{l}$ and $\sigma$ and is uniform in $N$, such that the second case above occurs with probability $<C N^{-1/8}$. This gives an additional factor $N^{-1/8}$. If the third case occurs, then there are $N^{1/4}$ choices of the jump time (and at most $|J|+|U|$ choices of the jump variable). Fixing the jump time, if the jump is a first or free jump, then we update the bound in Eqn.~(\ref{eq_taillemma2}) by decreasing one of the elements $1,\ldots,1,k_{1}',\ldots,k_{m}'$ and $\delta-\sum_{p=1}^{m}k'_{p}-2|U|-|J_{IV}|-|J_{m,VI}|-|J_{a,VI}|-|J_{m,A}|-|J_{a,A}|-|J_{B}|-|J_{C}|$ by 1. This gives an additional factor $M^{-1}\lesssim N^{-2/3}$. If the jump is a last jump, then its height is fixed by the walk $E(t)$ at the given jump time, and we decrease one of the $u'(\tilde{k})$ in Eqn.~(\ref{eq_step8factor}) by 1. This gives an additional factor $N^{-1/3}$. Together with the factor $N^{1/4}$, the total weight of the configurations with the third situation above is at most $O(N^{-1/12})$.

It then remains to consider the set of $(E,\vec{q})$ which have a collection of preimages $(E,\vec{q}')$ with mutually distinct $i_{1}',\ldots,i_{l}'$, such that, for $1\le j\le l$, $i_{j}'\in \llbracket1,N\rrbracket\setminus J$ if $i_{j}\notin \{a_{j}:\ j\in J\}$. For this $(E,\vec{q}')$, let $J'_{III}$ denote the set of its type III main variables. Then each element $(E,\vec{q}'')$ in  $R^{-1}(E,\vec{q})$ can be deformed from $(E,\vec{q}')$ by 
\begin{itemize}
    \item Take a subset $J''_{III}$ of $J'_{III}$.
    \item  Let $i$ denote the element in $J''_{III}$ that has the largest $a_{i}$. Let $j$ be the index in $J_{B}$ where $q_{j}(C_{a_{i}-1}+\tau_{i}-)>0$ and $q_{j}(t)=0$ for $t\ge C_{a_{i}-1}+\tau_{i}$, that is, the last jump of $x_{j}$ occurs at time $C_{a_{i}-1}+\tau_{i}$. Then let $i_{k}\ (1\le k\le l)$  be $j$ for all $i_{k}=i$, and let $q_{i}\equiv0$, that is, remove the only jump of $x_{j}$ at $C_{a_{i}-1}+\tau_{i}$.
    \item Repeat the above process until $|J''_{III}|=0$.
    \end{itemize}
Moreover, it is clear that $(E,\vec{q}'')$ is indexed by $J_{III}''$. Similarly to Eqn.~(\ref{eq_weightcancellation1}), for $E(t)$ such that $\ N^{-1/3}\max_{0\le t\le M}E(t)\in [h,h+1)$, we have
\begin{equation}
\resizebox{0.9\textwidth}{!}{$
    \begin{aligned}
        &\left|w(E,\vec{q})+\sum_{(E,\vec{q}''):\ R(E,\vec{q}'')=(E,\vec{q})}w(E,\vec{q}''
        )\right|\\
        \le & \prod_{p:\ j'_{p}\in J'_{III}}\Bigg[\prod_{\substack{1\le t\le \tau_{j'_{p}}\text{ s.t}\\ E(C_{a_{j_{p}}-1}+t)-E(C_{a_{j_{p}}-1}+t-1)=-1}}\left(1+\frac{E(C_{a_{j'_{p}}-1}+t)}{N}\right)\\-&\frac{N-|U|-2|J|+1}{N}\prod_{\substack{1\le t< \tau_{j'_{p}}\text{ s.t}\\ E(C_{a_{j_{p}}-1}+t)-E(C_{a_{j_{p}}-1}+t-1)=-1}}\left(\frac{N-E(C_{a_{j'_{p}}-1})}{N}+\frac{E(C_{a_{j'_{p}}-1}+t)-E(C_{a_{j'_{p}}-1})}{N}\right)\Bigg]
        \cdot |w(E,\vec{q}')|\\
        \le & \sum_{(E,\vec{q}'):\ i_{1}'\ne\ldots\ne i_{l}'}\prod_{p:\ j'_{p}\in J'_{III}}\left(C(\tau_{j'_{p}}N^{-2/3}+|U|N^{-1})(h+1)\exp(Ch)\right)\cdot |w(E,\vec{q}')|\\ \le& N^{|J_{III}'|}\prod_{p:\ j'_{p}\in J'_{III}}\left(C(\tau_{j'_{p}}N^{-2/3}+|U|N^{-1})(h+1)\exp(Ch)\right)\cdot |w(E,\vec{q}')|,
    \end{aligned}
    $}
\end{equation}
where in the last expression we pick an arbitrary $(E,\vec{q}')\in R^{-1}(E,\vec{q})$ with mutually distinct $i_{1}',\ldots,i_{l}'$ with $|J|$ main variables and $|U|$ auxiliary variables.
Summing over $(E,\vec{q})\in \tilde{\mathcal{B}}(i_{1},\ldots,i_{l})$ as in the  Lemma \ref{lem:cancellation1} and Corollary \ref{cor:typeIVandC}, we get an upper bound \footnote{Note that $l-|J|=|J_{III}'|$.}
\[C\cdot N^{-|J|}N^{-\frac{1}{3}(l-|J|)}\le C\cdot N^{-|J|}N^{-1/12},\]
where the last inequality holds if $i_{1},\ldots,i_{l}$ are not mutually distinct.
\end{proof}
\begin{cor}\label{cor:TypeIII}
    \begin{align*}
   & \frac{1}{\mu_{+}(N)^{M_{1}+\ldots+M_{l}}}\cdot \prod_{j=1}^{l}\left[\left(\frac{\mathcal{D}_{1}}{N}\right)^{M_{j}}+\ldots+\left(\frac{\mathcal{D}_{N}}{N}\right)^{M_{j}}\right]G_{N}(\vec{x};\beta)\Bigg|_{\vec{x}=0}\\
    =&\sum_{\substack{(i_{1},\ldots,i_{l})\\ \text{mutually\ distinct}}}\sum_{(E,\vec{q})\in \tilde{\mathcal{B}}(i_{1},\ldots,i_{l})\cap\mathcal{C}}w(E,\vec{q}) +o(1)=N(N-1)\cdots(N-l+1)\sum_{(E,\vec{q})\in \tilde{\mathcal{B}}(1,\ldots,l)\cap\mathcal{C}}w(E,\vec{q}) +o(1),
\end{align*}
where the $o(1)$ term vanishes as $N\to \infty$.
\end{cor}
\begin{proof}
    The left-hand side of the equality above is equal to 
    \begin{align*}
        \sum_{\substack{(i_{1},\ldots,i_{l})\\ \text{mutually\ distinct}}}&\sum_{(E,\vec{q})\in \tilde{\mathcal{B}}(i_{1},\ldots,i_{l})\cap\mathcal{C}}w(E,\vec{q})\\
        + & \sum_{\substack{(i_{1},\ldots,i_{l})\\ \text{ not\ mutually\ distinct}}}\sum_{(E,\vec{q})\in \mathcal{B}(i_{1},\ldots,i_{l})}\left(w(E,\vec{q})+\sum_{(E,\vec{q}'):\ R(E,\vec{q}')=(E,\vec{q})}w(E,\vec{q}')\right)
        .
    \end{align*}
    Summing over the $O(N^{|J|})$ choices of $(i_{1},\ldots,i_{l})$, the second line is of order $o(1)$ by Lemma \ref{lem:cancellation2}. This verifies the first equality. The second equality follows directly from the symmetry of the Dunkl operators and Bessel generating function.
\end{proof}

The main contribution in Eqn.~(\ref{eq_dunklpowersums}) is therefore given by configurations with mutually distinct $i_{1},\ldots,i_{l}$ and without main variables of Types III, VI, B, or C.

\subsection{Convergence of the functionals}\label{sec:convergence}

In this section, we prove that the Dunkl action Eqn.~(\ref{eq_dunklpowersums}) converges to the mixed moment of the Laplace transform of the $\mathrm{Airy}(\beta)$ process. 

By Corollary \ref{cor:TypeIII}, it suffices to consider configurations $(E,\vec{q})$ with the following properties. 
\begin{itemize}
    \item We take $J=\{1,2,\ldots,l\}$, $a_{i}=b_{i}=i$ for $i\in J$. $U=\{l+1,\ldots,l+u\}$ for some $u\in \Z_{\ge 0}$. In other words, we take the auxiliary variables to be $l+1,\ldots,l+u$ for notational simplicity.
    \item Let $\delta_{i,p}\in \Z_{\ge 0}$ denote the number of jumps of $x_{i}$ in $\llbracket C_{p-1}+1, C_{p}\rrbracket$. Then $\delta_{i,p}=0$ for $1\le i\le l,\ p\ge i$. Let $\vec{\delta}=(\delta_{i,p}:1\le i\le l+u,\ 1\le p\le l)$ be the tuple that records all the $\delta_{i,p}$’s (and hence encodes the number of auxiliary variables in $(E,\vec{q})$ via $u$).
\end{itemize}

Take $\epsilon>0$. Let $\tilde{\mathcal{B}}_{\epsilon}(1,2,\ldots,l)$ be the subset of $\mathcal{B}(1,2,\ldots,l)$ such that for each $p=2,\dots,l$, we have $\tau_p>\epsilon(C_p-C_{p-1})$ or $\tau_p=\varnothing$, and  $\tilde\tau_p>\epsilon(C_p-C_{p-1})$ or $\tilde\tau_p=\varnothing$.
\begin{lem}\label{lem:principalvalue}
There exists a constant $C>0$ uniform in $N$, such that
  \begin{align*}
   & \frac{1}{\mu_{+}(N)^{M_{1}+\ldots+M_{l}}}\cdot \prod_{j=1}^{l}\left[\left(\frac{\mathcal{D}_{1}}{N}\right)^{M_{j}}+\ldots+\left(\frac{\mathcal{D}_{N}}{N}\right)^{M_{j}}\right]G_{N}(\vec{x};\beta)\Bigg|_{\vec{x}=0}\\
    =&N(N-1)\cdots (N-l+1)\sum_{(E,\vec{q})\in \tilde{\mathcal{B}}_{\epsilon}(1,\ldots,l)\cap\mathcal{C}}w(E,\vec{q})+C\sqrt{\epsilon} +o(1).
    \end{align*}
\end{lem}
\begin{proof}
    By Corollary \ref{cor:TypeIII}, it remains to check that \[N(N-1)\cdots (N-l+1)\left|\sum_{(E,\vec{q})\in \tilde{\mathcal{B}}(1,\ldots,l)\setminus \tilde{\mathcal{B}}_{\epsilon}(1,2,\ldots,l)}w(E,\vec{q})\right|\le C\sqrt{\epsilon}.\]
    This follows from Lemma \ref{lem:cancellation1} and a similar summation over the characters as in the proof of Corollary \ref{cor:typeIVandC}.
\end{proof}

In the remainder of this section, we will classify the elements in $\tilde{\mathcal{B}}_{\epsilon}(1,2,\ldots,l)\cap \mathcal{C}$ by their data $\vec{\delta}$.
Given a fixed $\vec{\delta}$, in order to characterize the possible Type I local behavior of $(E,\vec{q})$, we introduce the \emph{virtual blocks}. These are a collection of $\nu_{p}\in \Z_{\ge 1}\cup \varnothing,\ 2\le p\le l$, such that\footnote{One can also define $\nu_{1}$, which by definition is always $\varnothing$.} 
\begin{itemize}
    \item $\nu_{p}=\varnothing$ if $\sum_{k=1}^{p-1}\delta_{p,k}=0$,
    \item if $\sum_{k=1}^{p-1}\delta_{p,k}>0$, $\nu_{p}\ne \varnothing$ if and only if $\tau_{p}\ne \varnothing$, and $\tilde{\tau}_{p}>\tau_{p}$ or $\tilde{\tau}_{p}=\varnothing$. In this case $\nu_{p}=\tau_{p}$. 
\end{itemize}

Given $\vec{\delta}$, $\Delta$ and $\vec{\nu}=(\nu_{1},\ldots,\nu_{l})$, let 
$$\resizebox{\textwidth}{!}{$(s_{0}<\ldots<s_{n})=(\Delta-1)\cup \Delta\cup  \{C_{p}:\ 0\le p\le l\}\cup  \{C_{p-1}+\nu_{p}-1:\ \nu_{p}\ne\varnothing,\ 2\le p\le l\}\cup \{C_{p-1}+\nu_{p}:\ \nu_{p}\ne\varnothing,\ 2\le p\le l\}.$}$$ 
Suppose $(E,\vec{q})$ is also given, let $\vec{H}=(E(s_{1}),\ldots,E(s_{n}))$ record the height of $E(t)$ at the points $s_i$, which split $E(t)$ as concatenation of segments where $q(t)$ is fixed on each segment $[s_{k},s_{k+1})$. 

The data $\vec{\delta}$, $\Delta$ and $\vec{\nu}$ gives us a way to determine $(E,\vec{q})$ by the following steps: 
\begin{itemize}
    \item Fix $\vec{\delta}$, $\Delta$ and $\vec{\nu}$.
    \item Determine  the $\sum_{i=1}^{l}\sum_{p=1}^{i-1}\delta_{i,p}+\sum_{i=l+1}^{l+u}\sum_{p=1}^{l}\delta_{i,p}$ jump times and $\sum_{i=1}^{l}\sum_{p=1}^{i-1}\delta_{i,p}+\sum_{i=l+1}^{l+u}\left(\sum_{p=1}^{l}\delta_{i,p}-1\right)$ jump heights of all the jumps  (the second sum has a term $-1$ because for $i\in U$, the last jump height is fixed by the previous jump heights). This completely determines $\vec{q}=(q_{1}(t),\ldots,q_{l+u}(t))$, and in particular $q(t)=\sum_{i\ne p}q_{i}(t)$ where $t\in [C_{p-1},C_{p})$. Moreover, each $\vec{q}$ is in an equivalence class of $u!$ elements upon swapping the indices $l+1,\ldots,l+u$, so we assume that $$\min\{t:q_{l+1}(t)>0\}<\ldots<\min\{t:q_{l+u}(t)>0\},$$ which gives a complete order $x_{l+1}\prec\cdots\prec x_{l+u}$.
    \item Determine $\vec{H}$ under the restrictions that: 
    
    1. $E(s_{k})\ge q(s_{k})$ for all $0\le k\le n$. 
    
    2. If $s_{k}=C_{p}$ for some $0\le p\le l$, then $E(s_{k})=q(s_{k})$;
    
    3. If $s_{k}=C_{p-1}+\nu_{p}-1$ and $s_{k+1}=C_{p-1}+\nu_{p}$ for some $2\le p\le l$, or $s_{k}\in \Delta-1$ and $s_{k+1}\in \Delta$ is the last jump time for some $x_{i}$, $l+1\le i\le l+u$, then $E(s_{k+1})=E(s_{k})-1$. 

    4. By Eqn.~(\ref{eq_jumprules}), if $s_{k}\in \Delta$, and $i$ is the index that $q_{i}(s_{k})\ne q_{i}(s_{k}-)$,
    \begin{equation}\label{eq_jumprules2}
      \begin{split}
    &|q_{i}(s_{k})-q_{i}(s_{k}-)| \le |E(s_{k}-)-q(s_{k}-)-q_i(s_{k}-)-\mathbf{1}_{q_{i}(s_{k})-q_{i}(s_{k}-)>0}|;\\
     &[q_{i}(s_{k})- q_{i}(s_{k}-)]\cdot (E(s_{k}-)-q(s_{k}-)-q_i(s_{k}-)) > 0.
    \end{split}
\end{equation}
    \item Denote $E_{k}(t)\ (0\le k\le n-1)$ as the restriction of $E(t)$ on $[s_{k},s_{k+1})$. Determine the rest of $E_{k}(t)$ under the restriction that all
$E_{k}(s_{k})$ are fixed by $\vec{H}$. $E_{k}(t)\ge q(t)$ for all $t\in [0,M]$, which implies:

    0. If $s_{k+1}\in \Delta\cup \{C_{p-1}+\nu_{p}:\ \nu_{p}\ne\varnothing,\ 1\le p\le l\}$, then $E(s_{k+1})=E(s_{k})-1$, and we denote this as $(s_{k},s_{k+1})\in \Xi_{0}$.

    1. If $s_{k}=C_{p-1}$ for some $2\le p\le l$, and \begin{itemize}
        \item $s_{k+1}=C_{p-1}+\nu_{p}-1$, or
        \item $s_{k+1}+1$ is the last jump time for some $x_{i}$, $l+1\le i\le l+u$, or  
        \item $s_{k+1}=C_{p}$,
    \end{itemize}  then $E_{k}(t)\ge E(s_{k})$ on $[s_{k},s_{k+1})$, $E(s_{k+1})=E(s_{k})$, and we denote this as $(s_{k},s_{k+1})\in \Xi_{1}$;

    2. If $s_{k}=C_{p}$ for some $1\le p\le l$ and $s_{k+1}+1$ is a first or free jump time, then $E_{k}(t)\ge E(s_{k})$ on $[s_{k},s_{k+1})$, and we denote this as $(s_{k},s_{k+1})\in \Xi_{2}$;

    3. If $s_{k}+1$ is a jump time and $s_{k+1}=C_{p}$  for some $1\le p\le l$, then $E_{k}(t)\ge E(s_{k+1})$, and we denote this as $(s_{k},s_{k+1})\in \Xi_{3}$;

    4. We write $(s_{k},s_{k+1})\in \Xi_{4}$ if none of $\Xi_{0},\ldots,\Xi_{3}$ occurs.

\end{itemize}

Therefore, by Eqn.~(\ref{eq_weight}),
\begin{equation}\label{eq_discretefunctional}
    \begin{split}
       N(N-1)&\cdots(N-l+1)\sum_{(E,\vec{q})\in \tilde{\mathcal{B}}_{\epsilon}(1,\ldots,l)\cap\mathcal{C}}w(E,\vec{q})
       =N(N-1)\cdots(N-l+1)\sum_{\vec{\delta}}\sum_{\vec{\nu}}\sum_{\vec{q}}\sum_{\vec{H}}w(E,\vec{q})\\
       =N(N-1)&\cdots(N-l+1)\sum_{\vec{\delta}}(N-l)\cdots (N-l-u+1)\sum_{\vec{\nu}}\sum_{\vec{q}}\sum_{\vec{H}}\prod_{k=1}^{n-1}\frac{\prod_{t\in \Delta}(-1)^{\mathbf{1}[q(t)<q(t-)]}}{N^{\delta}}\\
       \times&\left(\sum_{E_{k}}\mathbf{1}[E_{k}(t)-q(t)\ge 0]\frac{w(E_{k})}{\mu_{+}(N)^{s_{k+1}-s_{k}}}\prod_{t\in \llbracket s_{k}, s_{k+1}-1\rrbracket,\ t+1\in D(E)}\left(\frac{m_{p}(t)}{N}+2\frac{E_{k}(t)-q(t)}{\beta N}\right)\right),
    \end{split}
\end{equation}
where the index $p$ in $m_{p}(t)$  gives the main variable when $t\in (C_{p-1},C_{p}]$, and the factor $(N-l)\cdots (N-l-u+1)$ comes from choosing\footnote{Note that $u$ is determined by $\vec{\delta}$.} $u$ auxiliary variables from $x_{l+1},\ldots,x_{N}$ and giving them a complete order among the $u!$ choices.

Next, we recall the moment expression of the Laplace transform  of $\mathrm{Airy}(\beta)$ from \cite[Section 2]{GXZ}. The purpose of this is to identify the limit in Theorem \ref{thm:mainconvergence} with the object introduced in the concurrent work, while the convergence itself will follow from the arguments presented in this text. The interested reader can find a detailed exposition there.  Below we outline the necessary notation and formulas, with some small adaptations to fit our needs.

To begin we define the continuous analogues of the objects we introduced previously. Let $\bk_{1},\ldots,\bk_{l}\in \R_{>0}$, $\mathbf{Q}_{p}=\sum_{j=1}^{p}\bk_{j}$ for $p=1,2,\ldots,l$.  Let $\vec{\bp}=\{\bp_{j}\}_{j=1}^{l+u}$ be the \emph{block process} where each $\bp_{j}:[0,\mathbf{Q}_{l}]\rightarrow\R_{\ge0}$, $\bp_{j}(t)\equiv 0$ on  $(\mathbf{Q}_{j-1},\mathbf{Q}_{l}]$ for $1\le j\le l$. Let $\bp^{0}(t)=\sum_{j=1}^{l+u}\bp_{j}(t)$. Order $\bp_{j}:\ l+1\le j\le l+u$ so that $\min\{t:\bp_{l+1}(t)>0\}<\cdots<\min\{t:\bp_{l+u}(t)>0\}$. Let $\vec{\bb}=(\bb_{1},\ldots,\bb_{l})\in (\R_{>0}\cup\{\emptyset\})^{l}$ be the \emph{virtual block process}. Let $\mathbf{\Delta}$ be the set of discontinuity points of $\vec{\bp}$, which locate in $\bigcup_{p=2}^{l}(\mathbf{Q_{p-1}},\mathbf{Q}_{p})$, $\bdel$ be the size of $\mathbf{\Delta}$, $\vec{\bdel}=(\bdel_{j,p}:\ 1\le j\le l+u, 1\le p\le l)$ be the tuple, where $\bdel_{j,p}$ gives the number of discontinuity points of $\bp_{j}$ in $(\mathbf{Q}_{p-1},\mathbf{Q}_{p})$, and $\mathbf{H}:\mathbf{\Delta}\cup \{\mathbf{Q}_{p}:0\le p\le l\}\cup \{\mathbf{Q}_{p-1}+\bb_{p}:2\le p\le l, \bb_{p}\ne\varnothing,\ 1\le p\le l\}\rightarrow \R_{\ge 0}$ be a collection of continuous heights. We do not list the detailed definitions of $\vec{\bp},\vec{\bb}$ and $\mathbf{H}$, but simply note that all the conditions are parallel to those of our discrete objects $\vec{q}$, $\vec{\delta}$ and $\vec{H}$. Given the data $\vec{\bp},\vec{\bb}$, let $\mathcal{H}[\vec{\bk}]$ be the space of the tuple ($\vec{\bp}, \vec{\bb}, \mathbf{H})$, and $d(\vec{\bp}, \vec{\bb}, \mathbf{H})$ be the Lebesgue measure on $\mathcal{H}[\vec{\bk}]$\footnote{For clarity, here we use a different notation that we separate the discrete sums over $\vec{\bdel}$ and the events $\{\nu_{p}=\varnothing\}$ with the continuous integral, while in \cite{GXZ}, $d(\vec{\bp}, \vec{\bb}, \mathbf{H})$ denotes the combination of both parts.}. Then by \cite[Definition 2.11 and 2.12]{GXZ},

\begin{equation}
\resizebox{0.9\textwidth}{!}{$
    \bL_\beta(\vec\bk, \vec\ttt=0) =\sum\limits_{\vec{\bdel}}\sum\limits_{\mathbf{1}[\bb_{p}=\varnothing]\in \{0,1\},\ 2\le p\le l} \mathrm{p.v.}\int\limits_{\mathcal{H}[\vec\bk]} 2^{-\bdel-|\{p\in\llbracket 1,l\rrbracket: \bb_p\neq\varnothing,\ 1\le p\le l\}|}\prod\limits_{k\ge 0:\ (y_{k},y_{k+1})\in\mathbf{\Xi}} \mathbf{I}_{\beta,1}[y_{k},y_{k+1}]d(\vec\bp,  \vec\bb,\bH),
    $}
\end{equation}
where
\begin{itemize}
    \item $0\le y_{0}\le y_{1}\le \ldots\le y_{\delta+l+\# \bb_{p}:\ \bb_{p}\ne\varnothing}$ are the ordered elements in $\mathbf{\Delta}\cup \{\mathbf{Q}_{0},\ldots,\mathbf{Q}_{l}\}\cup \{\mathbf{Q}_{p-1}+\bb_{p}:\bb_{p}\ne\varnothing,\ 2\le p\le l\}$.
    \item Depending on $(y_{k},y_{k+1})\in \mathbf{\Xi}_{1},\ldots,\mathbf{\Xi}_{4}$, for\footnote{The dependence on $r$ is not presented in \cite{GXZ}, since in their context $r$ is always 1.} $r>0$,
\[
\resizebox{0.9\textwidth}{!}{$
   \mathbf{I}_{\beta,r}[y_{k},y_{k+1}]= \E\left[\exp\left(\frac{r}{\beta}\int_{y_{k}}^{y_{k+1}}(W_{k}(t)+[\mathbf{H}(y_{k})-\bp^{0}(y_{k})]\cdot \mathbf{1}[(y_{k},y_{k+1})\in \mathbf{\Xi}_{1}\ \text{or}\ \mathbf{\Xi}_{2}]dt\right)\right] \mathbf{f}(\mathbf{H}(y_{k}),\mathbf{H}(y_{k+1})),
   $}
\]
where (recalling Definition \ref{def:densities})
\begin{align*}
    \mathbf{f}(\mathbf{H}(y_{k}),\mathbf{H}(y_{k+1}))=\begin{cases}
        \mathbf{F}_{0,0}(y_{k+1}-y_{k}),\quad & (y_{k},y_{k+1})\in \mathbf{\Xi}_{1},\\
        \mathbf{F}_{0}(y_{k+1}-y_{k};\mathbf{H}(y_{k+1})-\mathbf{H}(y_{k})),\quad & (y_{k},y_{k+1})\in \mathbf{\Xi}_{2}, \\
         \mathbf{F}_{0}(y_{k+1}-y_{k};\mathbf{H}(y_{k})-\mathbf{H}(y_{k+1})),\quad & (y_{k},y_{k+1})\in \mathbf{\Xi}_{3},\\
         \mathbf{F}(y_{k+1}-y_{k};\mathbf{H}(y_{k})-\bp^{0}(y_{k}),\mathbf{H}(y_{k+1})-\bp^{0}(y_{k+1})),\quad & (y_{k},y_{k+1})\in \mathbf{\Xi}_{4},     
    \end{cases}
\end{align*}
and $W_{k}(\cdot)$ is
\begin{itemize}
    \item $B_{e}:\ [0,y_{k+1}-y_{k}]\rightarrow \R_{\ge 0}$, a Brownian excursion, if $(y_{k},y_{k+1})\in \mathbf{\Xi}_{1}$;
    \item $B_{3}:\ [0,y_{k+1}-y_{k}]\rightarrow \R_{\ge 0}$, a Bessel$_{3}$ process starting at $0$, conditioned on ending at $\mathbf{H}(y_{k+1})-\mathbf{H}(y_{k})$, if $(y_{k},y_{k+1})\in \mathbf{\Xi}_{2}$;
    \item $B_{3}(y_{k+1}-y_{k}-\cdot):\ [0,y_{k+1}-y_{k}]\rightarrow \R_{\ge 0}$, a time-reversed Bessel$_{3}$ process starting at $0$, conditioned on ending at $\mathbf{H}(y_{k+1})-\mathbf{H}(y_{k})$, if $(y_{k},y_{k+1})\in \mathbf{\Xi}_{3}$;
    \item $B: [0,y_{k+1}-y_{k}]\rightarrow\R_{\ge 0}$, a Brownian bridge from $\mathbf{H}(y_{k})$ to $\mathbf{H}(y_{k+1})$, conditioned on $B\ge 0$, if $(y_{k},y_{k+1})\in \mathbf{\Xi}_{4}$.
\end{itemize}

    \item The principal value integral is taken in the following sense:
for any small enough $\epsilon>0$, we let  $\mathcal{H}_\epsilon=\mathcal{H}_\epsilon[\vec\bk]$ denote the space of all tuples $(\vec\bp, \vec\bb, \bH)\in \mathcal{H}[\vec \bk]$, such that for each $p\in\llbracket 1,l\rrbracket$, $\bp^0$ is constant in $(\bQ_{p-1}, \bQ_{p-1}+\epsilon)$, and either $\bb_p>\epsilon$, or $\bb_p=\varnothing$.
We integrate over $\mathcal{H}_\epsilon$, then send $\epsilon\to 0+$.

In general $\bL_\beta(\vec\bk, \vec\ttt)$ gives the mixed moments of the Laplace transform of Airy$_\beta$ line ensemble at time $\ttt_{1},\ldots,\ttt_{l}$. Since $\mathrm{Airy}(\beta)$  is the one-time marginal of Airy$_{\beta}$ line ensemble, which is shift invariant in distribution, in our context we take $\ttt_{1}=\cdots=\ttt_{l}=0$.
\end{itemize}

\begin{thm}\label{thm:mainconvergence} 
Fix $\bk_{1},\ldots,\bk_{l}>0$. Let $M_{1},\ldots,M_{l}$ be sequences of positive integers depending on $N$ such that \[\left|M_{p}-\bk_{p}N^{2/3}\right|<C\] for some constant $C>0$ uniform in $N$ and $p=1,2,\ldots,l$. Take the sequence of ensembles in Eqn.~(\ref{eq_betaaddition}) satisfying the assumptions in Theorem \ref{thm:main}, with limiting parameters $P_{-1}$, $\sigma$. Then as $N\rightarrow\infty$,
 \begin{align*}
&\E\left[\prod_{p=1}^{l}\left(\sum_{i=1}^{N}\left(\frac{\lambda_{i}}{\mu_{+}(N)N}\right)^{M_{p}}\right)\right]
\rightarrow \bL_\beta(\left(2P_{-1}\sigma\right)^{-2/3}\vec\bk, \vec\ttt=0) .
 \end{align*}   
\end{thm}
\begin{proof}
    It remains to check that Eqn.~(\ref{eq_discretefunctional}) converges to the right-hand side above as $N\to\infty$ and $\epsilon\to0$. By Lemma \ref{lem:principalvalue}, we take $(E,\vec{q})\in \tilde{\mathcal{B}}_{\epsilon}(1,2,\ldots,l)\cap \mathcal{C}$ with fixed data $\vec{\delta}$ and $\mathbf{1}[\nu_{p}=\varnothing]$ for $2\le p\le l$. The principal value integral will be given by taking $\epsilon\to 0$. 
    
    Rescale each configuration $(E,\vec{q})$ horizontally by $ N^{2/3}$ and vertically by $\sigma N^{1/3}$, so that $\bdel=\delta$, $\mathbf{\Delta}=N^{-2/3}\Delta$, and
    \begin{itemize}
        \item For $0\le k\le n-1$, $E_{k}(t)= \sigma N^{1/3} W_{k}( tN^{2/3})$.
        \item For $1\le i\le l$, $q_{i}(t)= \sigma N^{1/3} \bp_{i}( tN^{2/3})$ where $t\in [0,C_{i-1}]$. For $l+1\le i\le l+u$, $q_{i}(t)=\sigma N^{1/3} \bp_{i}( tN^{2/3})$ where  $t\in [0,C_{l}]$. For $2\le p\le l$, $\nu_{p}=N^{2/3}\bb_{p}$ if $\nu_{p}\ne \varnothing$, and $\bb_{p}=\varnothing$ if $\nu_{p}=\varnothing$.  
        \item For $0\le k\le n$, $s_{k}=N^{2/3}y_{k}$, and $H(s_{k})=\sigma N^{1/3}\mathbf{H}(y_{k})$.
    \end{itemize}
    Then \begin{equation}\label{eq_minussigns}
    \begin{split}
       &N(N-1)\cdots (N-l+1)(N-l)\cdots (N-l-u+1) \frac{\prod_{t\in \Delta}(-1)^{\mathbf{1}[q(t)<q(t-)]}}{N^\delta}\\
       =&\frac{N^{-\delta+l+u}}{u!}(1+o(1))\prod_{(y_{k},y_{k+1})\in \Xi_{3}}(-1)^{\mathbf{1}[\bp^{0}(y_{k})<\bp^{0}(y_{k}-)]}\prod_{(y_{k},y_{k+1})\in \Xi_{4}}(-1)^{\mathbf{1}[\bp^{0}(y_{k})<\bp^{0}(y_{k}-)]}.
       \end{split}
    \end{equation}
    The mesh size tends to 0, and we can treat the sum over $\vec{q}$ and over those $\nu_{p}$ with $\nu_{p}\ne \varnothing$ as a Riemann sum, with an additional factor $N^{\frac{2}{3}(\delta+\nu')}(\sigma N^{\frac{1}{3}})^{\delta-u}$, where $$\nu':=|\{p\in\llbracket 2,l\rrbracket: \nu_p\neq\varnothing\}|=|\{p\in\llbracket 2,l\rrbracket: \bb_p\neq\varnothing\}|.$$ 
    In particular, the factor $\frac{1}{u!}$ gives an ordering  $\min\{t:\bp_{l+1}(t)>0\}<\dots<\min\{t:\bp_{l+u}(t)>0\}$. Similarly, the value of $H(s_{k})$ is not determined by $q(t)$ if and only if $s_{k}$ is the jump time of a first/free jump, so this gives a factor $(\sigma N^{\frac{1}{3}})^{\delta-u}$. These together give 
    \begin{equation}\label{eq_factor1}
        \sigma^{2(\delta-u)}N^{\frac{4}{3}\delta+\frac{2}{3}\nu'-\frac{2}{3}u}.
    \end{equation}

    Denote the partition function of all the possible walks $E_{k}$ on $[s_{k},s_{k+1}]$ by $Par(k)=\sum_{E_{k}} w(E_{k})$. For each $0\le k\le n-1$, by Lemma \ref{lem:conditionalwalk2},
    \begin{equation}\label{eq_transitiondensity}
    \begin{split}
        \frac{Par(k)}{\mu_{+}(N)^{s_{k+1}-s_{k}}} = \begin{cases}
        \frac{1}{\mu_{+}(N)}\rightarrow P_{-1}, \quad &(s_{k},s_{k+1})\in \Xi_{0};\\
        2^{-1}(P_{-1}\sigma)^{-1}N^{-1}    [\mathbf{F}_{0,0}(y_{k+1}-y_{k})+o(1)], \quad &(s_{k},s_{k+1})\in \Xi_{1};\\
          \sigma^{-2}N^{-2/3}  [\mathbf{F}_{0}(y_{k+1}-y_{k};\mathbf{H}(y_{k+1})-\mathbf{H}(y_{k}))+o(1)], \quad& (s_{k},s_{k+1})\in \Xi_{2};\\
           2^{-1}P_{-1}^{-1}N^{-2/3} [ \mathbf{F}_{0}(y_{k+1}-y_{k};\mathbf{H}(y_{k})-\mathbf{H}(y_{k+1}))+o(1)], \quad& (s_{k},s_{k+1})\in \Xi_{3};\\
             \sigma^{-1}N^{-1/3}[\mathbf{F}(y_{k+1}-y_{k};\mathbf{H}(y_{k}),\mathbf{H}(y_{k+1}))+o(1)],\quad &(s_{k},s_{k+1})\in \Xi_{4},\\
        \end{cases}
        \end{split}
    \end{equation}
    where $\mathbf{F}_{0,0}(x)$, $\mathbf{F}_{0}(x;h)$ and $\mathbf{F}(x;h,g)$ are defined as in Definition \ref{def:densities}.
    Denote the number of $(s_{k},s_{k+1})\in \Xi_{j}$ by $\mathcal{N}_{j} (j=0,..,4)$, then 
    $$\mathcal{N}_{0}=\delta+\nu',\quad\mathcal{N}_{2}=l-\mathcal{N}_{1},\quad \mathcal{N}_{3}=u+l+\nu'-\mathcal{N}_{1},\quad \mathcal{N}_{4}=\delta+l-(u+l)-l+\mathcal{N}_{3}$$
    give the multiplicities of the terms in Eqn.~(\ref{eq_transitiondensity}). In particular, the factor including $P_{-1}$, $\sigma$ and $N$ is 
    \begin{equation}\label{eq_factor2}
    \begin{aligned}
         P_{-1}^{\mathcal{N}_{0}}[(P_{-1}\sigma)^{-1}N^{-1}]^{\mathcal{N}_{1}}[ \sigma^{-2}N^{-2/3}]^{\mathcal{N}_{2}}[ P_{-1}^{-1}N^{-2/3} ]^{\mathcal{N}_{3}}[\sigma^{-1}N^{-1/3}]^{\mathcal{N}_{4}}\\
         \qquad =2^{-l-u-\nu'}P_{-1}^{-l-u+\delta}\sigma^{-l+u-\delta}N^{-\frac{1}{3}\delta-\frac{1}{3}u-\frac{2}{3}\nu'-l}.
    \end{aligned}
    \end{equation}
     
    Note that for all $1\le i\le l$ and $1\le t\le M$, $N-l-u\le m_{i}(t)\le N-1$. By Theorems \ref{thm:brownian1} and \ref{thm:brownian2}, each $E_{k}\in \Xi_{j}$, $j=1,\ldots,4$ converges weakly to a conditional Brownian bridge $W_{k}$. Moreover, by Lemma \ref{lem:timehom}, with high probability the down-steps $D(E)$ distribute homogeneously in each $E_{k}$ with proportion tending to $P_{-1}$, and Lemma \ref{lem:walkmax} gives the desired tightness of the walk bridges, which allows us to upgrade the weak convergence to the moment convergence of the functionals. Therefore,
    \begin{equation}\label{eq_brownianfunctional}
    \begin{split}
        &\left(\sum_{E_{k}}\mathbf{1}[E_{k}(t)-q(t)\ge 0]\frac{w(E_{k})}{\mu_{+}(N)^{s_{k+1}-s_{k}}}\prod_{t\in \llbracket s_{k}, s_{k+1}-1\rrbracket,\ t+1\in D(E)}\left(\frac{m_{i}(t)}{N}+2\frac{E_{k}(t)-q(t)}{\beta N}\right)\right)\\
        \rightarrow&\E\left[\exp\left(\frac{2P_{-1}\sigma}{\beta}\int_{y_{k}}^{y_{k+1}}(W_{k}(t)+[\mathbf{H}(y_{k})-\bp^{0}(y_{k})]\cdot \mathbf{1}[(y_{k},y_{k+1})\in \mathbf{\Xi}_{1}\ \text{or}\ \mathbf{\Xi}_{2}]dt\right)\right].
    \end{split}
    \end{equation}

Combining Eqn.~(\ref{eq_minussigns}), (\ref{eq_factor1}), (\ref{eq_transitiondensity}), (\ref{eq_factor2}) and (\ref{eq_brownianfunctional}), then take $\epsilon\to 0$ turns Eqn.~(\ref{eq_discretefunctional}) into 
\begin{equation}\label{eq_limitingexpression2}
\resizebox{0.9\textwidth}{!}{$
    \sum\limits_{\vec{\bdel}}\sum\limits_{\mathbf{1}[\bb_{p}=\varnothing]\in \{0,1\},\ 2\le p\le l}  \mathrm{p.v.}\int\limits_{\mathcal{H}[\vec\bk]} 2^{-(l+u+|\{p\in\llbracket 2,l\rrbracket: \bb_p\neq\varnothing\}|)}(P_{-1}\sigma)^{\bdel-l-u}\prod\limits_{k\ge 0:\ (y_{k},y_{k+1})\in\mathbf{\Xi}} \mathbf{I}_{\beta,2P_{-1}\sigma}[y_{k},y_{k+1}]d(\vec\bp,  \vec\bb,\bH).
$}
\end{equation}
To identify the above with $\bL_\beta(\left(2P_{-1}\sigma\right)^{-2/3}\vec\bk, \vec\ttt=0)$, first note that when\footnote{This is the case for a single G$\beta$E, considered in \cite{GXZ}.} $2P_{-1}\sigma=1$, this holds directly. Therefore it remains to match the exponent of the factors $(P_{-1}\sigma)$ and $\bk_{p}$. For the Brownian bridges $W_{k}(t)$ satisfy
$\bk^{1/2}W_{k}(t)\stackrel{d}{=}W_{k}(\bk t)$. Therefore in Eqn.~(\ref{eq_brownianfunctional})
 (or equivalently in $\prod_{k\ge 0:\ (y_{k},y_{k+1})\in\mathbf{\Xi}} \mathbf{I}_\beta[y_{k},y_{k+1}]$ above), for $[y_{k},y_{k+1}]\subset [\mathbf{Q}_{j-1},\mathbf{Q}_{j}]$,
\[\frac{2P_{-1}\sigma}{\beta}\int_{y_{k}}^{y_{k+1}}(W_{k}(t)+[\mathbf{H}(y_{k})-\bp^{0}(y_{k})]\cdot \mathbf{1}[(y_{k},y_{k+1})\in \mathbf{\Xi}_{1}\ \text{or}\ \mathbf{\Xi}_{2}]dt\sim (P_{-1}\sigma)\cdot \bk_{j}^{3/2} .\]
Similarly, \[\mathbf{F}_{0,0}(\bk y)\sim \bk^{-3/2},\]
so starting from the configurations where $\delta=0$, so that $W_{1}(t),\ldots,W_{l}(t)$ are $l$ Brownian excursions, Eqn.~(\ref{eq_transitiondensity}) gives a $\vec{\bk}$-dependent factor $\prod_{p=1}^{l}\bk_{p}^{-3/2}$, and the exponent of $P_{-1}\sigma $ in Eqn.~(\ref{eq_limitingexpression2}) is $\delta-l-u=-l$, so one copy of $\bk_{j}^{3/2}$ is paired with one copy of $(P_{-1}\sigma)$. When $\delta>0$, the measure $d(\vec\bp,  \vec\bb,\bH)$ gives a $\vec\bk$-dependent factor\footnote{The first term in the exponent counts the number of jumps in $[C_{p-1},C_{p})$, and the second term counts the number of auxiliary variables whose last jump time is in $[C_{p-1},C_{p})$.} \[\prod_{p=1}^{l}\bk_{p}^{\frac{3}{2}(\sum_{j=1}^{l+u}\bdel_{j,p}-\#\{l+1\le i\le l+u:\ \bdel_{j,p}>0,\ \bdel_{j,p+1}=\cdots=\bdel_{j,l}=0 \}) },\]
and the sum of the exponents above is $\delta-u$, so again one copy of $\bk_{p}^{3/2}$ is paired with one copy of $(P_{-1}\sigma)$. This finishes the proof.
\end{proof}

\section{Proof of the main theorem}\label{sec:momenttoweak}
In this section, we finish the proof of Theorem \ref{thm:main}, based on Theorem \ref{thm:mainconvergence}. We follow a three-step strategy: in Section \ref{sec:laplaceconv}, we identify the limiting Laplace transform expression of the values $\lambda_{i}'(N)$ along a subsequence; in Section \ref{sec:universality}, we identify the Laplace transform expression with that of $\mathrm{Airy}(\beta)$; then in Section \ref{sec:uniqueness}, we prove the uniqueness of the underlying point process.  
\subsection{Limit Laplace transform}\label{sec:laplaceconv}
The main result in this section is the following lemma.
\begin{lem}\label{lem:laplace}
 There exists an increasing sequence of positive integers $N_{1}<N_{2}<\cdots$, such that as $k\to\infty$, $\{\lambda_{i}'(N_{k})/\mu_{+}(N_{k})\}_{i=1}^{N_{k}}$ converges weakly to a point process $\eta_{1}\ge\eta_{2}\ge \cdots$ in $\R\cup \{-\infty\}$, and for any $\bk_{1},\ldots,\bk_{l}>0$,
\begin{equation}\label{eq_laplacefinal}
    \E\left[\prod_{j=1}^{l}\left(\sum_{i=1}^{\infty}\exp\left(\frac{\bk_{j}\eta_{i}}{\mu_{+}}\right)\right)\right]=\bL_\beta(\left(2P_{-1}\sigma\right)^{-2/3}\vec\bk, \vec\ttt=0).
\end{equation}    
\end{lem}

As preparation for the proof of Lemma \ref{lem:laplace}, we present a deterministic statement. Consider the space $\R\cup\{-\infty\}$, whose topology is generated by open intervals $(a, b)$ and $[-\infty, a)$ for all $a<b\in\R$.

For any $x\in\R$, $\alpha>0$, $N\in\N$, we denote
\[
M_+[x,\alpha,N]= \frac{1}{2}(1+N^{-2/3}x)^{2\lfloor \alpha N^{2/3}/2\rfloor},\quad
M_-[x,\alpha,N]= \frac{1}{2}(1+N^{-2/3}x)^{2\lfloor \alpha N^{2/3}/2\rfloor+1},
\]
and $M[x,\alpha,N]=M_+[x,\alpha,N]+M_-[x,\alpha,N]$.

\begin{prop}  \label{prop:analytic}
For integers $1\le i\le N$, let $x^N_i \in \R$ satisfy $x^N_i\ge x^N_{i+1}$ for each $1\le i \le N-1$.
Suppose there exists an increasing sequence of positive integers $N_1<N_2<\cdots$, such that for any $\alpha\in\N$ or $\alpha^{-1}\in\N$, as $k\to\infty$, the sum $\sum_{i=1}^{N_k} M[x_i^{N_k},\alpha,N_k]$ converges to a limit in $\R$, and $\sum_{i=1}^{N_k} M_+[x_i^{N_k},\alpha,N_k]$ is bounded by a constant independent of $k$ (but might depend on $\alpha$).
Then for each $i\in\N$, as $k\rightarrow\infty$, $x_i^{N_k}$ converges in $\R\cup\{-\infty\}$ to some limit $x_{i}$, and 
for any $\alpha>0$,
$$\sum_{i=1}^{N_k} M[x_i^{N_k},\alpha,N_k]\to \sum_{i=1}^\infty \exp(\alpha x_i).$$
\end{prop}
Proposition \ref{prop:analytic} can be viewed as a special case of \cite[Proposition 5.5]{GXZ}, and for completeness we present a proof in  Appendix B.
\vspace{0.2cm}

\noindent{\textbf{Proof of Lemma \ref{lem:laplace}.}} As a special case of Theorem \ref{thm:mainconvergence}, for any $\bk>0$, we have the convergence of 
\[
\E\left[\left(\sum_{i=1}^{N}M[\lambda_{i}'(N)/\mu_{+}(N),\bk,N]\right)^{2}\right] \quad \text{and} \quad \E\left[\left(\sum_{i=1}^{N}M_{+}[\lambda_{i}'(N)/\mu_{+}(N),\bk,N]\right)^{2}\right],
\]
which implies the tightness of 
\[
\sum_{i=1}^{N}M[\lambda_{i}'(N)/\mu_{+}(N),\bk,N] \quad \text{and} \sum_{i=1}^{N}M_{+}[\lambda_{i}'(N)/\mu_{+}(N),\bk,N].
\]
Therefore, there exists a subsequence $\{N_{k}\}_{k=1,2,\ldots}$ such that for all $\bk\in \N$ and $\bk^{-1}\in \N$, both 
\[
\sum_{i=1}^{N_{k}}M[\lambda_{i}'(N_{k})/\mu_{+}(N_{k}),\bk,N] \quad \text{and} \quad \sum_{i=1}^{N_{k}}M_{+}[\lambda_{i}'(N_{k})/\mu_{+}(N_{k}),\bk,N]
\]
converge weakly as $k\to\infty$. 

By the Skorokhod representation theorem, there exists a probability space on which the convergence holds almost surely, so 
\[
\sum_{i=1}^{N_{k}}M[\lambda_{i}'(N_{k})/\mu_{+}(N_{k}),\bk,N] \quad \text{and} \quad \sum_{i=1}^{N_{k}}M_{+}[\lambda_{i}'(N_{k})/\mu_{+}(N_{k}),\bk,N]
\]
satisfy the conditions in Proposition \ref{prop:analytic}. Then for each $i\in \N$, $\lambda_{i}'(N_{k})/\mu_{+}(N_{k})$ converges to a limit $\eta_{i}/\mu_{+}$ as $\mu_{+}(N_{k})\rightarrow\mu_{+}$, and for any $\bk_{1},\ldots,\bk_{l}>0$, as $k\to\infty$ \[\prod_{j=1}^{l}\left(\sum_{i=1}^{N_{k}}M[\lambda_{i}'(N_{k})/\mu_{+}(N_{k}),\bk_{j},N]\right)\rightarrow\prod_{j=1}^{l}\left(\sum_{i=1}^{\infty}\exp\left(\frac{\bk_{j}\eta_{i}}{\mu_{+}}\right)\right)\quad \text{almost surely}.\]
Again by Theorem \ref{thm:mainconvergence}, the left-hand side above is uniformly integrable, so the expectation of the right-hand side is equal to the limit expectation of the left-hand side, which is $\bL_\beta(\left(2P_{-1}\sigma\right)^{-2/3}\vec\bk, \vec\ttt=0)$ by Theorem \ref{thm:mainconvergence}.\qed

\subsection{Matching with Airy($\beta$)}\label{sec:universality}

\begin{thm}\label{thm:matchlaplace}
$\bL_\beta(\left(2P_{-1}\sigma\right)^{-2/3}\vec\bk, \vec\ttt=0) $ is equal to the mixed moment of Laplace transform of $\tilde{C}\cdot \mathrm{Airy}(\beta)$, where $\tilde{C}$ is defined in Proposition \ref{prop:main}.
\end{thm}
\begin{proof}
By Lemma \ref{lem:laplace}, $\bL_\beta(\left(2P_{-1}\sigma\right)^{-2/3}\vec\bk, \vec\ttt=0) $ is the Laplace transform expression of the limit point process $\{\eta_{i}\}_{i=1}^{\infty}$ in $\R\cup \{-\infty\}$. Since this limiting expression is universal up to a constant, and our ensembles include a single copy of L$\beta$E as a special case, by the classical result \cite[Theorem 1.4]{RRV}, $\bL_\beta(\left(2P_{-1}\sigma\right)^{-2/3}\vec\bk, \vec\ttt=0) $ is the Laplace transform expression of the $\mathrm{Airy}(\beta)$ process, rescaled by some constant. 

It remains to identify this constant with $\tilde{C}$.
For a single $N\times N$ $L\beta E$, $\kappa_{l}=1$ for all $l=1,2,\ldots$. So by  Eqn.~(\ref{eq_c0}) and (\ref{eq_sigma}) we have $\tilde{C}=2^{4/3}$, $\mu_{+}=4$, $\sigma=\sqrt{2}$, $P_{-1}=\frac{1}{2}$. 
Therefore, the result holds if and only if, by comparison with this special case, the Laplace transform expression is unchanged after replacing $\bk_{1},\ldots,\bk_{l}$ by $\frac{\mu_{+}}{4\tilde{C}}\bk_{1}$,\ldots,$\frac{\mu_{+}}{4\tilde{C}}\bk_{l}$ for a general ensemble. 
Again by Eqn.~(\ref{eq_c0}) and (\ref{eq_sigma}) we have
\[
 (2P_{-1}\sigma)^{2/3} \left(\frac{\mu_+}{4\tilde{C}}\right) = (2\times\frac{1}{2}\times \sqrt{2})^{2/3}\times 2^{-4/3},
\]
so this is indeed the case.

\end{proof}

\subsection{Tightness and uniqueness of the limiting point process}\label{sec:uniqueness}
Denote the random point configuration of $\{\frac{1}{\tilde{C}}\lambda_{i}'(N)\}_{i=1}^{N}$ by $\mathcal{X}_{N}$, a sequence of point processes on $\R$. By Lemma \ref{lem:laplace} and Theorem \ref{thm:matchlaplace}, the Laplace transform of $\mathcal{X}_{N}$ converges in moments to the Laplace transform of $\mathrm{Airy}(\beta)$. 

We claim that $\{\mathcal{X}_{N}\}_{N\in \Z_{\ge 100}}$ is tight in the $W^{\#}$ topology of boundedly finite Borel measures on $\R$, see \cite[A2.6]{DVe1} for 
the definition. By \cite[Proposition 11.1.VI]{DVe2}, this is equivalent to checking that for any compact interval $A\subset \R$, the random variables 
$$N_{A}:= \text{number\ of\ particles\ of\ } \mathcal{X}_{N}\ \text{in\ A}$$
(we omit the dependence on $N$ in the notation for simplicity) are uniformly tight in $N$. It suffices to check $A=[t,\infty)$ for any $t\in \R$ since such sets contain all compact intervals. We have 
$$\frac{N_{A}}{2}\exp\left(\frac{Tt}{\mu_{+}}\right)\le \sum_{i}\left(1+\frac{\lambda'_{i}(N)}{\mu_{+}(N)N^{\frac{2}{3}}}\right)^{\lfloor TN^{\frac{2}{3}}\rfloor} $$
when $N$ is large. The expectation of the right-hand side converges as $N\rightarrow\infty$ by Theorem \ref{thm:mainconvergence}. Hence the uniform tightness follows.

It remains to check that the moments of the Laplace transform of $\mathrm{Airy}(\beta)$ uniquely determine a point process on $\R$. If this holds, then
by the tightness of $\{\mathcal{X}_{N}\}$, every subsequence of $\{\mathcal{X}_{N}\}$ has a further subsequence that converges to $\mathrm{Airy}(\beta)$, and so does the whole sequence.

By definition of correlation measures $\rho_{l}$ of point processes (see \cite{DVe1}, \cite{DVe2}), for a point process $\mathcal{X}$ in $\R$ with particles $\Lambda_{1}\ge\Lambda_{2}\ge\ldots$
\begin{equation}\label{eq_laplace1}
\E\left[\prod_{j=1}^{l}\left(\sum_{i}\exp\left(\bk_{j}\Lambda_{i}\right)\right)\right]=\int\cdots\int \rho_{l}(dx_{1},\ldots,dx_{l})e^{\bk_{1}x_{1}+\ldots+\bk_{l}x_{l}}+\text{lower\ order\ terms},
\end{equation}
where the lower order terms give a linear combination of Laplace transforms of $\rho_{1},\ldots,\rho_{l-1}$ of similar form. 
\begin{lem}
    All $\rho_{l}$ are uniquely determined by the Laplace transform above.
\end{lem}
\begin{proof}
    By iterating from $l=1$ in Eqn.~(\ref{eq_laplace1}), it remains to show that for any $l=1,2,\ldots$, $\bk_{1},\ldots,\bk_{l}>0$, the Laplace transform 
    \begin{equation}\label{eq_laplace2}
    \int\cdots\int \rho_{l}(dx_{1},\ldots,dx_{l})e^{\bk_{1}x_{1}+\ldots+\bk_{l}x_{l}}
    \end{equation}
      uniquely determines $\rho_{l}$. It is clear that Eqn.~(\ref{eq_laplace2}) is finite for all $(\bk_{1},\ldots,\bk_{l})\in \R_{>0}^{l}$, so fixing some $(S_{1},\ldots,S_{l})\in \R_{>0}^{l}$, we have $\rho_{l}(dx_{1},\ldots,dx_{l})e^{S_{1}x_{1}+\ldots+S_{l}x_{l}}$ is in $L^{1}$, i.e., a finite measure, and its Fourier transform is fixed by Eqn.~(\ref{eq_laplace2}) under analytic continuation. Then $\rho_{l}(dx_{1},\ldots,dx_{l})e^{S_{1}x_{1}+\ldots+S_{l}x_{l}}$ and hence $\rho_{l}(dx_{1},\ldots,dx_{l})$ is also fixed.
\end{proof}
Recall that by definition, the distribution of $\mathcal{X}$ is determined by finite-dimensional distributions of $N_{A_{1}},\ldots,N_{A_{k}}$, for any compact sets $A_{1},\ldots,A_{k}$, and 
\begin{equation}
    \E\left[N_{A}(N_{A}-1)\cdots (N_{A}-k+1)\right]=\int_{A}\cdots\int_{A}\rho_{k}(dx_{1},\ldots,dx_{k}).
\end{equation}
So the (joint) moments of the observables $N_{A}$ are fixed by correlation measures. It remains to check that for any compact set $A$, the tail probability of $N_{A}$ is exponentially decaying, so the moments of $N_{A}$ determine the distribution of $N_{A}$. We check this for $A=(\lambda,\infty)$ for any $\lambda\in\R$, which will suffice.
\vspace{0.5cm}

Since $\mathcal{X}=\mathrm{Airy}(\beta)$, we use the following result from \cite[Theorem 1.2]{RRV}: for $\lambda\in \R$,
\begin{equation}\label{eq_rrvthm}
    \PP[N_{(\lambda,+\infty)}\ge k]=\PP[\eta_{k}>\lambda]=\int_{\R^{k}}\kappa(\lambda,dx_{1})\kappa(x_{1},dx_{2})\cdots \kappa(x_{k-1},dx_{k}),
\end{equation}
where $\kappa(x,\cdot)$ is the distribution of the first passage time to $-\infty$ of diffusion $p(t)$ ($B(t)$ is a standard Brownian motion)
\begin{equation}\label{eq_airydiffusion}
    dp(t)=\left(t-p(t)^{2}\right)dt+\frac{2}{\sqrt{\beta}}dB(t),
\end{equation}
that starts from $+\infty$ at time $x$.

Then Eqn.~(\ref{eq_rrvthm}) can be rewritten as 
$$\PP[N_{(\lambda,+\infty)}\ge k]=\PP[p(t)\ \text{reaches}\ -\infty\ \text{for\ at\ least}\ k\ \text{times\ when\ starting\ from\ time}\ \lambda].$$

We collect the facts about $\kappa(x,\cdot)$ we need in the next lemma.
\begin{lem}
    \noindent{(1).} $\kappa(x,y)=0$ for $y<x$.

    \noindent{(2).} For any $x\in \R$, $\kappa(x,\{+\infty\})>0$.

    \noindent{(3).} $\kappa(x,\{+\infty\})$ increases with $x$.
\end{lem}
\begin{proof}
    (1) is clear from the definition. To see (2) we set $k=1$ in Eqn.~(\ref{eq_rrvthm}) and get \begin{equation}\label{eq_explosiontime}
        \kappa(x,\{+\infty\})=1-\PP[\eta_{1}>x]=\PP[\eta_{1}\le x],
    \end{equation}
    which is clearly in $(0,1)$ for all $x\in \R$ (see \cite[Theorem 1.3]{RRV} for a tail estimate of $\eta_{1}$). (3) also follows directly from Eqn.~(\ref{eq_explosiontime}).
\end{proof}

It is clear from the definition that each run of $p(t)$ after the restart from $+\infty$ is independent. Hence by the above properties of transition measure $\kappa(x,\cdot)$, we have 
$$\PP[p(x)\ \text{reaches}\ -\infty\ \text{for\ at\ least}\ k+1\ \text{times\ when\ starting\ from\ time}\ \lambda]\le (1-\kappa(\lambda,\{+\infty\}))^{k+1},$$
which is exponentially decaying. This finishes the proof.\qed

\section{Further studies}\label{sec:further}
\subsection{Negative free cumulants}\label{sec:negativecumulants}
Recall that in this text we assume that for our $\beta$-additions, the  free cumulants $\kappa_{1}(N),\kappa_{2}(N),\ldots$ are all nonnegative, see Eqn.~(\ref{eq_cumulant}). This is necessary because the proofs presented from Section \ref{sec:technical} to Section \ref{sec:asymptotics} use a probabilistic point of view, that is, we interpret the huge sum arising from the actions of Dunkl operators on $G(z_{1},\ldots,z_{N};\beta)$ as a certain functional of a conditional random walk $E(t)$. The increments $X_{i}(N)$ of $E(t)$ are i.i.d. before conditioning and, as $N\rightarrow\infty$,
$$X_{1}(N)\implies X_{1},$$
whose distribution function is given in terms of the free cumulants $\kappa_{1},\kappa_{2},\ldots$ and the Voiculescu transform $V(z)$.

Recall also that this assumption includes all cases in which we add Gaussian ensembles with finitely many Laguerre ensembles, and a subset of the cases in which we have subtractions of Laguerre ensembles. We naturally conjecture that if we assume similar convergence for each summand as stated in Section \ref{sec:mainresult}, the edge limit of the additions will still be the $\mathrm{Airy}(\beta)$ process under proper rescaling, even when some $\kappa_{l}$ are negative.

It seems that our current proof cannot be directly applied to this case. For simplicity, let us only consider the asymptotics of \[\E\left[\sum_{i=1}^{N}\left(\frac{\lambda_{i}}{\mu_{+}(N)N}\right)^{M}\right],\] where $M\approx \bk N^{2/3}$. There are two main ingredients of our proof that do not work as before: the steepest descent calculation and the functional CLT in Section \ref{sec:technical}. 

Consider $\kappa_{l}(N)$ and $V_{N}(z)$ defined as in Eqn.~(\ref{eq_k2}), (\ref{eq_kl}) and (\ref{eq_vz}), and take
$\alpha_{1}>0$, $\alpha_{2}<0$.
For the steepest descent of such $V_{N}(z)$, we claim without proof that there is a unique root of $V^{'}_{N}(z)$ in $(-\frac{1}{\alpha_{2}},0)$ and another in $(0,\frac{1}{\alpha_{1}})$, which we denote by $z_{c,1}(N)$ and $z_{c,2}(N)$, respectively. By choosing a contour $C_{z}'(N)$ symmetric about the real line that crosses both $z_{c,1}(N)$ and $z_{c,2}(N)$ and only surrounds the pole $z=0$ of $V_{N}(z)$,  one can check that
\begin{equation}
|m_{M}(N)|=\frac{1}{M+1}\left|\frac{1}{2\pi i}\oint_{C_{z}'(N)}V_{N}(z)^{M+1}dz\right|\sim \frac{1}{\sqrt{M}}|V_{N}(z_{c}(N))|^{M+1},
\end{equation}
where $z_{c}(N)$ is whichever of $z_{c,1}(N)$ or $z_{c,2}(N)$ that maximizes $|V_{N}(z)|$ on $C_{z}'(N)$. Without loss of generality, suppose $V_{N}(z_{c}(N))>0$. Since the edge behavior is characterized by $m_{M}(N)$ for $M\sim N^{\frac{2}{3}}$, we believe that the critical point $z_{c}(N)$ in the steepest descent gives the upper edge $\mu_{+}(N)=V(z_{c}(N))$ as before.

\vspace{0.2cm}
\noindent{\textbf{Question:}} Is it always the case that $|z_{c,1}(N)|<\frac{1}{\alpha_{1}}$, so that $V_{N}(z_{c,1}(N))$ converges absolutely as $N\to \infty$? 
\vspace{0.2cm}

Issues also arise in the use of the functional CLT. If we take the same interpretation of conditional random walk with i.i.d. steps $X_{i}(N)$ as before but allow some cumulants to be negative, then $X_{i}(N)$ is a formal integer-valued random variable, such that for $l\ge 0$,
$$\PP(X_{i}(N)=l)=\frac{\kappa_{l+1}(N)z_{c}(N)^{l}}{V_{N}(z_{c}(N))},$$
which is negative for some $l$. Everything can be understood in a measure-theoretic sense by replacing the probability measure with a signed measure of total mass 1. However, such formal random objects behave poorly from a probabilistic point of view, with many elementary arguments breaking down. In particular, in Section \ref{sec:technical}, $P_{-1}=\PP(X_{i}=-1)$ gives the asymptotic proportion of down-steps, but our arguments no longer hold when $X_{i}$ is not supported by a probability measure, because even the classical weak LLN no longer holds. In addition, it is not clear to us whether one can state the analog for signed measures of Theorems \ref{thm:brownian1} and \ref{thm:brownian2} that deal with the weak convergence of conditional random walks to conditional Brownian motions. Another obstruction is showing that the tail probability of the random walk decays exponentially, as done in Section \ref{sec:tailrw}. We give a simple example to highlight the issue: if one instead considers the i.i.d. sum of formal steps $X_{i}$ such that $\E[X_{i}]=0$, with distribution say 
$$\PP(X_{i}=1)=\frac{3}{2},\ \PP(X_{i}=3)=-\frac{1}{2},$$
where $\PP$ denotes a signed measure of total mass 1, then the total variation of the tail ``probability" of $X_{1}+\ldots+X_{M}$ grows to infinity as $M\rightarrow\infty$. Therefore, the weak convergence of the signed measure does not imply the convergence of moments.
\vspace{0.2cm}

\noindent{\textbf{Question:}} Can one generalize the probabilistic interpretation of Dunkl action to formal random walks governed by signed measure, or is a different point of view needed?

\noindent{\textbf{Question:}} Suppose the probabilistic argument can be generalized. Recall also that the sum in Dunkl action does not directly involve the parameter $z_{c}(N)$ (the walk starts and ends at 0, so powers of $z_{c}(N)$ cancel out after conditioning). However, $z_{c}(N)$ and $V_{N}(z_{c}(N))$ determine the scaling factor/variance of the random walk. Is this still the case for negative cumulants? If yes, among the two candidates $z_{c,1}(N)$ and $z_{c,2}(N)$, what makes the functional ``choose" $z_{c}(N)$ that maximizes $|V_{N}(z)|$?   

\subsection{Edge universality for general $\beta$-additions}
Going beyond G$\beta$E and L$\beta$E, the $\beta$-additions of general ensembles have been introduced and studied in \cite{GM}, \cite{BCG}, \cite{Xu}, \cite{CX} and \cite{Y}. While all these works focus on the global limits of the empirical measure, it is natural as well to consider its edge limit in a proper sense. Inspired by the results in \cite{Ahn} and \cite{JP},  we state a conjecture that extends the Airy edge universality from $\beta=1,2$, as well as the finite additions of Gaussian and Laguerre ensembles, to a large class of general $\beta$-additions.

Let $a(N)=(a_{1}(N)\ge \ldots\ge a_{N}(N))$, $b(N)=(b_{1}(N)\ge \ldots\ge b_{N}(N))$ be two sequences of random $N$-tuples in the Weyl chamber $W_{N}$, with exponentially decaying distributions $\mu_{a(N)}$ and $\mu_{b(N)}$. 

It is known that for general $\beta>0$, the product \[G_N(x_1,\ldots,x_N;\beta,\mu_{c(N)}):=G_N(x_1,\ldots,x_N;\beta,\mu_{a(N)})G_N(x_1,\ldots,x_N;\beta,\mu_{b(N)})\]
defines a unique tempered distribution $\mu_{c(N)}$  on $W_N$, see e.g., \cite{A} for more details. Moreover, when $\beta=1,2,4$, $\mu_{c(N)}$ is a legitimate probability measure that gives the distribution of the eigenvalues of \[diag(a(N))+U_{N}diag(b(N))U_{N}^{*},\]
where $U_{N}$ is an $N\times N$ Haar random orthogonal/unitary/symplectic matrix, respectively. We view $\mu_{c(N)}$ as the distribution of a random $N$-tuple $c(N)$ in $W_{N}$, consisting of the eigenvalues of the $\beta$-additions of two virtual Hermitian matrices $A(N)$ and $B(N)$ with eigenvalues $a(N)$, $b(N)$ respectively.

We further assume that \begin{itemize}
    \item As $N\rightarrow\infty$, the empirical measures \[\mu_{A(N)}=\frac{1}{N}\sum_{i=1}^{N}\delta_{a_{i}(N)/N}\rightarrow\mu_{A},\quad \mu_{B(N)}=\frac{1}{N}\sum_{i=1}^{N}\delta_{b_{i}(N)/N}\rightarrow\mu_{B}\]
for some compactly supported deterministic measures $\mu_{A}$, $\mu_{B}$ on $\R$.
\item Let $\mu_{+}(A)$, $\mu_{+}(B)$ be the right endpoints of $\operatorname{supp}(\mu_{A})$, $\operatorname{supp}(\mu_{B})$, and $\mu_{+}(A),\mu_{+}(B)>0$. 
\end{itemize} Then by \cite[Theorem 1.3]{CX}, the empirical measure of $c(N)/N$ converges in moments to $\mu_{A}\boxplus\mu_{B}$, where $\boxplus$ is the free convolution. We then expect that the largest entries of $c(N)/N$ will also converge to $\mathrm{Airy}(\beta)$ under proper scaling. Since it is still open whether $c(N)$ is supported by a probability measure, we state our conjecture in terms of moments using Dunkl operators.
\vspace{0.2cm}

\noindent{\textbf{Conjecture.}} Under certain mild conditions on $\mu_{a(N)},\ \mu_{b(N)},\ \mu_{A}$ and $\mu_{B}$ (see \cite[Definition 2.2 and Assumption 2.3]{JP}), there exists a positive constant $\tilde{C}$ and a sequence of positive constants $\mu_{+}(N)$, such that for any $\vec{\bk}=(\bk_{1},\ldots,\bk_{l})\in \R_{>0}^{l}$, $M_{j}=\lfloor \bk_{j}N^{{2/3}}\rfloor$ for $1\le j\le l$, $M=\sum_{j=1}^{l}M_{j}$, as $N\rightarrow\infty$, $\mu_{+}(N)\rightarrow\mu_{+},$ the right endpoint of $\mu_{A}\boxplus\mu_{B}$, and
\begin{align*}
    \frac{1}{\mu_{+}(N)^{M}}\cdot \prod_{j=1}^{l}\left[\left(\frac{\mathcal{D}_{1}}{N}\right)^{M_{j}}+\ldots+\left(\frac{\mathcal{D}_{N}}{N}\right)^{M_{j}}\right]G_{N}(\vec{x};\beta,\mu_{c(N)})\rightarrow &\bL_\beta\left(\frac{\mu_{+}}{2\tilde{C}}\vec\bk, \vec\ttt=0\right).\\
    \end{align*}

For $\beta=1,2$, this is proved in \cite{Ahn} and \cite{JP} using different approaches. We note that in particular, one can take $a(N)$, $b(N)$ to be two deterministic sequences, that satisfy the assumptions in these works. This implies that the Airy statistic arises from mixing two Hermitian matrices, whose only randomness comes from the uniformly distributed eigenvectors.  As for the proofs, \cite{Ahn} also relies on BGF and moment methods, but the operators in that work are a version of Macdonald difference operators and are limited to the case $\beta=2$. \cite{JP} establishes a Green function comparison theorem for both $\beta=1,2$, and their proof relies on the Hermitian matrix structure, which no longer exists for general $\beta$.

Besides the difficulties discussed in Section \ref{sec:negativecumulants}, for general $\beta$-additions, $G_{N}(\vec{x};\beta,\mu_{c(N)})$ no longer has the explicit expression as in Eqn.~(\ref{eq_additionbgf}). Because of these limitations, we suspect that, in order to prove the full conjecture, it will be essential to acquire a deeper understanding of the combinatorics of multivariate Bessel functions. 

\appendix

\section{Generalized ballot problem}
Suppose we have a walk $W(t)$, $t=0,1,2,\ldots,L$, whose possible increments are $\{(1,k)\}_{k=-1,0,1,\ldots}$. We assign a weight of $w_k$ to a step of vertical increment $k$, and let the weight of each walk be the product of the weights of its increments. Fix the endpoints of all the walks to be starting at $(0,y_0)$ and ending at $(L,0)$. We make the following definitions:
\begin{itemize}
    \item Define $Z_{(y_0,0,L)}(\{w\})$ to be the partition function for these walks as a function of their weights $\{w\}=\{w_{-1},w_0,w_1,w_2,\ldots\}$.
    \item Define $Z_{(y_0,0,L)}^{good}(\{w\})$ to be the partition function for the walks constrained to stay strictly above the $x$-axis before the $L^{th}$ step. We call such walks ``good''.
    \item Define $Z_{(y_0,0,L)}^{bad}(\{w\})$ to be the partition function for the walks that are not ``good''. We call such walks ``bad''.
\end{itemize}

We are interested in computing the partition function of the ``good'' paths. In the literature these types of paths are known as \L{}ukasiewicz paths. See \cite{Lukasiewicz}, and references therein, for a study of these paths with uniform weights. 

\begin{prop} \label{prop:Ballot}
The partition function for the ``good'' paths can be written as
    \begin{equation}
Z_{(y_0,0,L)}^{good}(\{w\}) = \frac{y_0}{L} Z_{(y_0,0,L)}(\{w\}).
\end{equation}
\end{prop}

\begin{remark}
When $w_{-1}=w_1=1$ and all other $w_i$ are 0, this is equivalent to Bertrand's ballot problem \cite{Ber}. Many other variations of the problem have been solved. The variation presented here in Proposition \ref{prop:Ballot} is essentially done in \cite[Section 6]{GS} using a generalized version of the ``reflection principle". For completeness, we present a proof below.
\end{remark}

\begin{proof}
The last step of a ``bad'' walk can either be a down-, flat-, or up-step, where the up-step can be of height $k=1,2,\ldots$. Grouping these walks based on their last step, we may write their partition function as
\begin{equation}\label{eq:badPaths}
Z_{(y_0,0,L)}^{bad}(\{w\}) = Z_{(y_0,0,L)}^{bad,-}(\{w\}) +Z_{(y_0,0,L)}^{0}(\{w\}) + \sum_{k=1}^\infty Z_{(y_0,0,L)}^{+,k}(\{w\})
\end{equation}
where $Z_{(y_0,0,L)}^{bad,-}(\{w\})$ is the partition function for ``bad'' walks that end with a down-step, $Z_{(y_0,0,L)}^{0}(\{w\})$ is the partition function for walks that end with a flat step, and $Z_{(y_0,0,L)}^{+,k}(\{w\})$ is the partition function for walks that end with an up-step of height $k$. Note that walks whose last step is a flat step or an up-step are necessarily ``bad''.

Each bad walk that ends with a down-step must have some last time (before step $L$) at which it hit the $x$-axis. We can partition the walks in $Z_{(y_0,0,L)}^{bad,-}(\{w\})$ by the size $k$ of the up-step and the height $i=1,\ldots,k$ of the walk after that step. Let $Z_{(y_0,0,L)}^{bad,-,k,i}(\{w\})$ be the partition function for these walks with fixed $k$ and $i$.

For each such walk, we can take all the steps after (and including) the last time it hits the $x$-axis and rotate this piece of the walk by 180 degrees. Doing this results in a walk whose last step is an up-step of size $k$. Note that this is invertible: for any walk whose last step is an up-step of size $k$, we can find the last down-step starting from height $-(k-i)$ and rotate all of the walk starting from this step. See Example \ref{ex:Ballot}. Note that these rotations are weight-preserving. It follows that
\[
Z_{(y_0,0,L)}^{bad,-,k,i}(\{w\}) = Z_{(y_0,0,L)}^{+,k}(\{w\})
\]
for each $k$ and $i$. Since for each $k$ we have that $i$ ranges over $1,\ldots, k$, summing the above over all $i$ and $k$ results in
\begin{equation}\label{eq:badPath-}
Z_{(y_0,0,L)}^{bad,-}(\{w\}) = \sum_{k=1}^{\infty} k\;  Z_{(y_0,0,L)}^{+,k}(\{w\}).
\end{equation}

Let $Z^{+}$, $Z^{0}$ and $Z^{-}$ denote the partition functions of the walks with positive/zero/negative increment in the last step, respectively. From Eqn.~\eqref{eq:badPaths} and Eqn.~\eqref{eq:badPath-}, we see that
\begin{equation}
\begin{aligned}
Z_{(y_0,0,L)}^{good}(\{w\}) =&\; Z_{(y_0,0,L)}(\{w\}) -Z_{(y_0,0,L)}^{bad}(\{w\})  \\
=&\; Z_{(y_0,0,L)}^{+}(\{w\}) + Z_{(y_0,0,L)}^{0}(\{w\})+Z_{(y_0,0,L)}^{-}(\{w\}) - Z_{(y_0,0,L)}^{bad}(\{w\})  \\
=&\; Z_{(y_0,0,L)}^{-}(\{w\}) - \sum_{k=1}^{\infty} k\; Z_{(y_0,0,L)}^{+,k}(\{w\}).
\end{aligned}
\end{equation}

Now, suppose we fix the number of steps of size $i=-1,0,1,2,\ldots$ to be $n_i$. Because we fix the endpoints of the walks to be at $(0,y_0)$ and $(L,0)$, we must have
\[
\sum_i n_i = L \qquad \text{and} \qquad \sum_i i\; n_i = -y_0.
\]
The weight of a walk with this distribution of steps is given by $\prod_i w_i^{n_i}$. The contribution that such walks make to $Z_{(y_0,0,L)}^{-}(\{w\}) - \sum_k k\; Z_{(y_0,0,L)}^{+,k}(\{w\})$ can be written as
\[
\begin{aligned}
&\left(\binom{L-1}{n_{-1}-1,n_0,n_1,\ldots} - \sum_k k\; \binom{L-1}{n_{-1},n_1,\ldots,n_k-1,\ldots} \right) \prod_i w_i^{n_i} \\
&=\frac{n_{-1} - \sum_k k\; n_k}{L} \binom{L}{n_{-1},n_0,n_1,\ldots}\prod_i w_i^{n_i}  \\
&=\frac{y_0}{L} \binom{L}{n_{-1},n_0,n_1,\ldots}\prod_i w_i^{n_i}.
\end{aligned}
\]
Finally, summing over all choices of the $n_i$, we obtain
\[
Z_{(y_0,0,L)}^{good}(\{w\}) = \frac{y_0}{L} Z_{(y_0,0,L)}(\{w\}).
\]
\end{proof}

\begin{example}\label{ex:Ballot}
Here we have $y_0=2$ and $L=9$. Below is a walk in $(\mathrm{bad},-,3,2)$, that is the last step to hit the $x$-axis is a step of size 3 that ends at height 2:
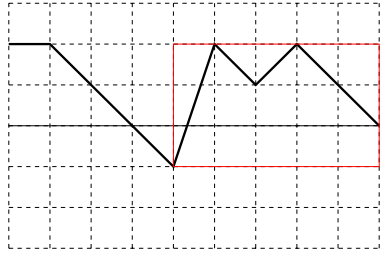
\[
\resizebox{0.3\textwidth}{!}{
\begin{tikzpicture}
    \draw[dashed] (1,-3) grid (10,3);
    \draw (1,0)--(10,0);
    \draw[ultra thick] (1,2)--(2,2)--(3,1)--(4,0)--(5,-1)--(6,2)--(7,1)--(8,2)--(9,1)--(10,0);
    \draw[red] (5,-1)--(10,-1)--(10,2)--(5,2)--(5,-1);
\end{tikzpicture}
}
\]
After the rotation, we get
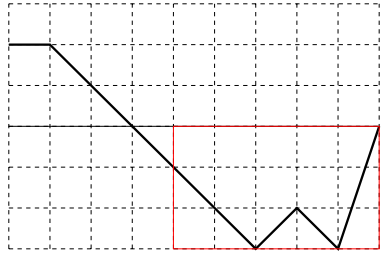
\[
\resizebox{0.3\textwidth}{!}{
\begin{tikzpicture}
    \draw[dashed] (1,-3) grid (10,3);
    \draw (1,0)--(10,0);
    \draw[ultra thick] (1,2)--(2,2)--(3,1)--(4,0)--(5,-1)--(6,-2)--(7,-3)--(8,-2)--(9,-3)--(10,0);
    \draw[red] (5,0)--(10,0)--(10,-3)--(5,-3)--(5,0);
\end{tikzpicture}
}
\]
a walk in $(+,3)$.

Note that the walks
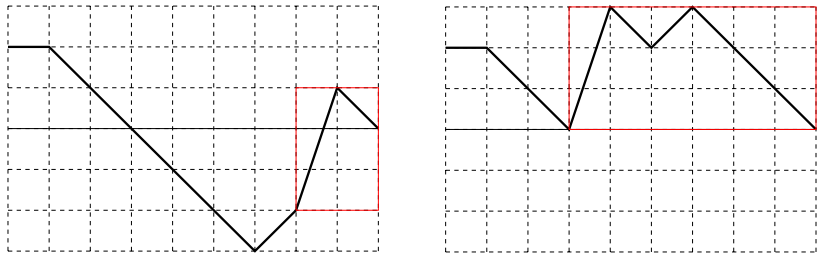
\[
\resizebox{0.3\textwidth}{!}{
\begin{tikzpicture}
    \draw[dashed] (1,-3) grid (10,3);
    \draw (1,0)--(10,0);
    \draw[ultra thick] (1,2)--(2,2)--(3,1)--(4,0)--(5,-1)--(6,-2)--(7,-3)--(8,-2)--(9,1)--(10,0);
    \draw[red] (8,-2)--(10,-2)--(10,1)--(8,1)--(8,-2);
\end{tikzpicture}
}
\qquad
\resizebox{0.3\textwidth}{!}{
\begin{tikzpicture}
    \draw[dashed] (1,-3) grid (10,3);
    \draw (1,0)--(10,0);
    \draw[ultra thick] (1,2)--(2,2)--(3,1)--(4,0)--(5,3)--(6,2)--(7,3)--(8,2)--(9,1)--(10,0);
    \draw[red] (4,0)--(10,0)--(10,3)--(4,3)--(4,0);
\end{tikzpicture}
}
\]
in $(\mathrm{bad},-,3,1)$ and $(\mathrm{bad},-,3,3)$, respectively, also map to the same walk in $(+,3)$ after rotating their highlighted blocks.
\end{example}

Let us now restrict to the types of walks defined in Section \ref{sec:walkbridges} in order to give a proof of Lemma \ref{lem:excursionpartition} which we restate below as a corollary to Proposition \ref{prop:Ballot}. Recall that we consider walks with weight given by
\[
w_{-1}(N) = \frac{1}{z_c(N)} \quad \text{and} \quad w_{l-1}(N) = \kappa_{l}(N)z_{c}(N)^{l-1}, \quad l\ge 1
\]
where the $\kappa_{l}(N)$ are the finite-$N$ version of the free cumulants of our matrix additions \eqref{eq:addition} and $z_c(N)$ is the unique positive zero of their Voiculescu transform
\[
V_{N}(z) = \frac{1}{z} + \sum_{l= 1}^\infty \kappa_l(N) z^{l-1}.
\]

\begin{cor}
    The partition function for the walk bridges  from $(0,H)$ to $(L,0)$, with weights given in Definition \ref{def:walkbridge}, is 
\begin{equation}
  \frac{H+1}{L+1} \frac{1}{2\pi i} \oint\; \left(\frac{z}{z_c(N)}\right)^H V_{N}(z)^{L+1} \; dz
\end{equation}
where the contour only encloses the pole at $0$. 
\end{cor}
\begin{proof}
We make a couple of observations. First, note that $V_N(z_c(N))^K$, $K\in \mathbb{Z}_{\ge0}$, can also be seen as the generating function for walks of length $K$, where the power of $z_c(N)$ corresponds to the total change in height of the walks. Second, rather than considering walks from $(0,H)$ to $(L,0)$ that stay weakly above $y=0$, we may consider walks from $(0,H+1)$ to $(L+1,0)$ that stay strictly above $y=0$. To do this we take each of the original walks, shift them up by one, and add an extra down-step to the end of the walk. This gives a bijection between the sets of walks with the difference in weight given by a factor of $\frac{1}{z_c(N)}.$

As these walks have a fixed change in height, the power of $z_c(N)$ is fixed and the same for each walk. It is given by $z_c(N)^{-(H+1)}$. By the residue theorem, the partition function for walks from $(0,H+1)$ to $(L+1,0)$ with no additional restrictions is given by
\[
\frac{1}{2\pi i} \oint\; \frac{z^H}{z_c(N)^{H+1}} V_{N}(z)^{L+1} \; dz.
\]
It follows from Proposition \ref{prop:Ballot} that the partition function for walks from $(0,H+1)$ to $(L+1,0)$ staying strictly above the $x$-axis is given by
\[
\frac{H+1}{L+1}\frac{1}{2\pi i} \oint\; \frac{z^H}{z_c(N)^{H+1}} V_{N}(z)^{L+1} \; dz.
\]
As mentioned above, the difference in weight between these walks and the walks we consider is a single factor of $\frac{1}{z_c(N)}$, and the result follows. 
\end{proof}

\section{ Proof of Proposition \ref{prop:analytic}}
\begin{proof}
We introduce an auxiliary approximation. For any $x\in\R$, $\alpha>0$, and $N\in\N$, define
\[
M_+'[x,\alpha,N] = \left(\max\{1+N^{-2/3}x, 0\}\right)^{2\lfloor \alpha N^{2/3}/2\rfloor}.
\]
First consider the case $x>-N^{2/3}$. We have
\[
|M_+'[x,\alpha,N]-M[x,\alpha,N]| = N^{-2/3}|x|M_+[x,\alpha,N].
\]
Fix $\delta_-, \delta_+>0$. Observe that $|x|\le \delta_+^{-1}(1+N^{-2/3}x)^{\delta_+ N^{2/3}}$ for $x>0$, and $|x|\le \delta_-^{-1}(1+N^{-2/3}x)^{-\delta_- N^{2/3}}$ for $-N^{2/3}<x<0$. Consequently, for $x>-N^{2/3}$, $\alpha>\delta_-$, and sufficiently large $N$,
\begin{equation}  \label{eq:mpmdif}
|M_+'[x,\alpha,N]-M[x,\alpha,N]| \le 2 N^{-2/3}(\delta_+^{-1}M_+[x,\alpha+\delta_+,N]+\delta_-^{-1} M_+[x,\alpha-\delta_-,N]).
\end{equation}
For the case $x < -N^{2/3}$, a symmetry argument applies. Let $x'=-2N^{2/3}-x > -N^{2/3}$. Then $|M_+' - M| = |2+xN^{-2/3}|M_+$. Since $|2+xN^{-2/3}| = |x'|N^{-2/3}$, and $M_+$ is symmetric with respect to the transformation $x \to x'$, the bound \eqref{eq:mpmdif} holds generally for all $x \in \R$ (noting that at $x=-N^{2/3}$, the difference is zero).

Take $N=N_{k}$. Since $\sum_{i=1}^{N_k} M_+[x_i^{N_k},\bk,N_k]$ is bounded for all relevant $\bk$, the error term on the right-hand side of Eqn.~\eqref{eq:mpmdif} vanishes as $k\to\infty$ due to the factor $N_k^{-2/3}$. Thus, for any $\alpha>0$, the sequence $\sum_{i=1}^{N_k} M_+'[x_i^{N_k},\alpha,N_k]$ converges if and only if $\sum_{i=1}^{N_k} M[x_i^{N_k},\alpha,N_k]$ converges, and their limits are identical.

We now prove the convergence of the first particle $x_1^{N_k}$.
Since $M_+[x_1^{N_k},1,N_k]$ is bounded uniformly in $k$, the sequence $x_1^{N_k}$ is bounded from above.
Suppose for the sake of contradiction that $x_1^{N_k}$ does not converge. Then there exist two subsequences of $\{N_k\}$, denoted by $\{P_k\}$ and $\{Q_k\}$, such that $x_1^{P_k} \to a$ and $x_1^{Q_k} \to b$ as $k\to\infty$, with $a, b \in \R\cup\{-\infty\}$ and $a<b$.

For any $\alpha\in\N$, utilizing the limit $b$ along the subsequence $Q_k$, we have
\[
\lim_{k\to\infty}\sum_{i=1}^{P_k} M_+'[x_i^{P_k},\alpha,P_k] = \lim_{k\to\infty}\sum_{i=1}^{Q_k} M_+'[x_i^{Q_k},\alpha,Q_k] \ge \lim_{k\to\infty} M_+'[x_1^{Q_k},\alpha,Q_k] = \exp(\alpha b).
\]
Conversely, considering the subsequence $P_k$ where $x_1^{P_k} \to a$, and noting that $x_1^{P_k} \ge x_i^{P_k}$ implies an exponential bound, we estimate:
\[
\lim_{k\to\infty} \sum_{i=1}^{P_k} M_+'[x_i^{P_k},\alpha,P_k] \le \exp((\alpha-1)a) \lim_{k\to\infty} \sum_{i=1}^{P_k} M_+'[x_i^{P_k},1,P_k].
\]
Combining these inequalities yields
\[
\exp(\alpha(b-a)+a) \le \lim_{k\to\infty} \sum_{i=1}^{P_k} M_+'[x_i^{P_k},1,P_k].
\]
Since $b > a$, the left-hand side grows exponentially with $\alpha$, while the right-hand side is a finite constant independent of $\alpha$. This is a contradiction. Thus, $x_1^{N_k}$ must converge.

Given the convergence of $x_1^{N_k}$, the first term in the sum converges. Subtracting this term implies that $\sum_{i=2}^{N_k} M[x_i^{N_k},\alpha,N_k]$ converges as $k\to\infty$ for any $\alpha\in \N$. By induction, we conclude that $x_i^{N_k}$ converges to some $x_i$ for every $i\in\N$.

Finally, for any $\alpha>0$, we prove the convergence of the sum to the exponential limit
\begin{equation} \label{eq:alphaconv}
    \sum_{i=1}^{N_k} M_+'[x_i^{N_k},\alpha,N_k]\to \sum_{i=1}^\infty \exp(\alpha x_i).
\end{equation}
Let $X=\limsup_{k\to\infty}\sum_{i=1}^{N_k} M_+'[x_i^{N_k},\alpha,N_k]$. Choose $\alpha_1, \alpha_2$ such that $\alpha_1^{-1}, \alpha_2 \in \N$ and $\alpha_1 < \alpha < \alpha_2$. The convexity inequality $y^\alpha \le y^{\alpha_1} + y^{\alpha_2}$ for $y\ge 0$ combined with the hypotheses ensures that $X < \infty$.
For any fixed $j$, since the terms are ordered ($x_i \ge x_{i+1}$), we see that for any $j' \ge j$ we have
\[
 M_+'[x_{j'}^{N_k},\alpha,N_k] \le \frac{1}{j}\sum_{i=1}^{N_k} M_+'[x_i^{N_k},\alpha,N_k].
\]
This tail decay implies
\[
\limsup_{k\to\infty}\sum_{i=j}^{N_k} M_+'[x_i^{N_k},\alpha,N_k]
\le (j^{-1}X)^{1-\alpha_1/\alpha} \lim_{k\to\infty} \sum_{i=1}^{N_k} M_+'[x_i^{N_k},\alpha_1,N_k].
\]
Consequently, the contribution of the tail tends to zero as $j\to\infty$:
\[
\lim_{j\to\infty}\limsup_{k\to\infty}\sum_{i=j}^{N_k} M_+'[x_i^{N_k},\alpha,N_k] = 0.
\]
Combining this tail property with the pointwise convergence $x_i^{N_k} \to x_i$, we obtain Eqn.~\eqref{eq:alphaconv}.
\end{proof}

\end{document}